%% file: main.tex
\documentclass[11pt]{amsart}

\input{amsmacs}

\input{psfig}

\begin{document}

\title[Geometry of the complex of curves II]{Geometry of the complex
                of curves II: hierarchical structure}
\author{Howard A. Masur}
\author{Yair N. Minsky}
\address{University of Illinois at Chicago}
\address{SUNY Stony Brook}
\date{July 14, 1998}
\thanks{
The first author was partially supported by NSF grant \#DMS
9201321. The second author was partially supported by a fellowship
from the Alfred P. Sloan Foundation and NSF grant \#DMS 9626233.}

\maketitle

\newcommand{\dist}{\operatorname{dist}}
\newcommand\bbar{\overline}
\newcommand{\I}{{\mathbf I}}
\newcommand{\T}{{\mathbf T}}
\newcommand{\tprec}{\prec_t}
\newcommand{\fprec}{\prec_f}
\newcommand{\bprec}{\prec_b}
\newcommand{\pprec}{\prec_p}
\newcommand{\sprec}{\prec_s}
\newcommand{\fsub}{\mathrel{\scriptstyle\searrow}}
\newcommand{\bsub}{\mathrel{\scriptstyle\swarrow}}
\newcommand{\fsubd}{\mathrel{{\scriptstyle\searrow}\kern-1ex^d\kern0.5ex}}
\newcommand{\bsubd}{\mathrel{{\scriptstyle\swarrow}\kern-1.6ex^d\kern0.8ex}}
\newcommand{\fsubeq}{\mathrel{\genfrac{}{}{0pt}{}{\searrow}{\raise.1ex\hbox{=}}}}
\newcommand{\bsubeq}{\mathrel{\genfrac{}{}{0pt}{}{\swarrow}{\raise.1ex\hbox{=}}}}
\newcommand{\tw}{\operatorname{tw}}
\newcommand{\base}{\operatorname{base}}
\newcommand{\trans}{\operatorname{trans}}
\newcommand{\rest}{|_}

\long\def\sotto#1{{\bf [[} {\em #1 } {\bf ]]}}

\renewcommand\marginpar[1]{}  

\setcounter{tocdepth}{1}
\tableofcontents
\input{intro}
\input{defs}

\input{projection}

\input{paths}

\input{slices}

\input{control}

\input{conjugacy}

\input{biblio}

\end{document}

%% file: amsmacs.tex
\catcode`\@=11

\long\def\@savemarbox#1#2{\global\setbox#1\vtop{\hsize\marginparwidth 
  \@parboxrestore\tiny\raggedright #2}}
\long\def\sotto#1{{\bf [[} {\em #1 } {\bf ]]}}

\newcommand\lref[1]{\ref{#1}%
\@ifundefined{r@DisplaY #1}{}{ (#1)}}

\newcommand\fakelabel[2]{\@bsphack\if@filesw {\let\thepage\relax
   \newcommand\protect{\noexpand\noexpand\noexpand}%
\xdef\@gtempa{\write\@auxout{\string
      \newlabel{#1}{{#2}{\thepage}}}}}\@gtempa
   \if@nobreak \ifvmode\nobreak\fi\fi\fi\@esphack}

\catcode`\@=12

\def\Empty{}
\newcommand\oplabel[1]{
  \def\OpArg{#1} \ifx \OpArg\Empty {} \else
        \label{#1}
  \fi}
\newtheorem{theoremSt}{Theorem}[section]

\newtheorem{exampleSt}[theoremSt]{Example}
\newtheorem{exerciseSt}[theoremSt]{Exercise}

\newcommand\MakeStEnv[1]{
  \newenvironment{#1}[1]{
  \begin{#1St} \oplabel{##1}%
  \global\def\CrntSt{\thetheoremSt}%
}{ 
  \end{#1St} }
  \newenvironment{#1+}[1]{
  \begin{#1St} \label{##1}%
  \label{DisplaY ##1}%
  \global\def\CrntSt{\thetheoremSt}%
  \def\Labl{##1}\ifx\Labl\Empty{} \else {\em (\Labl)\,}\fi%
}{ 
  \end{#1St} }
}
\MakeStEnv{theorem}
\MakeStEnv{corollary}
\MakeStEnv{proposition}
\MakeStEnv{lemma}
\MakeStEnv{definition}
\MakeStEnv{conjecture}

\long\def\Restate#1#2#3{
\medskip\par\noindent
{\bf #1 \ref{#2}} {\rm(#2)}
{\it #3}
\par\medskip
}

\long\def\realfig#1#2#3{
\begin{figure}[htbp]
\centerline{\psfig{figure=#2}}
\caption[#1]{#3}
\oplabel{#1}
\end{figure}}

\long\def\scalefig#1#2#3#4{
\begin{figure}[htbp]
\centerline{\psfig{figure=#3,height=#2}}
\caption[#1]{#4}
\oplabel{#1}
\end{figure}}
                                        
\newlength{\saveu}

\newenvironment{pf}{%
 \begin{proof}%
}{ 
 \end{proof}
}
\newenvironment{pf*}[1]{%
 \begin{proof}[#1]%
}{ 
 \end{proof}
}

\newcommand{\finishproof}[1]{ 
  \def\FPArg{#1}
  \ifx\FPArg\Empty
        \newcommand\FPArg{\CrntSt}  \fi
  \smallbreak\noindent\makebox[\textwidth]{\hfill\fbox{\FPArg}}
  \medbreak\noindent
}

\newcommand{\bfheading}[1]{\par\smallskip\noindent {\bf #1}}

\newcommand\CC{{\mathcal C}}

\newcommand\FF{{\mathcal F}}
\newcommand\GG{{\mathcal G}}

\newcommand\LL{{\mathcal L}}
\newcommand\MM{{\mathcal M}}
\newcommand\NN{{\mathcal N}}

\newcommand\PP{{\mathcal P}}

\newcommand\TT{{\mathcal T}}
\newcommand\UU{{\mathcal U}}

\newcommand\PMF{{\PP\kern-2pt\MM\FF}}

\newcommand\PML{{\PP\kern-2pt\MM\LL}}

\newcommand\Mod{\operatorname{Mod}}

\newcommand\ep{\epsilon}

\newcommand\hhat{\widehat}

\newcommand\union{\cup}
\newcommand\intersect{\cap}
\newcommand\bbR{{\mathord{\text{I\kern-2pt R}}}}        
\newcommand\bbH{{\mathord{\text{I\kern-2pt H}}}}        

\newcommand\C{{\mathbf C}}

\newcommand\Z{{\mathbf Z}}
\newcommand\R{{\mathbf R}}

\newcommand\Hyp{{\mathbf H}}

\newcommand\bigrightarrow[1]{\hbox to #1{\rightarrowfill}}
\newcommand\bigleftarrow[1]{\hbox to #1{\leftarrowfill}}

\newcommand\boundary{\partial}
\newcommand\semidir{\mathrel{\hbox{\vrule depth-.03ex height1.1ex\kern-0.15em$\times$}}}

\newcommand\til{\widetilde}

\newcommand{\diam}{\operatorname{diam}}

\numberwithin{equation}{section}

%% file: intro.tex
\section{Introduction}
\label{intro}

In this paper we continue our geometric study
of Harvey's Complex of Curves \cite{harvey:boundary}, a
finite dimensional and locally infinite complex $\CC(S)$ associated to
a surface $S$,
which admits an action by the mapping class group $\Mod(S)$. 
The geometry and combinatorics of $\CC(S)$ can be applied to study
group-theoretic properties of $\Mod(S)$, and the geometry of Kleinian
representations of $\pi_1(S)$. 

In \cite{masur-minsky:complex1} we showed that, endowed with a natural
metric, $\CC(S)$ is an infinite diameter $\delta$-hyperbolic space in
all but a small number of trivial cases (see Section
\ref{defs} for precise definitions). This result suggests that one try
to apply the techniques of hyperbolic spaces and groups
to study $\CC(S)$ and its $\Mod(S)$-action, considering for
example such questions as the word problem, conjugacy problem and
quasi-isometric rigidity. The barrier to doing this is that
the complex is locally infinite, and hence the distance bounds one
obtains in a typical geometric argument give little a-priori
information.

Our goal in this paper is to develop tools for lifting this barrier.
The organizing philosophy is roughly this: Links of vertices in
$\CC(S)$ are themselves complexes associated to subsurfaces. The
geometry of these links is tied to the geometry of $\CC(S)$ by a family of
{\em subsurface projection maps}, which are analogous to closest-point
projections to horoballs in classical hyperbolic space. 
This gives a layered structure to the complex, with
hyperbolicity at each level, and the main construction of our paper is
a combinatorial device used to tie these levels together, which we
call a  {\em hierarchy of tight geodesics}.

Using these constructions, we derive a number of properties of
$\CC(S)$ which are similar to those of locally finite complexes, such
as a finiteness result for geodesics with given endpoints (Theorem
\ref{Finite Geodesics}), and a convergence criterion for sequences of
geodesics (Theorem \ref{Convergence of Hierarchies}). We then apply
these ideas to study the 
conjugacy problem in $\Mod(S)$, deriving a
linear bound  on the shortest word conjugating two pseudo-Anosov
mapping classes (Theorem \ref{Conjugacy Bound}).
Along the way we describe a class of quasi-geodesic words in $\Mod(S)$
(Theorem \ref{Quasigeodesic Words}), whose lengths can be estimated using
the subsurface projection maps (Theorem \ref{Move distance
  and projections}).

The rest of Section \ref{intro} gives a  more detailed outline of 
our results, and works through some explicit examples that
motivate our constructions.
Section \ref{defs} presents our definitions
and notation, and proves some basic lemmas. Section \ref{projection}
proves our fundamental result on subsurface projections, Sections
\ref{hierarchies} and \ref{resolution} develop the machinery of
hierarchies and their resolutions into sequences of markings, Section
\ref{large link etc} proves our basic geometric control theorems, and
Section \ref{MCG} proves the conjugacy bound theorem for $\Mod(S)$.

\subsection{Subsurface Projections}
A basic analogy for thinking about $\CC(S)$ is provided by the
geometry of a family $\FF$ of disjoint, uniformly spaced  horoballs in
$\Hyp^n$, for example 
the uniform cusp horoballs of a Kleinian group. The non-proper metric
space $X_\FF$ 
obtained by collapsing each horoball to a point is itself
$\delta$-hyperbolic -- see Farb \cite{farb:relhyp} and Klarreich
\cite{klarreich:thesis} -- and the horoballs play a role similar to
links of vertices in $\CC(S)$.

If $B$
is a horoball and $L$ is a hyperbolic geodesic disjoint from $B$, then
the closest-point projection of $L$ to $B$ has uniformly bounded
diameter, independently of $L$ or $B$. Interestingly, one can sensibly
define a ``projection'' from the collapsed space $X_\FF$ to $B$ which
similarly sends $X_\FF$-geodesics avoiding $B$ to bounded sets. This
turns out to 
be a crucial property in understanding the geometry of $X_\FF$ and
its relation to the geometry of $\Hyp^n$.

In our context, vertices of $\CC(S)$ are simple closed curves in $S$ (see
\S\ref{basic defs}) and the link of a vertex $v$ is closely related to
the complexes  $\CC(Y)$ for the complementary subsurfaces $Y$ of  $v$.
We will define projections $\pi_Y$ from $\CC(S)$ to $\CC(Y)$ as follows:
given a simple closed curve on $S$ take its arcs of
intersection with $Y$ and perform a surgery on them to obtain closed
curves in $Y$. (More precisely $\pi_Y$ sends vertices in $\CC(S)$ to
finite sets in $\CC(Y)$).
We will prove the following analogue to the situation
with horoballs:

\Restate{Theorem}{Bounded Geodesic Image}{
If $Y$ is an essential subsurface of $Z$ and $g$ is a geodesic in
$\CC(Z)$ all of whose vertices intersect $Y$ nontrivially, then the
projected image of $g$ in $\CC(Y)$ has uniformly bounded diameter.
}

The family $\FF$ of horoballs also satisfies the closely related
``bounded coset penetration property'' of Farb \cite{farb:relhyp},
which roughly speaking 
is a stability property for paths in $\Hyp^n$ whose images in $X_\FF$ are
quasi-geodesics: if two such paths begin and end near each other, then
up to bounded error they penetrate through the same set of horoballs
in the same way. This property does not hold in our case but a certain
generalization of it does. This will be the content of Lemmas
\lref{Large Link} and \lref{Common Links}, which will be briefly
discussed in \S\ref{hierarchy summary} below.

\subsection{The conjugacy problem}
Fix a set of generators for $\Mod(S)$ and let $|\cdot|$ denote the
word metric. As one application of our techniques,
in Section \ref{MCG} we establish the following
bound: 

\Restate{Theorem}{Conjugacy Bound}{Fix a surface $S$ of finite type and a generating set
  for $\Mod(S)$. If $h_1, h_2$ are
words describing 
conjugate pseudo-Anosov elements, then the shortest conjugating
element $w$ has word length 
$$|w| \le C(|h_1|+|h_2|),$$
where the constant $C$ depends only on $S$ and the generating set.
}

This linear growth property for the
shortest conjugating word is shared with word-hyperbolic groups
(see Lys\"enok \cite[Lemma 10]{lysenok:hyperbolic}),
although except in a few low-genus cases
the mapping class  group
is not word hyperbolic since it contains abelian subgroups of
rank at least $2$ generated by Dehn twists about disjoint curves.
Our proof is based on a proof which works in the word-hyperbolic case.
The case of general elements of $\Mod(S)$ introduces complications
similar to those that occur for torsion elements of word-hyperbolic
groups. We hope to address the general case in a future paper. 

This bound is related to the question of solubility of the conjugacy
problem, since a computable bound on $w$ provides a bounded search
space for an algorithm seeking to establish or refute conjugacy.
Hemion \cite{hemion:conjugacy} proved that 
the conjugacy problem  for $\Mod(S)$ is soluble, and Mosher
\cite{mosher:conjugacy} gave an explicit algorithm for 
determining conjugacy
for pseudo-Anosovs. In both cases, no explicit bound on the complexity
was given (although Mosher's algorithm is fast in practice).
Theorem \ref{Conjugacy Bound} is still short of a good complexity
bound since we have not described an efficient way to search through
the possible conjugating words. However, we are hopeful that the
techniques of this paper can be extended to give a more complete
algorithmic approach.

\subsection{Finiteness results}
In a locally finite graph, there are finitely many geodesics
between any two points. In
$\CC(S)$ this is easily seen to be false even for geodesics of length
2. However we shall introduce a finer notion of {\em tight geodesic}
(\S\ref{hierarchy defs}), for which we can establish:

\Restate{Theorem}{Finite Geodesics}{Between any two vertices in $\CC(S)$
  there are finitely many tight geodesics.}

This is part of a collection of results showing that in several useful
ways $\CC(S)$ is like a locally finite complex. Another is Theorem
\lref{Convergence of Hierarchies}, which generalizes the property of a
locally finite complex that any sequence of
geodesics meeting a compact set has a convergent
subsequence. 
In the locally finite setting this involves a simple diagonalization
argument, and this is replaced here by an application of Theorem \ref{Bounded
  Geodesic Image} and the hierarchy construction.

The following is an application of this result:

\Restate{Proposition}{Axis}{Any pseudo-Anosov element $h\in\Mod(S)$ has a
  quasi-invariant axis in $\CC(S)$: that is, a bi-infinite geodesic
  $\beta$ such that $h^n(\beta)$ and $\beta$ are $2\delta$-fellow
  travelers for all $n\in\Z$.}

That a {\em quasi-geodesic} exists which fellow-travels its $h$-translates
is a consequence of work in \cite{masur-minsky:complex1}. The geodesic with
this property is
obtained by a limiting process using Theorem \ref{Convergence of Hierarchies}.

\subsection{Hierarchies of geodesics}
\label{hierarchy summary}
A geodesic in $\CC(S)$ is a sequence of curves in $S$, but words in $\Mod(S)$
are more closely related to sequences of pants decompositions separated
by elementary moves (replacement of one curve at a time). The
hierarchy construction is based on the idea that a geodesic can be
``thickened'' in a natural way to give a family of pants
decompositions. We will illustrate this in one of the simplest
examples, that of the five-holed sphere $S_{0,5}$, below. We will then
give a more general discussion of the construction and state some of
our main results about it. Finally in \S\ref{example} we will give a
more extended, but still relatively simple, collection of examples.

\medskip

A pants decomposition $P$ in $S=S_{0,5}$ 
is a pair of disjoint curves $\alpha,\beta$; i.e. an edge of $\CC(S)$.
An {\em elementary move} of pants $P\to P'$ fixes one of the
curves, say $\alpha$, and replaces $\beta$ with a curve $\beta'$
which intersects $\beta$ minimally and is disjoint from $\alpha$.

Now given $P = \{\alpha,\beta\}$, take some $\psi\in\Mod(S)$,  and
consider ways of connecting $P$ to 
$\psi(P)=\{\alpha',\beta'\}$.  Choose a geodesic in $\CC(S)$
whose vertices are
$\alpha=\alpha_0,\alpha_1,\ldots,\alpha_N=\alpha'$.
The subsurface $S\setminus \alpha_0$ has two
components, a three-holed sphere, and a four-holed sphere which
contains $\alpha_1$ and $\beta$. Let
$S_{\alpha_0}$ denote the four-holed sphere.
The complex of curves $\CC(S_{\alpha_0})$ is isomorphic to the Farey
graph (see \S\ref{example}),
so let us join $\beta$ to
$\alpha_1$ by a geodesic
$\beta=\gamma_0,\gamma_1,\ldots,\gamma_m=\alpha_1$ in
$\CC(S_{\alpha_0})$.  The transition  $(\alpha_0,\gamma_i)$ to
$(\alpha_0,\gamma_{i+1})$ is an elementary move in pants.  This
path concludes with the pants decomposition 
$\{\alpha_0,\alpha_1\}$ (see Figure \ref{fivehole pants}).
Now working in $S_{\alpha_1}$
join $\alpha_0$ to $\alpha_2$ by a geodesic, giving a path of
elementary moves ending with $\{\alpha_1,\alpha_2\}$.  We repeat
this procedure, eventually ending with $\psi(P_0)$.  

\realfig{fivehole pants}{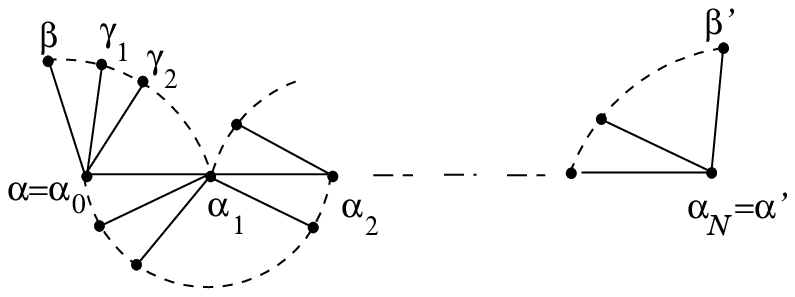}{A sequence of pants decompositions in the
  five-holed sphere. Each edge represents a pants
  decomposition. Transitions between edges are elementary moves.}

In this example, the same final pants decomposition $\psi(P_0)$ could have
been written as $\psi(\theta(P_0))$ where $\theta$ is any  product of
Dehn twists around $\alpha$ and $\beta$. 
Thus in order to keep track of elements in the mapping class group,
and not just pants decompositions,  we will also need to keep track of
twisting information around pants curves. In fact twisting data is
implicit everywhere in this example: in
the geodesic in $\CC( S_{\alpha_0})$, for each $i\in(0,m)$, 
$\gamma_{i-1}$ and
$\gamma_{i+1}$ both intersect $\gamma_i$ minimally, and hence differ by
a product of Dehn twists (or half-twists) about $\gamma_i$. Keeping
track of this information will
require the introduction of complexes associated to annular
subsurfaces (see \S\ref{annulus defs}).

Ultimately we will be considering sequences of {\em complete markings}, which
are pants decompositions together with twisting data,
such that successive markings are separated by
appropriately defined elementary moves. 
(More generally the markings need not be complete, but let us assume for
the rest of this discussion that they are.)
In considering such a sequence
carefully, one finds in it segments where some sub-marking is fixed
and all the elementary moves take place in a subsurface of $S$. Thus
one obtains some interlocking structure of paths in subcomplexes of
$\CC(S)$. 
  
{\em Hierarchies of geodesics} will be our method for constructing and
manipulating such structures. Roughly, a hierarchy is
a collection $H$ of
geodesics, each geodesic $h$ contained in a complex $\CC(Y)$ where $Y$ is
the {\em domain} of $h$. The geodesics will satisfy a technical
condition called {\em tightness}, which makes them easier to control.
There will be one ``main geodesic'' whose
domain is all of $S$, and in general the geodesics will interlock
via a relation called ``subordinacy'', which is related to the nesting
of their domains. There will be an {\em initial} and a {\em terminal}
marking, called $\I(H)$ and $\T(H)$, and a partial order on the
geodesics which is related to the linear order of a sequence of
elementary moves connecting $\I(H)$ to $\T(H)$ (the reason for a
partial rather than linear order is that some elementary moves commute
because they take place in disjoint subsurfaces).

Any hierarchy will admit a {\em resolution} of its partial order to a
linearly ordered sequence of markings, separated by elementary moves,
connecting $\I(H)$ to $\T(H)$. This resolution will be nonunique, but
efficient in the following sense.
Let $\til \MM$ be the graph whose
vertices are complete markings and whose edges are elementary
moves. We then have: 

\Restate{Theorem}{Efficiency of Hierarchies}{
Any resolution of a hierarchy $H$ into a sequence of
complete markings is a quasi-geodesic in $\til \MM$, with uniform constants.}

In the case where $\T(H) = \psi(\I(H))$ for some $\psi\in\Mod(S)$,
a resolution gives rise to a quasi-geodesic word in $\Mod(S)$ 
(Theorem \ref{Quasigeodesic Words}).

Hierarchies will be constructed inductively, with the main geodesic
chosen first and then further geodesics in subsurfaces determined by
vertices of the previous ones. At every stage a geodesic is not
uniquely determined, although hyperbolicity implies that all choices
are fellow travelers. This is a priori a fairly loose constraint, but
it has the following rigidity property, which is our generalization of
Farb's bounded coset penetration property.

\Restate{Lemma}{Common Links}{Suppose $H$ and $H'$ are hierarchies
  whose initial and 
  terminal markings differ by at most $K$ elementary moves. Then there
  is a number $M(K)$ such that, if a geodesic $h$ appears in $H$ and
  has length greater than $M$, then $H'$ contains a  geodesic $h'$
  with the same domain. Furthermore, $h$ and $h'$ are
  fellow-travelers with a uniform separation constant. }

This lemma is a consequence of the following lemma, which
characterizes in terms of the subsurface projections $\pi_Y$ when
long geodesics appear in a hierarchy. For two markings $\mu$ and $\mu'$
we let $d_Y(\mu,\mu')$ denote the distance in $\CC(Y)$ of their
projections by $\pi_Y$ (see \S\ref{markings}).

\Restate{Lemma}{Large Link}{
There exists $M=M(S)$ such that, 
if $H$ is any hierarchy in $S$ and   $d_Y(\I(H),\T(H)) \ge M$ for a
subsurface $Y$ in $S$, then $Y$ is the domain of a geodesic $h$ in
$H$.

Furthermore if $h$ is in $H$ with domain $Y$ then its length $|h|$ and
the projection distance $d_Y(\I(H),\T(H))$ are within a uniform additive
constant of each other.
}

Both of these results follow from Theorem \lref{Bounded Geodesic
  Image} together with the structural properties of hierarchies, which
  are summarized by Theorem \ref{Structure of Sigma}.

Applications of Lemmas \ref{Large Link} and \ref{Common Links}
are based on the idea that, whenever
geodesics in a hierarchy have short length, one can apply arguments
that work for locally finite complexes. Whenever geodesics become
long, one has this rigidity for all ``nearby'' hierarchies, and can
work inductively in the shared domains of the long geodesics.
Theorems \lref{Finite Geodesics}, \lref{Convergence of Hierarchies} and
\lref{Axis} are all consequences of this sort of argument.

\subsection{Motivating Examples}
\label{example}
To illustrate the above theorems, we will work through some 
more extended low-genus examples. 

Let us first take a closer look at the 
case where $S$ is a once-punctured torus or
four-times punctured sphere. Then $\CC(S)$ is the Farey graph
(See figure \ref{farey graph} and \S\ref{basic defs}), and in spite of
the fact that the link of every vertex is 
infinite we have fairly explicit and rigid control of geodesics. In
particular we note the following phenomenon. Let $h$ be a geodesic and
$v$ a vertex in $h$, preceded by $u$ and followed by $w$. The link of
the vertex $v$ can be identified with $\Z$, and we can measure the
distance between $u$ and $w$ in this link, an integer $d_v(u,w)$.
If $h'$ is a geodesic with the same endpoints as $h$, 
then $h'$ must pass through $v$ {\em provided $d_v(u,w)$ is suffiently
  large} (5 will do). In fact the
same holds if $h'$ has endpoints, say, distance 1 from those of $h$.
Furthermore, $h'$ must enter the link of
$v$ at a point within 1 of  $u$ and exit within 1 of $w$. All of these
claims are easy to show starting from the basic fact that any edge in
the Farey graph separates it.

\scalefig{farey graph}{3in}{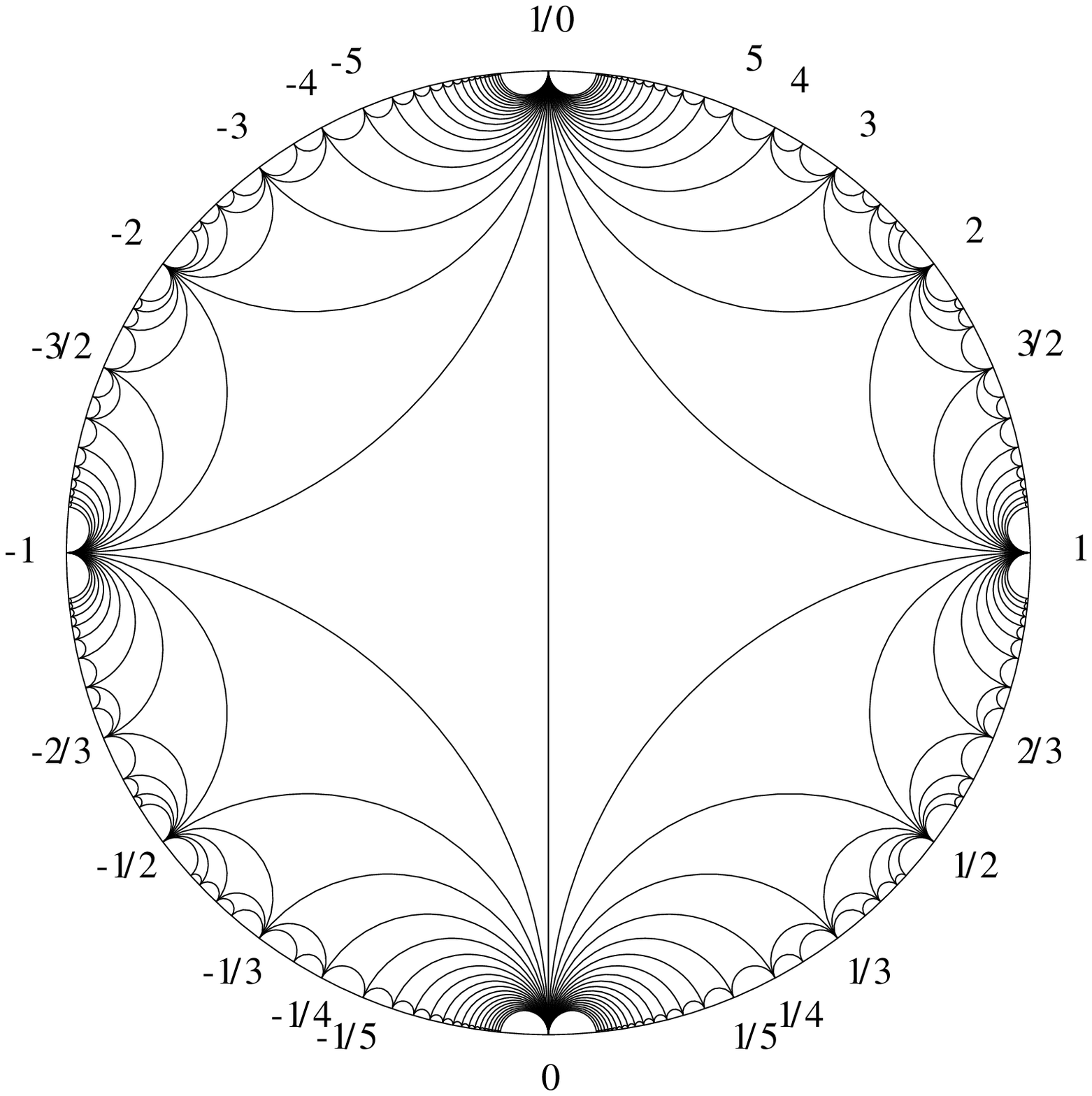}{The complex of curves for a
  once-punctured torus or
  4-times punctured sphere is the classical Farey graph. Vertices are
  labelled by slopes of the corresponding curves relative to some
  fixed homology basis.}

This phenomenon, that a large link distance generates strong constraints on
fellow traveling geodesics, persists in higher genus (even though the
separation property of edges does not generalize), and gives rise
to Lemmas \ref{Large Link} and \ref{Common Links}.
Let us now demonstrate this generalized phenomenon, together with the
main features of our hierarchy construction, in the case where $S$ is
a closed genus 2 surface.

\medskip

Let $h$ be a geodesic
in $\CC(S)$ with a segment $..,u,v,w,...$ occuring somewhere in $h$.
Let $h'$ be a fellow traveler of $h$ -- for concreteness suppose the
endpoints of $h$ and $h'$ are distance 1 or less apart, and occur at a
distance at least $2\delta+2$ from $u,v$ and $w$ (where $\delta$ is the
hyperbolicity constant of $\CC(S)$). Hyperbolicity of 
$\CC(S)$ implies that $h$ and $h'$ are
$2\delta+1$-fellow travelers. 

\realfig{genus 2 nonsep}{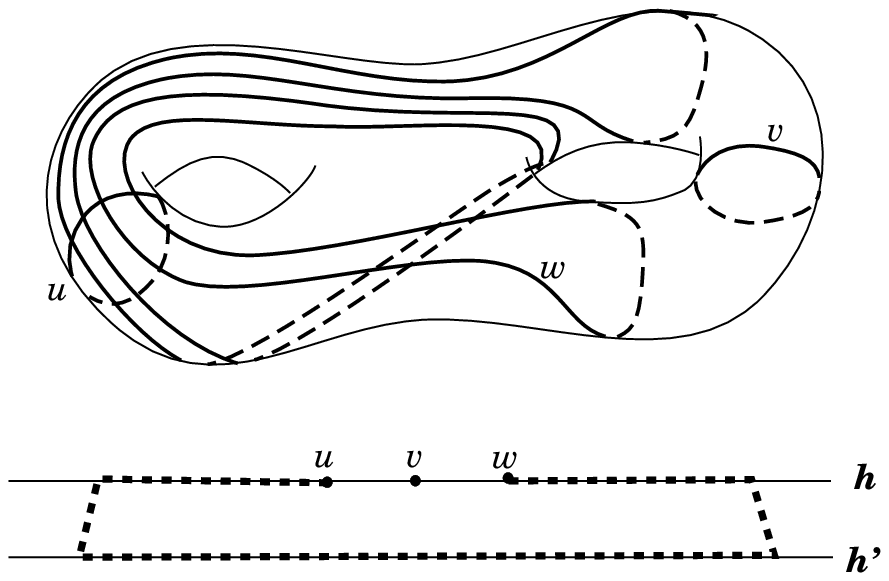}{The short-cut argument for a
  genus 2 surface, where $v$ is non-separating. If as shown $h'$ does
  not meet $v$ then the dotted rectangle bounds $d_Y(u,w)$.}

Suppose first that the subsurface $Y=S\setminus v$ is connected -- a
2-holed torus (figure \ref{genus 2 nonsep}).
Then $u$ and $w$ give points in $\CC(Y)$ and let us
denote their distance in $\CC(Y)$ by $d_Y(u,w)$. 
We can show the following statement:

{\em 
If 
the ``link distance'' $d_Y(u,w)$ is sufficiently large then
the fellow-traveler $h'$ must also pass through $v$.}

Suppose not -- then every vertex of $h'$ has
nontrivial intersection with $Y$.
Consider a path 
beginning at $w$, moving forward in $h$ a distance $2\delta+2$, across
to $h'$ by a path of length at most $2\delta+1$, back along $h'$ and over
to $h$ by another path of length at most $2\delta+1$, which lands at a
point $2\delta+2$ behind $u$, and from there back up to $u$ along $h$. 
By the triangle inequality, every point of this path not on $h'$ has
distance at least 2 from $v$ (except the endpoints $u$ and $w$ which
are in $Y$). Together with the assumption about $h'$
we have that every point on the path represents a
curve having nontrivial intersection with $Y$. The length of the
segment on $h'$ is bounded by $8\delta+8$ by the triangle inequality, so
the total length of
the path from $w$ to $u$ is at most $16\delta + 14$.
If we replace every curve with one arc of its intersection with $Y$,
we obtain a sequence of properly embedded arcs or curves in $Y$, each
disjoint from the previous. As we will see in Lemma
\ref{arcs to curves}, these
can each be replaced with a simple closed curve, so that each one is
distance at most 2 from its predecessor in $\CC(Y)$. (This is the
subsurface projection $\pi_Y$.)
We therefore obtain a
path in $\CC(Y)$ connecting $w$ to $u$, of length at most
$32\delta+28$. If we assumed
$d_Y(u,w)>32\delta+28$ this would be a contradiction, and then
$h'$ would have to pass through the vertex $v$.

In that case, we can say more.  Let $u'$ be the predecessor and $w'$
the successor  
of $v$ along $h'$. The same kind of argument, applied to the segments of
our path joining $u$ to $u'$ and $w$ to $w'$, 
gives an upper  bound 
for $d_Y(u,u')$ and $d_Y(w,w')$.
Joining $u$ and $w$ by a geodesic $k$ in $\CC(Y)$ and
$u'$ and $w'$ by a geodesic $k'$ in $\CC(Y)$, 
we now know by
hyperbolicity of $\CC(Y)$ that $k$ and $k'$ are fellow-travelers. 

Now suppose instead that $v$ divides $S$ into components $Y_1$
and $Y_2$, each necessarily a one-holed torus (see figure 
\ref{genus 2 sep}).
Since $u,v,w$ is a
geodesic, $u$ and $w$ must intersect nontrivially and hence belong to
the same component, say 
$Y_1$. The previous ``short-cut'' argument now implies that,  if
$d_{Y_1}(u,w)>32\delta+28$,  some 
curve $v'$ of $h'$ must miss $Y_1$.  We could again have $v'=v$,
or now the additional possibility that $v'$ lies (nonperipherally) in $Y_2$. 
Suppose this case  happens. Set $Y'=S\setminus v'$, noting that
it must be a single two-holed torus containing $Y_1$,
and let $u'$ and $w'$ be
the  predecessor and successor of $v'$ in $h'$.
We again apply the short-cut  argument to conclude that
$d_{Y'}(u',w')\leq 32\delta+28$; for otherwise $v'$ would
appear in $h$, but it is not $u,v$ or $w$ and is distance 1 from $v$,
so this contradicts the fact that $h$ is a geodesic. Let 
$m'$ be a geodesic in $\CC(Y')$ joining $u'$ and $w'$.
If every vertex of $m'$ intersects $Y_2$, then using $m'$, and
thus bypassing $v'$, we can find a path of some bounded length joining
$u$ and 
$w$, such that every point on it represents a curve that meets $Y_1$.
Thus assuming 
$d_{Y_1}(u,w)$ is sufficiently large, $m'$ must pass through a curve
missing $Y_1$. Since it is an essential curve in $Y'$, this curve in
fact can only be $v$ itself.
Let $y'$ and
$z'$ be the predecessor and successor of $v$ along $m'$.  They must
lie in $Y_1$ and now in fact the same argument gives  an upper 
bound for the distance in $\CC(Y_1)$ between $y'$ and $u$ and
between $z'$ and $w$.  Again by hyperbolicity any geodesic $k'$
in $\CC(Y_1)$ joining $y'$ and $z'$ fellow travels the geodesic
$k$ joining $u$ and $w$. This is 
essentially the content of Theorem \lref{Common Links} in this case.

\realfig{genus 2 sep}{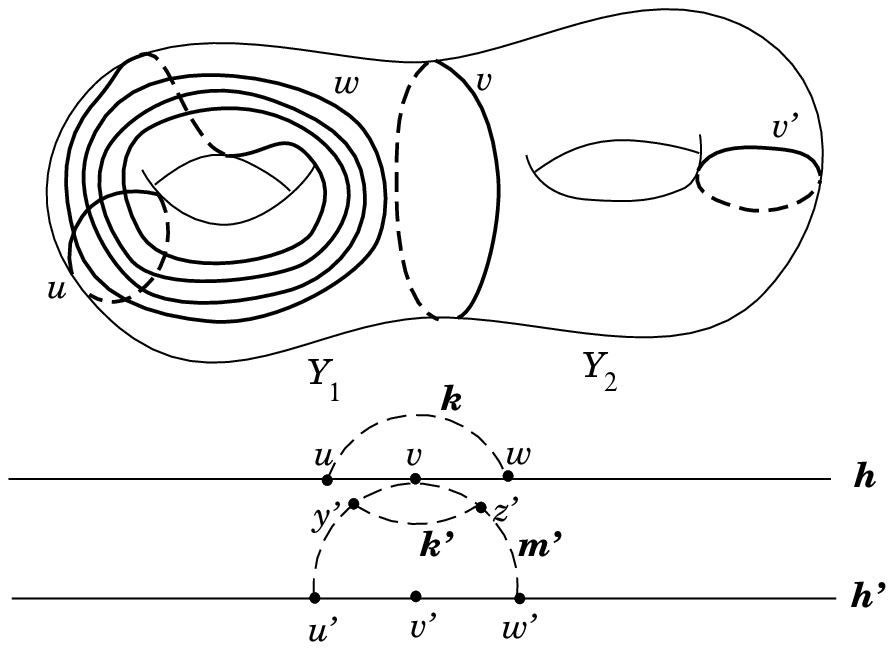}{When $v$ separates $S$ into $Y_1$
  and $Y_2$, $h'$ can pass through $v'$ in $Y_2$, But if $d_{Y_1}(u,w)$ is
  large then $m'$, supported in $Y'=S\setminus v'$, must pass through $v$.}

So far, we have constructed over $h'$ a ``hierarchy'' of geodesics:
$m'$ is obtained as a geodesic in the link of $v'$, joining 
its predecessor and its successor in $h'$. $k'$ is obtained in the
link of $v$, appearing in $m'$, in the same way. We say that $m'$ is
{\em subordinate} to $h'$, and $k'$ to $m'$. 

For the hierarchy over $h$ we have something similar, with the
geodesic $k$ supported in one of the complementary domains $Y_1$ of $v$,
and hence subordinate to $h$, but we have not constructed anything in
the domain $Y_2$. A geodesic in $Y_2$ does arise naturally, in the
following way. Let $U=S\setminus u$ and $W=S\setminus w$, noting that
both of these are two-holed tori containing $Y_2$. There are geodesics
$p$  supported in $U$ and $q$ supported in $W$, so that $p$ joins the
predecessor of $u$ to its successor $v$, and $q$ joins the predecessor
$v$ of $w$ to its successor (see figure \ref{r in Y2} for
schematic). Let $s$ be the vertex of $p$ preceding 
$v$, and let $t$ be the vertex of $q$ following $v$. Each is disjoint
from $v$, and therefore must lie in $Y_2$. We therefore may join $s$
to $t$ by a geodesic $r$ in $\CC(Y_2)$. In the notation we will later
develop, $r$ is {\em forward subordinate} to $q$, since it is
supported in the domain of $q$ minus the vertex $v$, and its last
vertex is the successor of $v$. Similarly $r$ is {\em backward
  subordinate} to $p$.

\realfig{r in Y2}{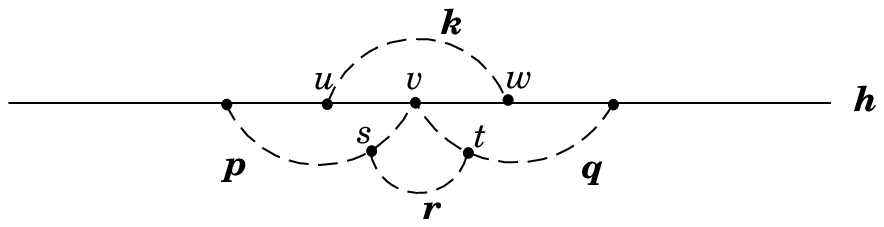}{The geodesic $r$, supported in $Y_2$, arises
  naturally after the geodesics $p$ and $q$ are constructed in the
  links of $u$ and $w$.}

Let us see how pants decompositions arise in this structure. In the
hierarchy over $h'$, the vertices $v'$ at bottom level (in $h'$), $v$ on
the next level (in $m'$), and any vertex $x$ in the
geodesic $k'$, form a pants decomposition, which we also call a {\em
  slice} of the hierarchy. If $x'$ is the successor of $x$ in $k'$
(so $x$ and $x'$ are neighbors in the Farey graph $\CC(Y_1)$),
the transition from $(v',v,x)$ to $(v',v,x')$ is an {\em elementary
  move}. 
 In the hierarchy over $h$ we can see a slice with different
organization: starting with $v$ at bottom level, we take any vertex
$a$ in $k$ and $b$ in $r$, and the triple $(v,a,b)$ make a pants
decomposition. We can move $a$ and $b$ independently in their
respective geodesics, since their domains ($Y_1$ and $Y_2$) are disjoint.
This kind of idea will give a way to ``resolve'' a hierarchy
(non-uniquely) into a sequence of slices, or markings, which will then
enable us to describe a useful class of words in the mapping class group.

In these examples we have only produced pants decompositions, but
in our final construction there will be complete markings, which
include twisting data around each pants 
curve. This will be done  using ``annulus complexes,'' which  are
analogous to the links of vertices in the Farey graph.

\subsection{Other applications and directions}
We hope that the tools developed here can be used to give an
algorithmic approach to $\Mod(S)$ in which the complexity of
problems such as the conjugacy problem can be computed. 
In particular, the conjugacy bound of Theorem \ref{Conjugacy Bound},
together with the quasi-geodesic words constructed from hierarchies,
are a good start provided that one can give an effective algorithm to {\em
construct} hierarchies with a Turing machine. 

The word problem, by comparison, admits a quadratic-time solution
because $\Mod(S)$ is known to have an  {\em automatic structure} (see Mosher
\cite{mosher:automatic}). A stronger condition known as a {\em
  biautomatic structure} (see \cite{epstein-et-al} for definitions of
these terms) would give bounds on the conjugacy problem,
but whether one exists remains open. Finding a biautomatic structure
was an initial motivation for this paper, but significant problems
remain. In particular the paths obtained from resolutions of
hierarchies are {\em not} a bicombing of $\Mod(S)$, because of the
presence of disjoint domains in $S$, whose order of traversal can
differ in different paths. The standard ``diagonalization'' method of moving in
both domains at once runs into some significant technical problems in
our setting. However, we believe that the hierarchy structure 
should be powerful enough by itself to give algorithmic results.

A rather different application of our ideas is to questions of
rigidity and classification for hyperbolic 3-manifolds. In
\cite{minsky:torus}, Kleinian representations of the fundamental group
of the punctured torus were studied via the length functions they
induce on its curve complex, the Farey graph. A connection between the
combinatorics of this graph and the geometry of the corresponding
3-manifolds was established, which was a primary ingredient in the
proof of Thurston's Ending Lamination Conjecture in that case. In
general, given a representation $\rho:\pi_1(S)\to PSL(2,\C)$ one can
study the complex translation lengths of conjugacy classes of simple
curves, viewed as a function on $\CC(S)$. In \cite{minsky:lengthfunctions}
some preliminary convexity properties are established for these
functions, which we hope will prove useful in studying the
general classification problem for Kleinian groups.

%% file: defs.tex
\section{Complexes and subcomplexes of curves}
\label{defs}

We review here the definitions of the various complexes of curves,
paying particular attention to the way in which subsurfaces of a given
surface give rise to sub-complexes. We will prove Lemma \ref{arcs
  to curves} relating arc complexes to curve complexes, define 
projections from a complex to its sub-complexes
and  prove Lemma \lref{Lipschitz Projection}.

We will also treat the
case of annuli, which are exceptional in various respects, and
conclude with a discussion of markings and elementary moves.

\subsection{Basic definitions and notation}
\label{basic defs}
Let $S=S_{\gamma,p}$ be an orientable surface of finite type, with genus
$\gamma(S)$ and $p(S)$ punctures.  It will be convenient to measure
the complexity of $S$ by $\xi(S) = 3 \gamma(S) + p(S)$.  Note that
$\xi$ is not equivalent to Euler characteristic, but has the property
that if $T\subset S$ is an incompressible proper subsurface then
$\xi(T)$ is strictly smaller than $\xi(S)$.  We will only consider
surfaces with $\xi > 1$, thus excluding the sphere and disk. We will
also exclude the standard torus (which does not arise as a subsurface
of a hyperbolic surface), so that from now on $\xi(S) = 3$ implies $S$
is the thrice-punctured sphere.

The {\em complex of curves} $\CC(S)$, introduced by Harvey
in \cite{harvey:boundary},
is a finite-dimensional and
usually locally infinite simplicial complex defined as follows:
A {\em curve} in $S$ is by definition a nontrivial homotopy class of simple
closed curves, not homotopic into a puncture.
If $\xi(S)>3$ then the set of curves is non-empty, and we let these be
the vertices of $\CC(S)$.
If $\xi(S)>4$ then the $k$-simplices are the sets $\{v_0,\ldots,v_k\}$
of distinct curves that have pairwise disjoint representatives.
One easily checks that $\dim(\CC(S)) = \xi(S) - 4$.

When $\xi(S)=4$, $S$ is either a once-punctured torus $S_{1,1}$ or four times
punctured sphere $S_{0,4}$, and the complex as defined above has
dimension 0. In 
this case we make an alternate definition: an edge in $\CC(S)$ is a
pair $\{v,w\}$ where $v$ and $w$ have representatives that intersect
once (for $S_{1,1}$) or 
twice (for $S_{0,4}$). Thus $\CC(S)$ is a graph, and in fact is
isomorphic to the familiar Farey graph in the plane (see e.g.
Bowditch-Epstein \cite{bowditch-epstein:triang}, 
Bowditch \cite{bowditch:markoff}, 
Hatcher-Thurston \cite{hatcher-thurston}, and
Series \cite{series:cfrac}). In particular this graph is a
triangulation of the 2-disk with vertices on the boundary, and the link
of each vertex can be identified with the integers, on which Dehn
twists (or half-twists for $S_{0,4}$) act by translation (see Figure
\ref{farey graph}).

When $\xi(S)=3$, $\CC(S)$ is empty (recall we have excluded the
regular torus). When $\xi(S)=2$, $S$ is the annulus and this case is
of interest when $S$ appears as a subsurface of a larger surface. We
consider this case further in \S\ref{annulus defs}.

\subsection{Distance geometry and hyperbolicity}
Let $\CC_k(S)$ denote the $k$-skeleton of $\CC(S)$. It is easy to show
that $\CC_k$ is connected for $k\ge 1$, 
see e.g. \cite[Lemma 2.1]{masur-minsky:complex1}. 
We can make $\CC_k(S)$ into a complete geodesic metric space by giving each
simplex the metric of a regular Euclidean simplex with side-length
1 (see Bridson \cite{bridson:simplicial}). It is easy to see that the
resulting spaces are quasi-isometric for all $k>0$. In
\cite{masur-minsky:complex1} we showed
\begin{theorem+}{Hyperbolicity}
If $\xi(S)\ge 4$ and $k>0$ then $\CC_k(S)$ is an infinite-diameter
$\delta$-hyperbolic metric space  for some $\delta>0$.
\end{theorem+}
See e.g.
\cite{cannon:negative,gromov:hypgroups,bowditch:hyperbolicity,ghys-harpe,short:notes}
for background on $\delta$-hyperbolic metric spaces. We recall here
just the definition that a geodesic metric space is $\delta$-hyperbolic if
for any geodesic triangle each edge is in a $\delta $-neighborhood of
the union of the other two edges.

We will usually consider just distances between {\em vertices} in
$\CC(S)$, i.e. points in $\CC_0(S)$, 
for which it suffices to consider distances in the graph $\CC_1(S)$,
which we note are integers. Thus by the notation $d_{\CC(S)}(v,w)$, or
even $d_S(v,w)$, we will always mean distances as measured in
$\CC_1(S)$. 
Writing $\diam_S$ to mean diameter in $\CC_1(S)$, 
we define for subsets $A,B\subset \CC_0(S)$
\begin{equation}
  \label{set distance is max}
d_S(A,B) = \diam_S(A\union B).  
\end{equation}

We will also usually think of a {\em geodesic} in $\CC_1(S)$ as a
sequence of vertices $\{v_i\}$ in $\CC_0(S)$, such that $d_S(v_i,v_j)
= |i-j|$. In particular $v_i$ and $v_{i+1}$ are always disjoint (when
$\xi(S) > 4$) and $v_i$ and $v_{i+3}$ always fill $S$, in the sense
that the union of the curves they represent, in minimal position, cuts
$S$ into a union of disks and once-punctured disks.

A final abuse of notation throughout the paper is in the usage of the
term ``vertex'':  when we introduce the notion of tight geodesics in
\S\ref{hierarchy defs} we will use ``vertex of a geodesic'' to denote
something more general than a point of $\CC_0(S)$, namely a simplex of
$\CC(S)$, representing a multi-component curve (or multicurve).
(One can think of this as a vertex of the first barycentric subdivision).
We will also go back and forth freely between vertices or simplices
and the (multi)curves they represent.

\subsection{Subdomains, links, arc complexes}
\label{subsurfaces}
A {\em domain} (or {\em subdomain}) $Y$ in $S$ will always be taken to mean
an (isotopy class of an) incompressible, non-peripheral, connected
open subsurface.  Unless we say {\em proper} subdomain, we include the
possibility that $Y=S$.
We usually omit the mention of isotopy classes for both surfaces and
curves, and to make the discussion clear one might fix a
complete hyperbolic metric on $S$ and consider geodesic
representatives of curves, and surfaces bounded by them.
We also take the word  
``intersection'' to mean {\em transverse} intersection.
However, annuli are an exceptional case in several ways;  see below.

In particular note that the boundary curves of a surface do not intersect
it. 

We immediately obtain an embedding $\CC(Y)\subset \CC(S)$ except when
$\xi(Y)\le 4$. Another complex of interest is the {\em arc complex}
$\CC'(Y)$,
which we define as follows: Suppose again that $\xi(Y)>3$.
An {\em arc} in $Y$ is 
a homotopy class of properly embedded paths in $Y$, which cannot be
deformed rel punctures to a point or a puncture. 
The vertices of $\CC'(Y)$ are both the arcs and the curves,
and simplices as before are sets of vertices that can be
realized disjointly. 
The complex $\CC'(Y)$  naturally arises when we try to ``project''
$\CC(S)$ into $\CC(Y)$ by taking intersections with $Y$ of curves in
$S$.  

{\bf Remark:} The punctures of $Y$ can come from either punctures of
$S$ or from boundary components of $Y$ in $S$. In fact,
it is often useful to think of all the punctures of $Y$ as
boundary components, in which case we consider arcs up to homotopy
which allows the endpoints to move on the boundary. These points of
view are equivalent, and we shall go back and forth between them for
convenience. 

The next elementary observation is that $\CC(Y)$ embeds in $\CC'(Y)$
as a co-bounded set. More precisely, letting
$\PP(X)$ denote the set of finite subsets of $X$, we have:

\begin{lemma}{arcs to curves}
Let $\xi(Y)>3$. There is a map $\psi=\psi_Y:\CC'_0(Y)\to \PP(\CC_0(Y))$
such that:
\begin{itemize}
\item
$\psi(v) = \{v\}$ for $v\in\CC_0(Y)$,
\item
$d_{\CC'(Y)}(\alpha,\psi(\alpha))\le 1$, and 
\item
if $d_{\CC'(Y)}(\alpha,\beta) \le 1$ then
$d_{\CC(Y)}(\psi(\alpha),\psi(\beta)) \le 2.$
\end{itemize}
\end{lemma}

\begin{pf}
If $\alpha$ is an arc,
let $\NN$ be a regular neighborhood in $Y$ of the union of $\alpha$
with the component(s) of $\boundary Y$ on which its endpoints lie, and
consider the frontier of $\NN$ in $Y$. This has either one or two
components, and at least one of them must be both nontrivial and
nonperipheral, since otherwise $Y$ is a disk, annulus or thrice-punctured
sphere, contradicting $\xi(Y)>3$. We let $\psi(\alpha)$ be the union
of the (at most two) nontrivial  components  
(see figure \ref{psi def fig}). If $\alpha$ is a curve (vertex of
$\CC_0(Y)$), we define $\psi(\alpha) = \{\alpha\}$.

\realfig{psi def fig}{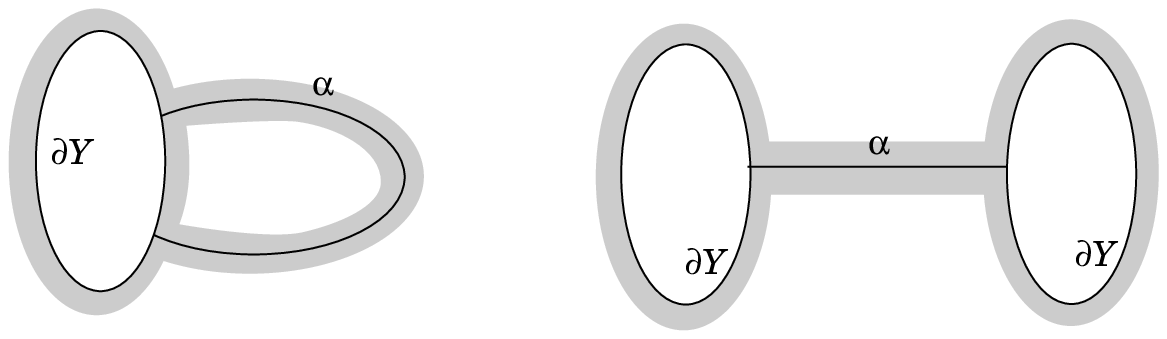}{The neighborhood $\NN$ is shaded. Note
that its frontier in $Y$ has two components in the first case and one
in the second}



Let $\alpha $ and $\beta$ be adjacent in $\CC'(Y)$, so they have
disjoint representatives. If either of them is a closed curve
then automatically $d(\psi(\alpha),\psi(\beta)) \le 1$, so assume both
are arcs. Similarly if their endpoints lie on disjoint boundary
components of $Y$ then $\psi(\alpha)$ and $\psi(\beta)$ have disjoint
representatives, so we can 
assume from now on that there is at least one boundary component which
touches both of them.
Suppose that the complement of $\alpha \union \beta$ in $Y$ contains a
non-trivial, non-peripheral
simple closed curve $\gamma$. Then $\gamma$ is also
disjoint from $\psi(\alpha)$ and $\psi(\beta)$, and we conclude
$d_{\CC(Y)}(\psi(\alpha),\psi(\beta)) \le 2$. 

If there is no such $\gamma$,
then $\alpha$ and $\beta$ cut $Y$ into a union
of (at most 3) disks or punctured disks. The possible cases 
can therefore be enumerated explicitly. 

\realfig{psi cases}{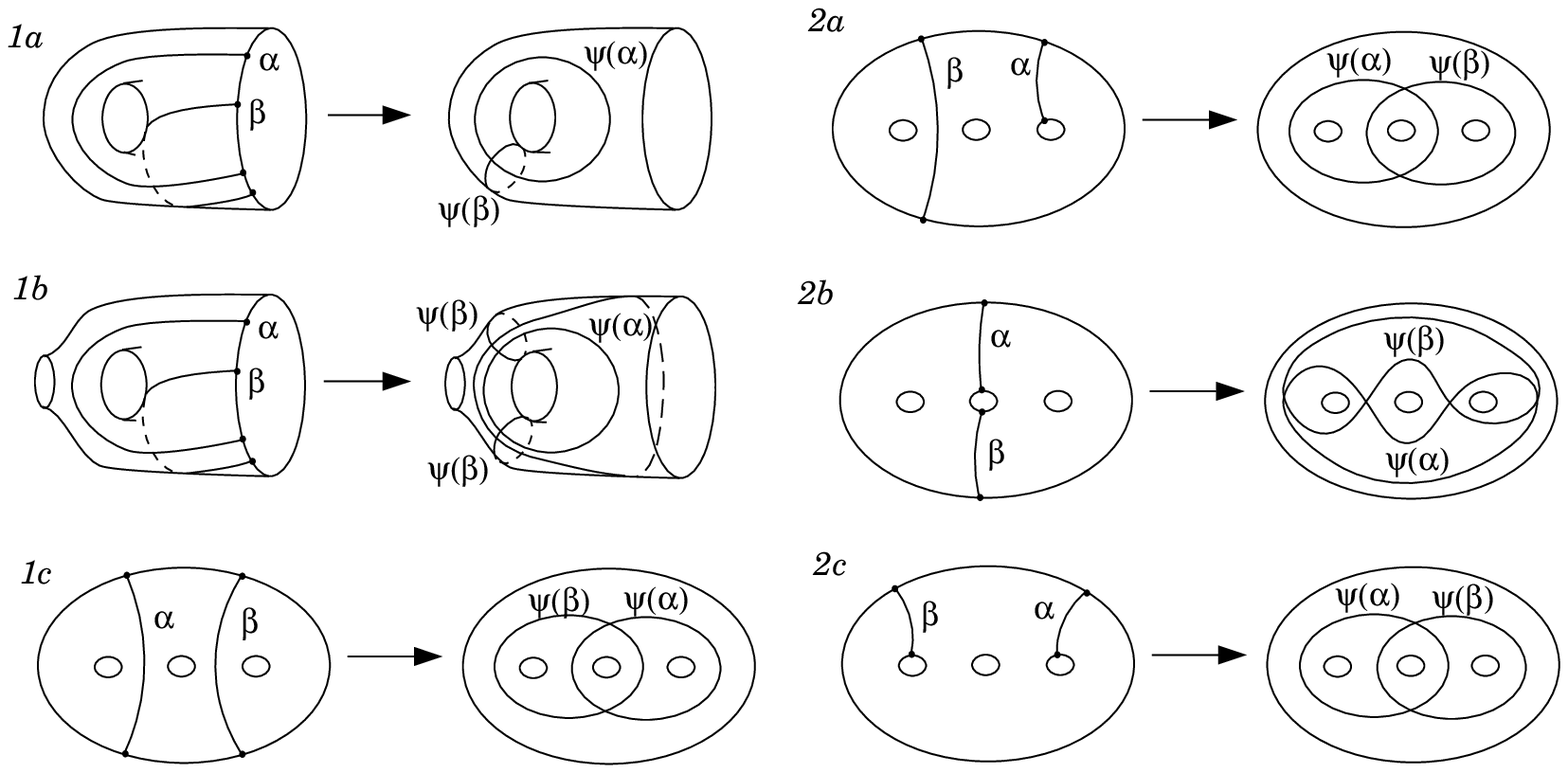}{The different cases in the proof
that $d_{\CC(Y)}(\psi(\alpha),\psi(\beta)) \le 2$.}

Let $C$ be a
boundary component of $Y$ meeting both $\alpha$ and $\beta$. 
If $C$ meets all the endpoints then there are two possibilities, 
according to whether the endpoints
separate each other on $C$. If they separate each other, 
$Y$ must be a once or twice
punctured torus, as in cases 1a and 1b of
Figure \ref{psi cases}. In case 1a we have
$d_{\CC(Y)}(\psi(\alpha),\psi(\beta)) = 1$, 
and in case 1b, $d_{\CC(Y)}(\psi(\alpha),\psi(\beta)) = 2$, as shown 
(note in this case that $\psi(\alpha)$ and $\psi(\beta)$ each have two
components). 
If they do not separate then 
$Y$ must be a quadruply-punctured sphere (case 1c) and
$d_{\CC(Y)}(\psi(\alpha),\psi(\beta)) = 1$.  
(Recall that in the cases where $\xi(Y)=4$, the definition of $d_{\CC(Y)}$ is
slightly different).

Suppose that $\alpha$ has one endpoint on $C$ and one on another
boundary $C'$. In all these cases $Y$ turns out to be a
quadruply-punctured sphere. 
If both $\beta$'s endpoints are on $C$ we get case 2a, where
$d_{\CC(Y)}(\psi(\alpha),\psi(\beta)) = 1$. If $\beta$'s other endpoint is on
$C'$ we get 2b, where $d_{\CC(Y)}(\psi(\alpha),\psi(\beta)) = 2$. If 
$\beta$'s other endpoint is on a third component $C''$, we get case 2c where
again $d_{\CC(Y)}(\psi(\alpha),\psi(\beta)) = 1$.
\end{pf}

\bfheading{Projections to subsurfaces:}
If $Y$ is a proper subdomain in $S$ with $\xi(Y)\ge 4$ we can define a map
$\pi'_Y : \CC_0(S) \to \PP(\CC'_0(Y))$, simply by taking
for any curve $\alpha$ the union of (homotopy classes of) its
essential arcs of intersection with $Y$. 
If $\alpha$ does not meet $Y$ essentially then
$\pi'_Y(\alpha)=\emptyset$, and otherwise it is always a simplex
of $\CC'(Y)$.

Adopting the convention for set-valued maps that $f(A) = \union_{a\in
  A}f(a)$, we define $\pi_Y$ by $\pi_Y(\alpha) = \psi_Y(\pi'_Y(\alpha))$.
We also define
\begin{equation}
  \label{dY convention}
  d_Y(A,B) \equiv d_Y(\pi_Y(A),\pi_Y(B))
\end{equation}
For sets or elements $A$ and $B$ in $\CC_0(S)$, and similarly
we let $diam_Y(A)$ denote $\diam_{\CC(Y)}(\pi_Y(A))$. 

\subsection{Annular domains} 
\label{annulus defs}
An annular domain is an annulus $Y$ with incompressible boundary in $S$,
which is not homotopic into a puncture of $S$. 
The purpose of defining complexes for such annuli is to keep track
of Dehn twisting around their cores; hence one would like $\CC(Y)$ to
be $\Z$. However there seems to be no natural way to do
this, and we will be content with something more complicated which is
nevertheless quasi-isometric to $\Z$. The
statements made in this subsection are all elementary, and we only
sketch the proofs.

Let $\til Y$ be the annular cover of $S$ to which $Y$ lifts
homeomorphically. There is a natural compactification of $\til Y$ to a
closed annulus $\hhat Y$, obtained in the usual way
from the compactification of the universal
cover  $\til S = \Hyp^2$ by the closed disk.
Define the vertices of $\CC(Y)$ to be
the paths connecting the two boundary components
of $\hhat Y$, modulo homotopies that {\em fix the endpoints}.
Put an edge between any two elements of $\CC_0(Y)$ which have
representatives with disjoint interiors. As before we can make
$\CC(Y)$ into a metric space with edge lengths 1.
If $\alpha\in\CC_0(S)$  is  the core curve of $Y$ we also write
$\CC(\alpha)=\CC(Y)$, and similarly $d_Y = d_{\alpha  }$.

Fixing an orientation on $S$ and
an ordering on the components of $\boundary \hhat Y$, we can define
algebraic intersection number $\alpha\cdot \beta$
for $\alpha,\beta\in\CC_0(Y)$ (only interior intersections count).
It is easy to see by an inductive argument that
\begin{equation}
  \label{annulus distance and intersection}
  d_Y(\alpha,\beta) = 1 + |\alpha\cdot \beta|
\end{equation}
whenever $\alpha\ne\beta$.
Let us also observe the convenient identity
\begin{equation}
  \label{adding twists}
  \gamma\cdot\alpha = \gamma\cdot\beta + \beta\cdot\alpha + j
\end{equation}
where $j=0,1$ or $-1$ (the value of $j$ depends on the exact
arrangement of endpoints on $\boundary \hhat Y$). 

We claim that $\CC(Y)$ is quasi-isometric to $\Z$ with the standard
metric. In fact define a map $f:\CC_0(Y) \to \Z$ by fixing some
$\alpha\in\CC_0(Y) $ and letting $f(\beta) = \beta\cdot\alpha$. Then
(\ref{adding twists}) and (\ref{annulus distance and intersection}) imply
\begin{equation}
  \label{f quasiisometry}
  |f(\gamma)-f(\beta)| \le d_Y(\gamma,\beta) \le |f(\gamma)-f(\beta)|
  + 2. 
\end{equation}
In particular this implies $\CC(Y)$ is hyperbolic so Theorem
\ref{Hyperbolicity} holds for this complex as well. 

\bfheading{Projections to annuli:}
We can define $\pi_Y:\CC_0(S) \to \PP(\CC_0(Y))$ as follows: 
If $\gamma$ is a simple closed curve in $S$ crossing the core of
$Y$ transversely, then the lift of $\gamma$ to $\til Y$ has at least
one component that connects the two boundaries of $\hhat Y$, and
together these  
components make up a (finite) set of diameter 1 in $\CC(Y)$.
Let $\pi_Y(\gamma) $ be this set.
If $\gamma$ does not intersect $Y$ essentially (including the case
that $\gamma$ is the core of $Y$!) then $\pi_Y(\gamma) = \emptyset$,
as in the previous section. 

Finally, for consistency we also define $\pi_Y:\CC_0(Y) \to \PP(\CC_0(Y))$
by $v\mapsto \{v\}$, and define $d_Y(A,B)$ and $\diam_Y(A)$ using the
same conventions 
(e.g. (\ref{set distance is max}) and (\ref{dY convention})) as for
larger subdomains. If $\alpha$ is the core of $Y$ we also write
$\diam_{\alpha}$ and $\pi_{\alpha}$.

We remark that $\CC(Y)$ is not a subcomplex of $\CC(S)$, but just as
for larger subdomains, any $f\in\Mod(S)$ acts by isomorphism
$f : \CC(Y) \to \CC(f(Y))$, and this fits naturally with the action on
$\CC(S)$ via $\pi_{f(Y)} \circ f  = f \circ \pi_Y$.

With these definitions in place we have the following:

\begin{lemma+}{Lipschitz Projection}
Let $Y$ be a subdomain of $Z$.
For any simplex $\rho$ in $\CC(Z)$, if $\pi_Y(\rho)\ne\emptyset$ then
$\diam_Y(\rho) \le 2$. If $Y$ is an annulus and $\xi(Z)>4$ then the
bound is 1.  
\end{lemma+}
\begin{proof}
  For an annulus $Y$, if $\xi(Z)>4$ the bound is immediate, since any
  two disjoint curves in   $Z$ lift to disjoint arcs in $\til Y$.
  If $\xi(Z) = 4$, one easily checks that Farey neighbors in $\CC(Z)$
  lift to curves that intersect at most once in any annulus cover.

  For $\xi(Y)\ge 4$, the bound follows from Lemma \ref{arcs to
  curves}.
\end{proof}

\bfheading{Dehn twists:}
Let $Y$ be an annulus with core $\alpha$. Let $D_\alpha$ be a positive
Dehn twist in $S$ about $\alpha$, and let $\hhat D_\alpha$ be a
positive Dehn twist in the covering annulus $\hhat
Y$ about its core.
Then $\hhat D_\alpha$ acts on $\CC(Y)$ and it is immediate 
for any $t\in\CC_0(Y)$ that $(\hhat D^n_\alpha t)\cdot t = n-1$ if $n>0$
and $n+1$ if $n<0$. Thus we obtain from
(\ref{annulus distance and intersection}) that
$d_Y(\hhat D_\alpha^n(t),t) = |n|$ for all $n\in\Z$

With a little more thought one can see that, for any curve $\beta$
intersecting $\alpha$ transversely, 
\begin{equation}
  \label{twist distance}
  d_Y(D^n_\alpha(\beta),\beta) = 2+|n|
\end{equation}
for $n\ne 0$. This is because the Dehn twist in $S$ affects every
intersection of the lift of $\beta$ with lifts of $\alpha$ in $\til
Y$, and this shifts the endpoints on $\boundary \hhat Y$ enough to
enable components of $\pi_Y(D^n_\alpha(\beta))$ and
$\pi_Y(\beta)$ to intersect an additional two times.

If $\beta$ intersects $\alpha$ exactly 2 times with opposite
  orientation, one can apply a {\em 
  half twist} to $\beta$ to obtain a curve $H_\alpha(\beta)$, which is
  equivalent to taking $\alpha\union\beta$ and resolving the
  intersections in a way consistent with orientation (see
  \cite{luo:intersectionnum}
  for a generalization). Then $H^2_\alpha(\beta)=D_\alpha(\beta)$, and
  one can also see for $n\ne 0$ that
  \begin{equation}
    \label{half-twist distance}
    d_Y(H^n_\alpha(\beta),\beta) = 2 + \left\lfloor\frac{|n|}{2}\right\rfloor.
  \end{equation}

\subsection{Markings}
\label{markings}
Assume $\xi(S)\ge 4$ and let
$\{\alpha_1,\ldots,\alpha_k\}$  be  some simplex in
$\CC(S)$. A {\em marking} in $S$ is a set $\mu=\{p_1,\ldots,p_k\}$,
where each $p_i$ is either just $\alpha_i$, or a pair $(\alpha_i,t_i)$
such that $t_i$ is a diameter-1 set of vertices of
the annular complex $\CC(\alpha_i)$.
The $\alpha_i$ are called the base curves and the simplex
$\{\alpha_i\}$ is denoted
$\base(\mu)$. The (possibly empty) set $\{t_i\}$ is called the
set of {\em transversals} and denoted $\trans(\mu)$.
Thus a special case is when $\trans(\mu)=\emptyset$ and then
$\mu=\base(\mu)$.

If $\base(\mu)$ is contained in $\CC(Y)$ for some non-annular
subsurface in $Y$, we call $\mu$ a marking in $Y$.
If $Y$ is an essential annulus in $S$ then a marking $\mu$ in $Y$ is any set
of diameter 1 in $\CC_0(Y)$ (typically these sets will have at most
two elements), and we have $\mu=\base(\mu)$ in this case.

If $\base(\mu)$ is maximal  and 
every curve has a transversal,
the marking is called {\em complete.}

Markings can be very complicated objects, because the transversals,
being arcs in annular covers, can have complicated images in $S$. 
Let us therefore define something called a {\em clean marking}:

Given $\alpha\in\CC_0(S)$ a {\em clean transverse curve} for $\alpha$
is a curve $\beta\in\CC_0(S)$ such that
a regular neighborhood of $\alpha\union \beta$  (in minimal position) is a
surface $F$ with $\xi(F)=4$, in which $\alpha$ and $\beta$ are
$\CC(F)$-neighbors (note there are only two possible configurations,
corresponding to $F$ being a 1-holed torus or 4-holed sphere, and
$\alpha$ and $\beta$ intersect once or twice, respectively).

A marking $\mu$ is called {\em clean} if every pair in $\mu$ is of the
form $(\alpha_i,\pi_{\alpha_i}(\beta_i))$ where 
$\beta_i$ is a clean transverse curve for $\alpha_i$,
which also misses the other curves in $\base(\mu)$.
Note that if $\mu$ is clean then the curves $\beta_i$ are uniquely
determined by the transversals $t_i=\pi_{\alpha_i}(\beta_i)$.

We note that, up to homeomorphisms of $S$, there are only a finite
number of clean markings.

If $\mu$ is a complete marking, there is an almost
canonical way to select a related clean marking. Let us say
that a clean marking $\mu'$ is {\em compatible} with  a marking $\mu$ 
provided $\base(\mu) = \base(\mu')$, a base curve $\alpha$ has a
transversal $t'=\pi_Y(\beta)$ in $\mu'$ if and only if it has a
transversal $t$ in 
$\mu$, and $d_{\alpha}(t,t')$ is minimal among all possible choices of $t'$.

\def\ncompat{n_0}
\def\dcompat{n_1}
\begin{lemma}{clean markings}
Let $\mu$ be a complete marking of $Y\subset S$. Then there exist at
least 1 and at most $\ncompat^b$ complete clean markings $\mu'$
compatible with $\mu$, where $b$ is the number of base curves of
$\mu$, and $\ncompat$ is a universal constant.
Furthermore,  
for each $(\alpha,t)\in \mu$ and
$(\alpha,t')\in \mu'$ we have $d_\alpha(t,t') \le \dcompat$, where
$\dcompat$ is a universal constant.
\end{lemma}

\begin{proof} 
Fix one clean marking $\mu_0$ with $\base(\mu_0)=\base(\mu)$.
All other clean markings with this base are obtained from $\mu_0$ by
twists and half-twists, so it follows immediately from 
(\ref{twist distance},\ref{half-twist distance}) and the
quasi-isometry (\ref{f quasiisometry}) of an annular complex to $\Z$ 
that for each
$\alpha\in\base(\mu)$ there is a choice 
of clean transversal $\beta$  that minimizes
$d_\alpha(t,\pi_\alpha(\beta))$, and that there is a uniform bound on
this minimum.
The fact that the number of choices of $\beta$ are uniformly bounded
for each base curve also follows from (\ref{twist
  distance}) and (\ref{half-twist distance}).
\end{proof}
(One can in fact show that $\ncompat\le 4$ and $\dcompat=3$, but we
will not need this).

\bfheading{Projections of markings:}
If $Y$ is any subdomain of $S$ and $\mu$ any marking in $S$ 
we can define $\pi_Y(\mu)$ as follows: 
If $Y$ is an annulus whose core is some $\alpha\in\base(\mu)$, and
$\alpha$ has a transversal $t$, we
define $\pi_Y(\mu) = t$. If $\alpha$ has no transversal $\pi_Y(\mu) =
\emptyset$. In all other cases, $\pi_Y(\mu) = \pi_Y(\base(\mu))$.

\bfheading{Elementary moves on clean markings:}
Let $\mu$ be a complete clean marking, with 
pairs $(\alpha_i,\pi_{\alpha_i}(\beta_i))$  as above.
There are two types of elementary moves that transform $\mu$ into a
new clean marking. 

\begin{enumerate}
\item Twist: Replace $\beta_i$ by
   $\beta'_i$, where $\beta'_i$ is obtained from $\beta_i$
  by a Dehn twist or half-twist around $\alpha_i$.
\item Flip: Replace $(\alpha_i,\pi_{\alpha_i}(\beta_i))\in\mu$ by
  $(\beta_i,\pi_{\beta_i}(\alpha_i))$ to get a 
  non-clean marking $\mu''$.
 Then replace $\mu''$ by a compatible clean marking $\mu'$.
\end{enumerate}
In the first move a twist can be positive or negative. A
half-twist is  possible when $\alpha_i$ and $\beta_i$ intersect twice.

The replacement part of the Flip move requires further discussion:
The surface $F$ filled by $\alpha_i$ and $\beta_i$ has $\xi(F)=4$, and its
(non-puncture) boundary components are other elements of
$\base(\mu)$. For each such element $\alpha_j$ there is a transverse 
$\beta_j$ which misses $\alpha_i$ but hits $\beta_i$. Thus after interchanging
$\alpha_i$ and $\beta_i$ the marking is no longer clean.
We must therefore replace $\beta_j$ by $\beta'_j$ which
misses $\beta_i$, subject 
to the condition that $d_{\alpha_j}(\beta_j,\beta'_j)$ is as small as possible.
Lemma \ref{clean markings} says that this distance is at most
$\dcompat$, and there are $\ncompat$ possible choices for each $\beta_j$.
(Actually this is a more special case than Lemma \ref{clean markings}
and one can get a distance bound of 2).

Thus, given $\mu$ there is a finite number of possible elementary
moves on it, depending only on the topological type of $S$.

We conclude with an extension of Lemma \lref{Lipschitz Projection}.
\begin{lemma+}{Elementary Move Projections}
If $\mu,\mu'$ are complete clean markings differing by one elementary
move, then for any domain $Y$ in $S$ with $\xi(Y)\ne 3$,
$$
d_Y(\mu,\mu') \le 4
$$
If $Y$ is an annulus the bound is 3.
\end{lemma+}
\begin{proof}
If $Y$ is an annulus  with core curve $\alpha\in\base (\mu)$, then
$\mu $ contains $(\alpha,\pi_\alpha(\beta))$ for a clean transversal curve
$\beta$, and $\pi_Y(\mu) = \pi_\alpha(\beta)$. Then if
$\mu'$ is obtained 
by a twist or half-twist on $\alpha$,  a bound of 3 follows from
(\ref{twist distance}) and (\ref{half-twist distance}). If $\mu'$ is
obtained by a Flip move, replacing $(\alpha,\pi_\alpha(\beta))$ by
$(\beta,\pi_\beta(\alpha))$, then
$\pi_Y(\mu')=\pi_Y(\base(\mu'))=\pi_\alpha(\beta)$, so the distance is 0.

A similar analysis holds if $Y$ is an annulus with core curve in
$\base(\mu')$. 

In all other cases, $\pi_Y(\mu)=\pi_Y(\base(\mu))$ and
$\pi_Y(\mu')=\pi_Y(\base(\mu'))$,
and by definition 
$d_Y(\mu,\mu') = \diam_Y(\pi_Y(\base(\mu))\union\pi_Y(\base(\mu'))$.
If $\pi_Y(\base(\mu))$ and $\pi_Y(\base(\mu'))$ have at least one curve in
common, the bound of 4 follows from Lemma \ref{Lipschitz Projection}.
If not, then the move must be a Flip move, and $Y$ meets only the two
base curves $\alpha,\alpha'$ involved in the Flip. Let $F$ be the
surface of $\xi=4$ 
filled by these curves, which are neighbors in $\CC(F)$.
If $\xi(Y)=4$ then $Y=F$ and we are done, with a bound of 1.
The remaining possibility is that $Y$ is an essential annulus in $F$ meeting
both curves, and then any two lifts of $\alpha$ and $\alpha'$ to 
$\til Y$ intersect at most once, giving a bound of 2.
\end{proof}

%% file: projection.tex
\section{Projection bounds}
\label{projection}

Our goal in this section will be to 
prove Theorem \ref{Bounded Geodesic Image},
which gives strong 
contraction properties for the subsurface projections $\pi_Y$.

\begin{theorem+}{Bounded Geodesic Image}
Let $Y$ be a proper subdomain of $Z$ with $\xi(Y)\ne 3$
and let $g$ be a geodesic segment, ray,
or biinifinite line in $\CC(Z)$, such that $\pi_Y(v)\ne\emptyset$ for
every vertex $v$ of $g$.

There is a constant $M$ depending only on $\xi(Z)$ so
that 
$$\diam_Y(g) \le M.$$
\end{theorem+}

The intuition behind the statement is this: as we move in one
direction in $g = 
\{...v_1,v_2,...\}$, we expect the 
vertices to converge to some foliation in $Z$. Hence their projections
to $Y$ should converge to the intersection with $Y$ of the foliation
leaves. Recalling that $\pi_Y$ identifies parallel arcs, it should
follow that eventually $\pi_Y(v_i)$ should stabilize to a finite
collection of possible arcs. 
To make this precise we have to re-introduce the tools of
Teichm\"uller geometry from \cite{masur-minsky:complex1}.
We also emphasize that the statements we prove will be strictly weaker
than this intuitive description, but will suffice for the diameter bound.

\subsection{Quadratic differentials, vertical and horizontal}
Given a finite-type complex structure on $Z$, recall that
a holomorphic quadratic differential $q$ on $Z$ is a
tensor of the form $\varphi(z)dz^2$ in local coordinates, with
$\varphi$ holomorphic. Away from zeroes, a coordinate $\zeta$ can be
chosen so that $q = d\zeta^2$, which determines a Euclidean metric
$|d\zeta^2|$ together with a pair of orthogonal foliations parallel to
the real and imaginary axes in the $\zeta$ plane. These are
well-defined globally and are called the {\em horizontal} and {\em
vertical} foliations, respectively. The zeroes of $q$ are cone points
with cone angle $n\pi$, $n\in\Z$, $n\ge 2$.
(See Gardiner \cite{gardiner} or
Strebel \cite{strebel}.) 

For a closed curve or arc $\alpha$ in $Z$, denote by
$|\alpha|_q$ its length in the $q$ metric. Let $|\alpha|_{q,h}$
and $|\alpha|_{q,v}$ denote its horizontal and vertical lengths,
respectively, by which we mean the total lengths  of the (locally
defined) projections of $\alpha$ to the horizontal and vertical
directions of $q$. 

Henceforth assume $q$ has finite area, which means that at the
punctures it has poles of order 1 or less, and equivalently that its
metric completion gives a surface $\hat Z$ which is $Z$ with a cone point
added at each puncture, with cone angle $n\pi$, $n\in\Z$, $n\ge 1$.

Define a {\em straight segment} to be a path in $Z$ which meets no
punctures or zeroes of $q$, and is a straight line in the Euclidean metric.
A geodesic is composed of a finite number of straight segments,
meeting at zeroes with a certain angle condition.
We must slightly generalize the notion of ``geodesic representative''
as follows:
If $Z$ has punctures, the incompleteness of $|q|$ means that
a non-peripheral homotopy class $\alpha$ may not have a geodesic
representative. However, there is a representative in $\hat Z$ which
goes through 
the punctures some finite number of times and is geodesic elsewhere,
which we can think of as a limit of geodesic representatives in the
compact surfaces obtained by deleting open disks of $q$-radius $r$
around the punctures, for $r\to 0$. Thus by ``geodesic
representatives'' we will in fact mean representatives in this sense.

Let $\ep,\theta>0$ be some fixed (small) constants. We say that a
straight segment $\alpha$ is {\em almost vertical} with respect to $q$
if if it makes an angle of at most $\theta$ with the vertical
foliation. We say a geodesic is almost vertical  if it is composed of
straight segments meeting at punctures or zeros, each of which is
either almost vertical, or has 
length at most $\ep$. We define {\em almost horizontal} in the
analogous way.

\begin{lemma}{almost verticals close}
There is a choice of $\ep,\theta$ depending only on $\xi(Z)$ such that
the following holds. 
Let $Y$ be a domain in $Z$ with $\xi(Y)\ne 3$, $q$ a unit-area quadratic
differential on $Z$, and 
$\alpha$ a boundary component of $Y$ whose $q$-geodesic contains an
almost-horizontal segment $\sigma$ of horizontal length 1. Then
if $\beta$ and $\gamma$ are two almost-vertical curves intersecting
$Y$,
$$
d_Y(\beta,\gamma) \le 4.
$$
\end{lemma}

\begin{pf}
We begin with the case 
where $Y$ is not an annulus. 

For simplicity, suppose first that $Y$ is isotopic to an embedded surface with
$q$-geodesic boundary. Thus we may assume that $\alpha$ is already geodesic.
Consider the flow starting from $\sigma$ and
moving along the vertical foliation into $Y$ until it returns to
$\sigma$ or meets a singularity. The points corresponding to flow
lines that meet singularities divide $\sigma$ into at most $k_0$
intervals $\{I_j\}$, where $k_0$ depends only on $\xi(Z)$. Each $I_j$
determines a ``flow rectangle'', which is actually a Euclidean
trapezoid or parallelogram with two vertical sides and two almost-horizontal
sides which have slope at most $\tan \theta$. The interior of the
rectangle is embedded, though its top and bottom edges are segments of
$\sigma$ that may overlap.
Since $\sigma $ has horizontal length 1 there must be an interval $I_j$
of horizontal length at 
least $1/{k_0}$. Let $R$ denote the corresponding flow rectangle, and
$h$ the average height of $R$. Thus $R$ has
area  at least $h/k_0$, and
since $q$ is unit-area, $h\le k_0$.

Suppose that $\ep < 1/k_0$ and that $\theta$ is sufficiently small
that $\cot \theta > k_0^2 + \tan \theta$.

With these choices, we claim that
an almost-vertical geodesic $\beta$ cannot cross $R$ from one
vertical  side to the other: Since $R$ has no singularities in its
interior, such 
a crossing would have to be a straight segment $\tau$, and
the slope of $\tau$ would be at most 
$h k_0  + \tan \theta \le k_0^2 + \tan \theta $, which is less than
$\cot\theta$ by the choice of $\theta$.
Hence $\tau$ could not be almost vertical. Thus it would have
to have length bounded by $\ep$, and hence be shorter than the width
of $R$, again a contradiction.

We conclude that $\beta$ is disjoint from the interior of some arc $a$
in $R$ connecting the top edge
and the bottom edge. Thus, any component of $a\intersect Y$ gives
an element of $\CC'(Y)$ which is distance 1 from each vertex of
$\pi'_Y(\beta)$. 
The same argument 
applies to $\gamma$, with $a$ in the same homotopy class, 
and we conclude $d_{\CC'(Y)}(\pi'_{Y}(\beta),\pi'_{Y}(\gamma))\le 2.$
Lemma \ref{arcs to curves} then gives the desired bound.

Now consider the possiblity that $Y$ is not homotopic to an embedded surface
with geodesic boundary $\alpha$. In particular the geodesic
representative of $\alpha$ may traverse one or more geodesic segments
more than once, producing arcs of self-tangency.
However even in this case we obtain a map of $Y$ into $\hat Z$
($Z$ union its punctures)
which is homotopic to the inclusion by a homotopy 
that is an embedding until the last moment.
At that moment families of arcs in $Y$ or its complement, 
with endpoints on $\boundary Y$, are collapsed to points, producing the
arcs of self-tangency.  
It is easy to see that the same argument holds except that the
rectangles $R$ in question may have height zero, with horizontal arcs
on the self-tangencies and vertical arcs collapsed.

This concludes the case where $Y$ is not an annulus. 

\medskip

When $Y$ is an annulus, $\alpha$ is in the homotopy class of its core.

Lift $\alpha$ to $\til\alpha$ in the universal cover $\til Z$. We
remind the reader that again the geodesic representative of $\alpha$
may pass through punctures, and as the universal covering is
infinitely branched around punctures the topology is easier to keep
track of if we keep $\alpha$ outside a small neighborhood of the
punctures. At any rate our segment $\sigma$ can be assumed disjoint
from the punctures so we need not worry about this. 

If we consider the lines of the vertical flow which start at $\sigma$ and
go in both directions until they hit $\sigma$ again, we obtain at most $2k_0$
rectangles composed of vertical flow lines with $\sigma$ passing
through them, and we choose $R$ to be one which has width
at least $1/2k_0$. 
Let $\{R_n\}_{n\in\Z}$ denote its lifts
corresponding to the lift of $\alpha$ to $\til\alpha$, so that
$\til\alpha $ passes through the interior of each $R_n$, and the top
and bottom edges of $R_n$ lie on translates of $\alpha$ called
$\til\alpha_n$ and $\til\alpha'_n$, respectively. 
(We include also the degenerate possibility that $R$ has height 0 on one
side or the other of $\alpha$
and so the
$\til\alpha_n$, or $\til\alpha'_n$, are each tangent to $\til\alpha$ along
a segment.)

\realfig{twist projection}{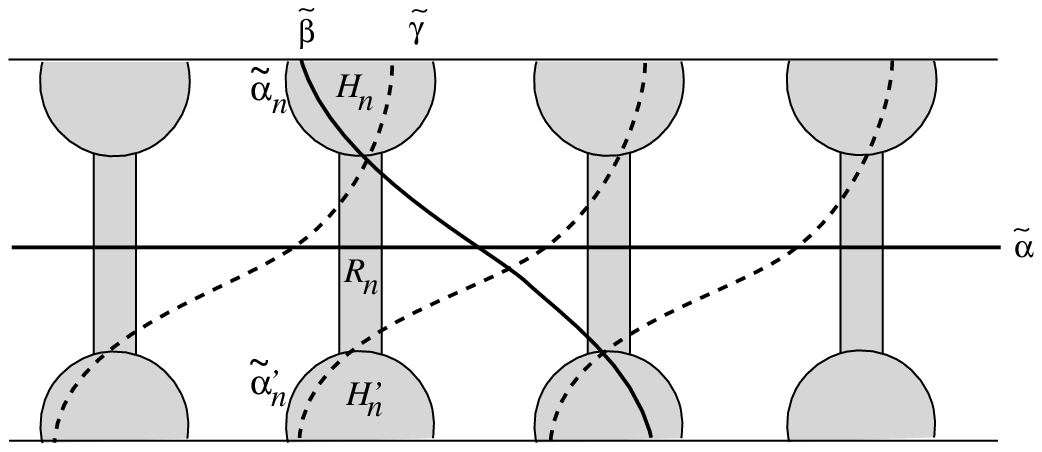}{}

After an arbitrary choice of orientation for $\alpha$, each $R_n$
has a left and a right vertical 
edge. Let $\til\beta$ and $\til\gamma$ be components of the lifts
of $\beta $ 
and $\gamma$ which cross $\til\alpha$. As argued before, and with
appropriate choice of $\ep,\theta$, neither
$\til\beta$ nor $\til\gamma$ can cross a rectangle $R_n$ from left to
right.
Let $H_n$ and $H'_n$ denote
the halfplanes bounded by
$\til\alpha_n$ and $\til\alpha'_n$, respectively, whose interiors  are
disjoint from 
$\til\alpha$ (see figure \ref{twist projection}).
Let $U_n$ denote $R_n\union H_n \union H'_n$.
Then neither $\til\beta$ nor $\til\gamma$
can cross through $U_n$ from left to right, because this would involve
either crossing $R_n$ from left to right, or entering and exiting the
interior of $H_n$ or $H'_n$, which a geodesic cannot do.

Thus if $\rho_n$ is the right-hand boundary of $U_n$, $\til\beta$
and $\til\gamma$ can each
cross at most one of the $\rho_n$. If $\rho$, $\hhat\beta$ and
$\hhat \gamma$ are the covering projections of $\rho_n$, $\til\beta$ and
$\til\gamma$, respectively, to the annulus $\hhat Y$, then we obtain
$|\hhat\beta\cdot\rho| \le 1 $ and $|\hhat\gamma\cdot\rho| \le 1 $.
It follows by (\ref{adding twists}) that $|\hhat\beta\cdot\hhat\gamma|\le 3$
and hence $d_Y(\beta,\gamma) \le 4$ by 
(\ref{annulus distance and intersection}). 
\end{pf}

\subsection{Teichm\"uller geodesics and balancing.}
A Teichm\"uller geodesic in $\TT(Z)$ ``shadows'' a $\CC(Z)$-geodesic
in the following 
specific sense, which played a crucial role in
\cite{masur-minsky:complex1}.

Recall that a
Teichm\"uller geodesic $L:\R \to \TT(Z)$ can be described in terms of
a family of quadratic differentials $q_t$ holomorphic on $L(t)$: Each
$q_t$ is obtained from $q_0$ by scaling the horizontal directions by
$e^t$ and the vertical by $e^{-t}$. This determines the conformal
structure $L(t)$. 

In \cite{masur-minsky:complex1}, we associate to the geodesic $L$ a
map $F:\R\to \CC_0(Z)$ by letting $F(t)$ be any simple curve of minimal
extremal length with respect to  $L(t)$. Furthermore we define a
map $\pi=\pi_q: \CC_0(Z) \to \R\union\{\pm\infty\}$, called a
``balancing projection,'' as follows:
Given any $\alpha\in\CC_0(Z)$, its horizontal length $|\alpha|_{q_t,h}$
has the form $|\alpha|_{q_0,h} e^t$ and
its vertical length $|\alpha|_{q_t,v}$
has the form $|\alpha|_{q_0,v} e^{-t}$. Thus if both of these are non
zero there is a unique point $t$ where they are equal, and we say
$\alpha$ is {\em balanced} at $t$, and set
$\pi_q(\alpha)=t$. If the horizontal lengths are 0 ($\alpha$ is parallel
to the vertical foliation) then we let
$\pi_q(\alpha) = +\infty$, and if the vertical lengths are $0$ we let
$\pi_q(\alpha) = -\infty$.

Now suppose we are given $v,w\in\CC_0(Z)$ with $d_{\CC(Z)}(v,w)\ge
3$. Then $v$ and $w$ fill $Z$, and so there is a conformal structure
and a quadratic differential
$q_0$ for which the  horizontal and vertical foliations
have closed nonsingular leaves which are isotopic to $v$ and $w$,
respectively. The corresponding Teichm\"uller geodesic is called the
Teichm\"uller geodesic associated to $(v,w)$. We note immediately
that $\pi_q(\alpha) = -\infty$ for $d(\alpha,v)\le 1$
and similarly that 
that $\pi_q(\alpha) = +\infty$ for $d(\alpha,w)\le 1$.

Some basic properties of this projection map are outlined in the
following lemma. Here $d()$ 
and $\diam()$ refer to 
distance and diameter in
$\CC_1(Z)$. The $K_i$ are constants depending only on $\xi(Z)$.
The notation $[s,t]$ refers to the interval with endpoints $s$ and
$t$, regardless of order.

\begin{lemma}{teichmuller properties}
Let $g = \{v_i\}_{i=M}^{N}$ be a geodesic segment in $\CC(Z)$ with
$M-N\ge 3$, and let $L:\R\to\TT(Z)$  be the Teichm\"uller geodesic
associated to $(v_M,v_N)$,
$F:\R\to\CC_0(Z)$ its associated map and $\pi:\CC_0(Z)\to\R$ the
associated projection. There are constants $K_0,K_1,K_2,m_0>0$, depending
only on the surface $Z$, such that:
\begin{enumerate}
\item \label{lipschitz}
(Lipschitz)
If $d(v,w) \le 1$ for $v,w\in\CC_0(Z)$ then 
$$\diam(F([\pi(v),\pi(w)])) \le K_0.$$
\item \label{Fellow traveling 1}
(Fellow traveling 1) for any $v_i$ in $g$, 
$$d(v_i,F(\pi(v_i))) \le K_1$$
\item \label{Fellow traveling 2}
(Fellow traveling 2)
For all $t\in\R$, there exists some $v_i\in g$ such that
$$\diam(F([t,\pi(v_i)])) \le K_2$$
\item \label{coarse monotonicity}
(Coarse monotonicity) 
Whenever $v_i$, $v_j$ are in $g$ with $j> i+m_0$, 
$$\pi(v_j)>\pi(v_i)$$
\end{enumerate}
\end{lemma}

\begin{pf}
Part (\ref{lipschitz}) is part of Theorem 2.6 of 
\cite{masur-minsky:complex1}, and parts
(\ref{Fellow traveling 1}) and 
(\ref{Fellow traveling 2}) 
follow from Theorem 2.6 together with 
the proof of Lemma 6.1 of \cite{masur-minsky:complex1}.

Part (\ref{coarse monotonicity}) is a consequence of parts
(\ref{lipschitz}) and (\ref{Fellow traveling 1}):
Since the last vertex $v_N$ has $\pi(v_N) = +\infty$ by definition,
if we have $\pi(v_{i+m}) < \pi(v_i)$ 
then there is some $m'\ge m$ for which $\pi(v_i) \in
[\pi(v_{m'}),\pi(v_{m'+1})]$. 
By (\ref{lipschitz}), we then have
$d(F(\pi(v_i)),F(\pi(v_{i+m'}))) \le K_0$. However since 
$d(v_i,v_{i+m'}) = m'$, this implies together with (\ref{Fellow
traveling 1}) and the triangle inequality that $m'\le K_0 + 2K_1$.
Setting $m_0 = K_0 + 2K_1$, we have part (\ref{coarse monotonicity}).
\end{pf}

\subsection{Proof of Theorem \ref{Bounded Geodesic Image}}
Let us first consider the case where $g$ is a finite segment
$\{v_i\}_{i=M}^N$.
Note that at most 3 of the vertices can actually be
contained $Y$, since they would all be $\CC(Z)$-distance 1 from
$\boundary Y$.
We may assume without loss of generality that $|g|=N-M \ge 3$.

Select a Teichm\"uller geodesic $L: \R \to \TT(Z)$ as above,
associated to $(v_M,v_N)$, as well as the associated map $F$,
family of quadratic differentials $q_t$, and balancing map $\pi$.

Let $\alpha$ be any boundary component of $Y$ (non-peripheral in $Z$).
Let $s_0 = \pi(\alpha)$. Note that possibly $s_0 =
-\infty$ or $+\infty$, if 
$\alpha$ is disjoint from  $v_M$ or $v_N$ (but not both).

If $s_0\ne \pm \infty$, 
then $\alpha$ is balanced at $s_0$. If $s_0=-\infty$ then it is
horizontal at any $q_t$. In either case, 
Lemmas 5.3 and 5.6 of
\cite{masur-minsky:complex1} imply that, for $K_3>0$ depending only on
$Z$, there is $s_1\ge s_0$ with 
\begin{equation}
\label{s0 s1 bound}
\diam(F[s_0,s_1]) \le K_3,
\end{equation}
such that $\alpha$ is almost-horizontal
with respect to $q_{s}$ whenever $s\ge s_1$, and contains an almost
horizontal segment of 
horizontal length $\ep_1$, for some fixed $\ep_1>0$. In fact
we may assume $\ep_1=1$, because horizontal length expands at a
definite exponential rate with distance along the Teichm\"uller
geodesic $L$, and the map $F$ is quasi-Lipschitz
by Lemma 5.1 of \cite{masur-minsky:complex1}.
(The case $s_0=\infty$ is treated similarly, interchanging horizontal
and vertical).

%

Lemma 5.7 of \cite{masur-minsky:complex1} implies that, for $K_4>0$
depending only on $Z$,
there exists $s_2>s_1$ such that 
\begin{equation}
\label{s1 s2 bound}
\diam(F[s_1,s_2])\le K_4
\end{equation}
and, 
for any $\gamma\in\CC(Z)$, if $\pi(\gamma)>s_2$ then $\gamma$ is
almost vertical with respect to $q_{s_1}$. 
Again, possibly $s_2 = \infty$.

Let $j_0$ be the index of the vertex of $g$ for which part 
(\ref{Fellow traveling 2}) of Lemma \ref{teichmuller properties}
gives
\begin{equation}
\label{j0 s0 bound}
\diam(F([s_0,\pi(v_{j_0})])) \le K_2.
\end{equation}
We will now show that for $i>j_0$ sufficiently large, 
$\pi(v_i) > s_2$. 

By the coarse monotonicity (\ref{coarse monotonicity}) of Lemma
\ref{teichmuller properties}, if $i>j_0+m_0$
then $\pi(v_i) > \pi(v_{j_0})$.
Thus if $\pi(v_i)\le s_2$, we have $d(F(\pi(v_i)),F(\pi(v_{j_0})))\le
\diam(F([\pi(v_{j_0}),s_2]))$, and the latter is bounded by
$K_2 + K_3 + K_4$ because of the bounds (\ref{s0 s1 bound}), (\ref{s1
s2 bound}) and (\ref{j0 s0 bound}).
Thus, $i-j_0 = d(v_i,v_{j_0}) \le K_2 + K_3 + K_4 + 2K_1$ by (\ref{Fellow
traveling 1}) of Lemma \ref{teichmuller properties} and the triangle
inequality. Letting 
$m_1 = 1+ \max (m_0,K_2 + K_3 + K_4 + 2K_1)$, we are therefore assured
that if $i\ge j_0 + m_1$ then $\pi(v_i)>s_2$.

Thus, if $i\ge j_0 + m_1$, then $v_i$ is almost vertical with
respect to $q_{s_1}$. 

We can now apply Lemma \ref{almost verticals close} using the
quadratic differential $q_{s_2}$ and the boundary component
$\alpha$. If $j_0+m_1 \le N$ then
for any $i,i'\in [j_0+m_1,N]$ we have by the above that both $v_i$
and $v_{i'}$ are almost vertical with respect to $q_{s_2}$, and thus
$d_Y(v_i,v_i')\le 4$.

The same argument, with horizontal and vertical interchanged, applies
to give a bound for $i,i'\in[M,j_0 - m_1]$, if $M\le j_0-m_1$. The
remaining segment between
$\max(M,j_0-m_1)$ and $\min(N,j_0+m_1)$ has a diameter bound of
$2m_1$, so its $\pi_Y$-image has diameter at most $4m_1$ by Lemma
\lref{Lipschitz Projection}. Thus the image 
of the full segment $g$ is bounded by $4m_1 + 8$.

Since this bound is independent of $N$ and $M$, it implies a bound
also in the infinite cases, via an exhaustion of $g$ by finite subsegments.
This concludes the proof of Theorem \ref{Bounded Geodesic Image}.

%% file: paths.tex
\section{Tight geodesics and hierarchies}
\label{hierarchies}

This section describes the main construction of our paper, hierarchies
of tight geodesics. After defining these notions in \S\ref{hierarchy
  defs}, we prove some existence results, Lemma \ref{Tight geodesics
  exist} and Theorem \ref{Hierarchies exist}, in \S\ref{existence results}.

Hierarchies give us the combinatorial framework in which to carry out
the link projection arguments first outlined in the examples in
\S\ref{example} (and done in generality in Section \ref{large link
  etc}). 
The main ingredient in this is the {\em backward and forward
  sequences} $\Sigma^\pm$, whose basic structural properties are
stated in Theorem \ref{Structure of Sigma}. The proof of this theorem
takes up the rest of Section \ref{hierarchies}, and along the way we
will develop a number of results, notably Theorem \ref{Completeness}, 
which describes when a hierarchy is {\em complete}.
We will also define a ``time order'', which is a partial order on a
hierarchy, generalizing the linear order on vertices of a single
geodesic, that will serve as a basic organizational principle in the
proofs here and in later sections.

\subsection{Definitions}\label{hierarchy defs}
\mbox{}

\bfheading{Tight geodesics.}
The non-uniqueness of geodesics in $\CC(S)$ is already manifested at a
local level, where typically, if $d_\CC(\alpha,\gamma) = 2$
there can be infinitely many choices for a curve $\beta$ disjoint
from both. The notion of
{\em tightness}, defined  below, addresses this local problem, but
more importantly introduces a crucial ingredient of control that makes
our combinatorial description of hierarchies possible. It is worth
noting that the only place where we make direct use of tightness is in
Lemma \ref{contiguous footprint}.

A pair of curves or curve systems $\alpha,\beta$ in a surface $Y$ are
said to {\em fill } $Y$ if all non-trivial non-peripheral curves in $Y$ 
intersect at least one of $\alpha$ or $\beta$. If $Y$ is a subdomain
of $S$ then it also holds that any curve $\gamma$ in $S$ which intersects a
boundary component of $Y$ must intersect one of $\alpha$ or $\beta$.

Given arbitrary curve systems $\alpha,\beta$ in $\CC(S)$, there
is a unique subsurface $F(\alpha,\beta)$ which they fill: Namely,
thicken the union of the geodesic representatives,  and fill in all
disks and once-punctured disks. Note that $F$ is connected if and only
if the union of geodesic representatives is connected. 

For a subdomain $X\subseteq Z$ let $\boundary_Z(X)$ denote the
{\em relative boundary} of $X$ in $Z$, i.e. those boundary components
of $X$ that are non-peripheral in $Z$.

\begin{definition}{tight seq def}
  Let $Y$ be a domain in $S$. If $\xi(Y)>4$,  a sequence of simplices
  $\{v_0,\ldots,v_N\}$ in $\CC(Y)$ is called {\em tight} if
\begin{enumerate}
\item For any vertices $w_i$ of $v_i$ and $w_j$ of $v_j$ where $i\ne
  j$, $d_{\CC(Y)}(w_i,w_j) = |i-j|$,
\item For each $1\le i \le N-1$, $v_i$ represents the relative
  boundary $\boundary_Y F(v_{i-1},v_{i+1})$.
\end{enumerate}

If $\xi(Y)=4$ then a tight sequence is just the vertex sequence of any
geodesic.

If $\xi(Y)=2$ then a tight sequence is the vertex sequence of any
geodesic, with the added condition that the set of endpoints on
$\boundary\hhat Y$ of arcs representing the vertices equals the set of
endpoints of the first and last arc.
\end{definition}
Note that condition (1) of the definition specifies that given any
choice of components $w_i$ of $v_i$ the sequence $\{w_i\}$ is a
geodesic in the original sense. It also implies that $v_{i-1}$ and
$v_{i+1}$ always have connected union.

In the annulus case, the restriction on endpoints of arcs is of little
importance, serving mainly to guarantee that there between any two
vertices there are only finitely many tight sequences.

With this in mind, a {\em tight geodesic} will
be a tight sequence together with some additional data:

\begin{definition}{tight geod def}
A {\em tight geodesic} $g$ in $\CC(Y)$ 
consists of a tight sequence
$\{v_0,\ldots,v_N\}$, and two markings $\I=\I(g)$ and $\T=\T(g)$
(in the sense of \S \ref{markings}), called its {\em
  initial} and {\em terminal} markings, such that
$v_0$ is
a vertex  of $\base(\I)$ and $v_N$ is a vertex of
$\base(\T)$.

The number $N$ 
is called the length of $g$, usually written $|g|=N$.
We refer to each of 
the $v_i$ as {\em vertices} of $g$ (by a slight abuse of notation).
$Y$ is called the {\em domain} or {\em support of $g$} and we write $Y=D(g)$.
We also say that $g$ is {\em supported in $D(g)$}.
\end{definition}

\medskip

Finally we will also, occasionally, allow tight geodesics to be
infinite, in one or both directions. If
a tight geodesic $g$ is infinite in the forward direction then $\T(g)$
is not defined, and if it is infinite in the backward direction then
$\I(g)$ is not defined.

\medskip
\bfheading{Subordinacy.}
We first saw the relations of 
{\em forward subordinacy} and {\em backward subordinacy} in the 
simple examples in Section \ref{example}. Let us now introduce a bit
more notation and give the general definitions.

\medskip
\noindent{\em Restrictions of markings:}
If $W$ is a domain in $S$ and $\mu$ is a marking in $S$, then 
the {\em restriction}  of $\mu$ to $W$, which we write
$\mu\rest W$, is constructed from $\mu$ in the following way:
Suppose first that
$\xi(W)\ge 4$. Recall that for every $p\in\mu$, either
$p=\alpha\in\base(\mu)$ or 
$p=(\alpha,t)$ with $t$ a transversal to $\alpha$. We let $\mu\rest W$
be the set of those $p$ whose base curve $\alpha$ meets $W$ essentially.
(Recall that $\alpha$ meets $W$ essentially if it cannot be deformed
away from $W$ -- in particular if $\alpha\subset W$ it must be
non-peripheral). 

If $W$ is an annulus ($\xi(W)=2$) then $\mu\rest W$ is just $\pi_W(\mu)$.

Note in particular that, 
if all the base curves of $\mu$ which meet $W$ essentially are
actually contained 
in $W$, then $\mu\rest W$ is in fact a marking of $W$.
If $W$ is an annulus then $\mu\rest W$ is a marking of $W$ whenever it
is non-empty.

\medskip
\noindent{\em Component domains:}
Given a surface $W$ with $\xi(W)\ge 4$ and a curve system $v$ in $W$ we say
that $Y$ is a 
{\em component domain of $(W,v)$} if either: $Y$ is a component of
$W\setminus v$, or $Y$ is an annulus with core a component of $v$.
Note that in the latter case $Y$ is non-peripheral, and thus satisfies
our definition of ``domain''.

Call a subsurface $Y\subset S$ a {\em component domain of $g$} if
for some vertex $v_j$ of $g$, $Y$ is a component domain of 
$(D(g), v_j)$. 
We note that this determines $v_j$ uniquely.
In such a case, let
$$\I(Y,g) = \left\{
\begin{array}{ll}
v_{j-1}\rest Y & v_j \text{ is not the first vertex} \\
\I(g)\rest Y & v_j \text{ is the first vertex} 
\end{array}
\right.
$$
be the {\em initial marking} of $Y$ relative to $g$. 
Similarly let

$$\T(Y,g) = \left\{
\begin{array}{ll}
v_{j+1}\rest Y & v_j \text{ is not the last vertex} \\
\T(g)\rest Y & v_j \text{ is the last vertex} 
\end{array}
\right.
$$
denote the {\em terminal marking}.
Note in particular that these are indeed markings.

\medskip
\noindent{\em Special cases:}
\begin{enumerate}
\item
The motivating case is that in which $v_j$ is neither first nor last,
and $\xi(D(g))>4$.
If $Y$ is the component of $D(g)\setminus v_j$ which is filled by
$v_{j-1}$ and $v_{j+1}$, then $\I(Y,g) = v_{j-1}$ and $\T(Y,g) =
v_{j+1}$.
If $Y$ is any other component domain of $(D(g),v_j)$ then
$\I(Y,g)=\T(Y,g)       = \emptyset$.
\item
If $Y$ is a thrice punctured sphere  ($\xi(Y)=3$) then
always $\I(Y,g) = \T(Y,g) = \emptyset$. 
\item
If $\xi(D(g))>4$ and $Y$ is an annulus (whose core curve is a component
of $v_j$), then unless $j=0$ or $j=|g|$, we must have
$\I(Y,g) = \T(Y,g) = \emptyset$, since successive curves in $g$ are
disjoint. If e.g. $j=0$, then the core of $Y$ is a base curve of
$\I(g)$, so if this curve has a transversal in the marking
$\I(g)$ then $\I(Y,g)$ is nonempty.
\item
If $\xi(D(g)) = 4$ then $Y$ must be an annulus, and now
$\I(Y,g)$ and $\T(Y,g)$ may be nonempty
because successive curves in $g$ do intersect.
\end{enumerate}


If $Y$ is a component domain of $g$ and $\T(Y,g)\ne\emptyset$ then we
say that $Y$ is {\em directly forward subordinate} to $g$, or $Y\fsubd g$.
Similarly if
$\I(Y,g)\ne\emptyset$ we say that $Y$ is {\em directly backward
subordinate} to $g$, or $g\bsubd Y$.

\medskip

We can now define subordinacy for geodesics:

\begin{definition}{subordinate def}
If $k$ and $g$ are tight geodesics,
we say that $k$ is {\em directly forward subordinate} to $g$,
or $k\fsubd g$, provided
$D(k)\fsubd g$ and $\T(k) =
\T(D(k),g)$.  
Similarly we define  $g\bsubd k$ to mean $g \bsubd D(k)$ and $\I(k) =
\I(D(k),g)$. 
\end{definition}

We denote by {\em forward-subordinate}, or $\fsub$, 
the transitive closure of $\fsubd$,
and similarly for $\bsub$.
We let $h\fsubeq k$  denote the condition that $h=k$ or
$h\fsub k$, and similarly for $k\bsubeq h$.
We include the notation $Y\fsub f$ where $Y$ is a domain
to mean $Y\fsubd f'$ for some $f'$ such that
$f'\fsubeq f$, and similarly define
$b\bsub Y$.

\bfheading{Hierarchies.}

\begin{definition}{hierarchy def}
A {\em hierarchy of geodesics} is a collection $H$ of
tight geodesics
in $S$ with the following properties:
\begin{enumerate}
\item
There is a distinguished {\em main geodesic} $g_H$ with domain $D(g_H)
= S$. The initial and terminal markings of $g_H$ are 
denoted also $\I(H), \T(H)$.

\item
Suppose $b,f\in H$, and  $Y\subset S$ is a domain such that  $b\bsubd
Y$ and $Y\fsubd f$. Then $H$ contains a 
unique tight geodesic $k$ such that $D(k)=Y$, $b\bsubd k$ and $k\fsubd f$.

\item
For every geodesic $k$ in $H$ other than $g_H$, there are $b,f\in H$
such that  $b\bsubd k \fsubd f$.
\end{enumerate}
\end{definition}

Condition (3) implies that for any $k$ in $H$, there is a sequence
$k=f_0\fsubd\ldots\fsubd f_m=g_H$, and similarly
$g_H=b_n\bsubd\ldots\bsubd b_0=k$.
Later we will prove these sequences are unique.  

\medskip\par\noindent{\em Infinite hierarchies.} 
An infinite hierarchy is one in  which the main geodesic $g_H$ is
allowed to be an infinite ray or a line. Note that in this case
$\I(H)$ and/or $\T(H)$ may not be defined.
Typically a hierarchy will be finite, but
most of the machinery of the paper will work for infinite hierarchies,
so we will indicate
where relevant how each proof works in the infinite case.
Infinite
hierarchies will arise, as limits, in 
\S\ref{limits of hierarchies}, and will be used in Section \ref{MCG}.

\subsection{Existence.}\label{existence results}
In this section we will prove that hierarchies exist. The first step
is the following:

\begin{lemma+}{Tight geodesics exist}
Let $u$ and $v$ be two vertices in $\CC(Y)$. There
exists a tight 
sequence $v_0,\ldots,v_N$ such that $v_0=u$ and $v_N=v$.
\end{lemma+}

(Note that whereas $u$ and $v$ are single vertices in the complex
$\CC(Y)$, the interior vertices of the sequence may actually be curve
systems, i.e. simplices of $\CC(Y)$.)
\begin{pf}
If $\xi(Y)= 4$ then the vertex sequence of
any geodesic is tight, by definition. 
If $\xi(Y)=2$ then the proof is an easy exercise. For example one can
start with $u$ and apply Dehn twists in the covering annulus $\hhat Y$
to obtain a sequence of curves with the same endpoints as $u$ on $\boundary
\hhat Y$, arriving at one which has one intersection with $v$ and making
one final step.

We now assume $\xi(Y)>4$.
To begin, let $h=\{u=u_0,\ldots,u_N=v\}$ be a regular geodesic connecting
$u$ and $v$. We will describe a process that adjusts $h$ until a tight
sequence is obtained. 

Let $v_1,v_2,v_3,v_4$ be any sequence of simplices
satisfying condition (1) of Definition \ref{tight seq def}, and suppose also
that $v_3$ is the boundary of $F(v_2,v_4)$, so that (2) holds for
$v_3$. If we now replace $v_2$ by 
$v'_2 = \boundary F(v_1,v_3)$, we want to show that $v_3$ is still
$\boundary F(v'_2,v_4)$. In other words, ``fixing'' $v_2$ so that
Condition (2) holds for it will not spoil condition (2) for $v_3$.

Note that (1) still holds for $v_1,v'_2,v_3,v_4$, by the triangle
inequality. In particular each 
component of $v'_2$ intersects each component of $v_4$, so that their
union is connected and so is $F(v'_2,v_4)$. Since $v'_2$ is disjoint
from $v_3=\boundary F(v_2,v_4)$, $v'_2$ must be contained in
$F(v_2,v_4)$ and in particular $F(v'_2,v_4) \subseteq
F(v_2,v_4)$. Thus it 
suffices to show that $v'_2$ and $v_4$ fill $F(v_2,v_4)$. Let $\alpha$
be any curve in $F(v_2,v_4)$. If $\alpha$ doesn't intersect $v_4$ then
it must intersect $v_2$, and also $v_1$ since $v_1$ and $v_4$ fill
$S$. But since $v_2$ is not contained in $F(v_1,v_3)$, $\alpha$ must
cross $\boundary F(v_1,v_3)$, which is just
$v'_2$.  We conclude that $F(v'_2,v_4) = F(v_2,v_4)$.

Now we can adjust the vertices of $h$ in any order: For any
$i\in[1,N-1]$ replace $u_i$ by $\boundary F(u_{i-1},u_{i+1})$.
For the new sequence, condition (2) holds for the $i$-th vertex. 
Repeating the process for a new value of $i$ in $ [1,N-1]$, the
previous argument assures us that the condition persists for
previously adjusted values of $i$. Thus after $N-1$ steps we obtain a
tight sequence.

Note that there is no reason to expect a unique tight sequence -- the
process seems to depend on the order in which the indices are chosen.
\end{pf}

We will now show, starting with any two markings in a surface $S$,
how to build a hierarchy connecting them. That is, 

\begin{theorem+}{Hierarchies exist}
Let $P$ and $Q$ be two markings in a surface $S$. There exists a
hierarchy $H$ of tight geodesics such that $\I(H)=P$ and $\T(H)=Q$.
\end{theorem+}

\begin{pf}
We say that $H$ is a {\em partial hierarchy} if it satisfies
properties (1) and (3) of Definition \ref{hierarchy def}, and the
uniqueness part of (2), but not 
necessarily existence. That is: 
\begin{enumerate}
\item[(2')]
Suppose $b,f\in H$, and  $Y\subset S$ is a domain such that  $b\bsubd
Y$ and $Y\fsubd f$. Then $H$ contains {\em at most one}
tight geodesic $k$ such that $D(k)=Y$, $b\bsubd k $ and $k\fsubd f$.

\end{enumerate}
Of course every hierarchy is also a partial hierarchy. 

We begin by choosing vertices $v\in \base(P)$ and $w\in\base(Q)$,
and connecting them with a tight sequence, which exists by Lemma
\ref{Tight geodesics exist}. Define a tight geodesic $g$ 
by letting its sequence be this one, and setting
$\I(g) = P$ and $\T(g)=Q$.

Let $H_0 $ be the partial hierarchy $\{g\}$, and let us construct a finite
sequence of partial hierarchies $H_n$,
the last of which is a hierarchy.

Call a triple $(Y,b,f)$ with domain $Y$ and $b,f\in H_n$ an {\em
  unutilized configuration} if
$b\bsubd Y \fsubd f$ 
but $Y$ is not the support of any geodesic $k\in H_n$ 
such that $b\bsubd k \fsubd f$.

Choose $(Y_n,b_n,f_n)$ to be any unutilized configuration in $H_n$.
Again use Lemma \ref{Tight geodesics exist} to construct a 
tight geodesic $h_n$ supported in $Y_n$, with $\I(h_n)=\I(Y_n,b_n)$ and
$\T(h_n)=\T(Y_n,f_n)$. Let $H_{n+1} = H_n \union \{h_n\}$.

The only thing to check is that the sequence terminates. 
Define a sequence of tuples $M_n = (M_n(1),M_n(2),\ldots,M_n(\xi(S)-2))$ by
letting $M_n(j)$ denote the number of unutilized configurations
$(Y_n,b_n,f_n)$ in $H_n$ with $\xi(Y_n) = \xi(S)-j$. 
Then, since for
each $Y_n$ in the above step, all component domains occurring in the
geodesic $h_n$ have complexity $\xi$ strictly smaller than $\xi(Y_n)$,
it follows immediately that the sequence $M_n$ is strictly decreasing in
lexicographic order as $n$ increases.
(Recall that in lexicographic order $(x_1,...,x_k)<(y_1,...,y_k)$
when for some $j\le k$, $x_i=y_i$ for all $i<j$ and $x_j<y_j$.)
Therefore the
sequence terminates in a partial hierarchy with no unutilized configurations
-- that is, a hierarchy.
\end{pf}

Note that the uniqueness part of property (2) holds automatically: although
the choice of $h_n$ at each stage was arbitrary, we
never put in more than one geodesic for a given configuration $b\bsubd
Y \fsubd f$.

\subsection{Forward and backward sequences.}
Given a domain $Y\subset S$ and a hierarchy $H$, define 
$$
\Sigma^+_H(Y) = \{ f\in H: Y\subseteq D(f) \ \ \text{and}\ \ 
\T(f)\rest Y\ne \emptyset\}
$$
and similarly
$$
\Sigma^-_H(Y) = \{ b\in H: Y\subseteq D(b) \ \ \text{and}\ \ 
\I(b)\rest Y\ne \emptyset\}
$$
which we also abbreviate by omitting the $H$ or $Y$ when they 
are understood.
For infinite hierarchies, we alter the definition by also admitting
$g_H$ into $\Sigma^+$ whenever $g_H$ is infinite in the forward direction,
and into $\Sigma^-$ whenever it is infinite in the backward direction.

We will call $\Sigma^+(Y)$ the {\em forward sequence of $Y$} and
$\Sigma^-(Y)$ the {\em backward sequence of $Y$}. These names will be
justified by the following theorem, 
which is perhaps the main point of our construction. 

\begin{theorem+}{Structure of Sigma}
Let $H$ be a hierarchy, and $Y$ any domain in its support $S$.

\begin{enumerate}
\item If $\Sigma^+_H(Y)$ is nonempty then it has the form of a
sequence 
$$f_0\fsubd\cdots\fsubd f_n=g_H,$$ 
where $n\ge 0$. Similarly, 
if $\Sigma^-_H(Y)$ is nonempty then it has the form of a
sequence 
$$g_H=b_m\bsubd\cdots\bsubd b_0,$$
where $m\ge 0$.

\item If 
$\Sigma^\pm(Y)$ are both nonempty, then $b_0 = f_0$, and
$Y$ intersects every vertex of $f_0$ nontrivially.
\item If $Y$ is a component domain in any geodesic $k\in H$ and
  $\xi(Y)\ne 3$, then
$$f\in \Sigma^+(Y) \ \ \iff \ \ Y\fsub f,$$
and similarly, 
$$b\in \Sigma^-(Y) \ \ \iff \ \ b\bsub Y.$$

If, furthermore, $\Sigma^\pm(Y)$ are both nonempty,
then in fact $Y$ is the support of $b_0=f_0$. 

\item Geodesics in $H$ are determined by their supports. That is, if
  $D(h)=D(h')$ for $h,h'\in H$ then $h=h'$.

\end{enumerate}
\end{theorem+}

The ingredients for the proof of 
Theorem \ref{Structure of Sigma} will be developed throughout the rest
of the section, and the proof will be completed in
\S\ref{structure sigma proof}.

In Section \ref{large link etc}, $\Sigma^\pm$ will be converted into
forward and backward {\em paths} in $\CC(S)$ which will enable us to
generalize the projection arguments in the examples of
\S\ref{example}, and prove Lemmas \ref{Large Link}, \ref{Common Links}
and their relatives.

\subsection{Footprints and subordinacy}
\label{footprints and subordinacy}
We begin with the following basic lemma, which gives one direction of
Part (3) of Lemma \ref{Structure of Sigma}:
\begin{lemma+}{Subordinate Intersection 1}
Let $H$ be a hierarchy in a surface $S$ and $Y$ a domain in
$S$. Let $h$ and $f$ denote geodesics in $H$.
\begin{enumerate}
\item
If $Y\fsubd h$ then $h\in\Sigma_H^+(Y)$. 
\item
If $h\in\Sigma_H^+(Y)$ and $h\fsub f$, then
$f\in\Sigma_H^+(Y)$.
\item
If $Y\fsub f$  then $f\in\Sigma^+_H(Y)$.
\end{enumerate}
The same holds with $\fsub$ replaced by $\bsub$, and $\Sigma^+$
replaced by  $\Sigma^-$.
\end{lemma+}

\bfheading{Footprints.}
We will first need one new definition, which will be a basic tool in
all that follows:
\begin{definition}{footprint def}
For a domain $Y\subset S$ and a
tight geodesic $g$ with non-annular support $D(g)\subset S$, let $\phi_g(Y)$
be the set of vertices of $g$ disjoint 
 from $Y$. We call this the {\em footprint} of $Y$ on $g$.
\end{definition}

If $Y$ is a subdomain of $D(g)$, then immediately
\begin{equation}\label{footprint diam bound}
\diam(\phi_g(Y)) \le 2
\end{equation}
in the curve complex of $D(g)$ (note if $Y$ is an annulus then by
definition of subdomain it is nonperipheral in $D(g)$).
It is also an immediate consequence of the definition that 
\begin{equation}\label{footprint containment}
Y\subseteq Z  \implies \phi_g(Z)\subseteq\phi_g(Y).
\end{equation}

Let us record the following elementary but crucial property of
footprints, which 
is the only place where the tightness
property is used directly.
\begin{lemma}{contiguous footprint}
If $g$ is a tight geodesic and  $Y\subset D(g)$ is a proper subdomain,
then $\phi_g(Y)$ is a sequence of 0,1, 2 or 3 contiguous vertices of
$g$.
\end{lemma}

\begin{pf}
When $\xi(D(g))=4$, $\phi_g(Y)$ is empty except when $Y$ is an annulus
whose core is some vertex $v$ of $g$. In that case every other vertex
intersects $Y$, so $\phi_g(Y)$ is the single vertex $v$.

Now assume $\xi(D(g))>4$.
The diameter bound (\ref{footprint diam bound}) implies that 
the only possibility for $\phi_g(Y)$ other than those mentioned in the
lemma is
that $\phi_g(Y)$ contains some $v_j$ and $v_{j+2}$ but not
$v_{j+1}$. However, since $g$ is a tight geodesic, 
if $Y$ intersects $v_{j+1}$ it either intersects $v_j$ or $v_{j+2}$,
since $v_{j+1}=\boundary_{D(g)}F(v_j,v_{j+2})$.
\end{pf}

Denote by $\min\phi_g(Y)$ and $\max\phi_g(Y)$ the 
vertices of $\phi_g(Y)$ with lowest and highest index, respectively.

\begin{pf*}{Proof of Lemma \ref{Subordinate Intersection  1}}
Clearly (3) is a consequence of (1) and (2), so we prove them.
We will prove the forward-subordinate case. The backward-subordinate
case proceeds similarly.

To see (1), suppose $Y\fsubd h$. Then by definition $Y$ is a component
domain of 
$(D(h),v_i)$ for some vertex $v_i$ of $h$, and $\T(Y,h)\ne\emptyset$.
If $v_i$ is the last vertex then $\T(Y,h) = \T(h)\rest Y$, so
this is nonempty and $h\in\Sigma^+(Y)$.

If $v_i$ is not the last vertex, we note that $v_i\in\phi_h(Y)$ and
$v_{i+1}$ is not in $\phi_h(Y)$. 
It follows, since the footprint is contiguous
(Lemma \ref{contiguous footprint}), that
the last vertex is not in $\phi_h(Y)$, hence 
$\T(h)\rest Y\ne\emptyset$, and again $h\in\Sigma^+(Y)$.
If $h$ is infinite in the forward direction ($h=g_H$, and $\T(h)$
undefined) then automatically $h\in\Sigma^+(Y)$. 

Now to prove (2), if $h\in\Sigma^+(Y)$ we have $Y\subset D(h)$ and
$ \T(h) \rest Y\ne \emptyset$.
Suppose first that $h\fsubd f$ -- then
$D(h)$ is a component domain of $(D(f),v_j)$ for
some vertex $v_j$ of $f$. If $v_j$ is the last vertex then $\T(h) =
 \T(f)\rest {D(h)}$, and it follows that 
$\T(f)\rest Y\ne\emptyset$.
If $v_j$ is not the last vertex then $\T(h) =  v_{j+1}\rest {D(h)}$
and since this intersects $Y$, $v_{j+1}\notin \phi_f(Y)$.
Thus also the last vertex of $f$ is not in $\phi_f(Y)$, and we may again
conclude that $ \T(f)\rest Y\ne\emptyset$ (or $f$ is infinite in the
forward direction). In each case we have $f\in\Sigma^+(Y)$. 
Part (2) now follows by induction.
\end{pf*}

The proof also gives the following slightly finer statement:

\begin{corollary+}{Footprints}
Let $H$ be a hierarchy, geodesics $h,f,b\in H$ and domain
$Y\subset D(h)$.
If $h\in\Sigma^+_H(Y)$ and $h \fsub f$, then 
$$\max \phi_f(D(h)) = \max \phi_f(Y).$$
Similarly if $h\in\Sigma^-_H(Y)$ and $b\bsub h$ then 
$$\min \phi_b(D(h)) = \min \phi_b(Y).$$
\end{corollary+}

Note, a special case of this is that if $h_1\fsub h_2 \fsub f$ then 
$\max\phi_f(D(h_1)) = \max\phi_f(D(h_2))$ by letting $Y=D(h_1)$.
(The condition $h_2\in\Sigma^+(D(h_1))$ follows from Lemma
\ref{Subordinate Intersection 1}.)

\begin{pf}
Since $h\fsub f$ there exists $h'$ such that $h\fsubeq h'\fsubd f$. 
Examining the proof of part (2) in 
Lemma \ref{Subordinate Intersection 1} above,  
we note that it shows that $D(h')$ is a component domain of $(D(f),v)$
where $v = \max\phi_f(D(h'))$, and if $v$ is not the last vertex in
$f$ then its successor intersects $Y$ and $D(h)$. Hence
$v=\max\phi_f(Y) =\max\phi_f(D(h))$. If $v$ is the last vertex then 
automatically $v=\max\phi_f(Y) =\max\phi_f(D(h))$, since
$\phi_f(D(h'))\subseteq\phi_f(D(h))\subseteq\phi_f(Y)$
by (\ref{footprint containment}).
The backward case proceeds similarly.
\end{pf}

\subsection{Uniqueness of descent}
\label{sigma basics}
By virtue of lemma \ref{Subordinate Intersection 1}, we know that
$\Sigma^+(Y)$ contains any sequence of geodesics $f_0,\ldots,f_n$ 
satisfying $f_0\in\Sigma^+(Y)$ and $f_i\fsubd f_{i+1}$
(and similarly for $\Sigma^-$).
The goal of the next lemma is to show that in fact $\Sigma^+$ and
$\Sigma^-$ are each just one such sequence, and as a consequence to prove
that geodesics in $H$ are determined by their domains.

\begin{lemma+}{Uniqueness of Descent}
Let $H$ be a hierarchy, and $Y$ any domain in its support $S$.

\begin{enumerate}
\item If $\Sigma^+_H(Y)$ is nonempty then it has the form of a
sequence $f_0\fsubd\cdots\fsubd f_n=g_H$, where $n\ge 0$.
Similarly 
if $\Sigma^-_H(Y)$ is nonempty then it has the form of a
sequence $g_H=b_m\bsubd\cdots\bsubd b_0$, where $m\ge 0$.

\item If there is some $h\in H$ with $D(h) = Y$, then there is exactly
one such $h$, and $h=f_0=b_0$.
\end{enumerate}
\end{lemma+}

In particular, this gives parts (1) and (4) of Theorem
\lref{Structure of Sigma}.

\begin{pf}
Note that (2) is a consequence of (1) for any given $Y$, since
if $Y=D(h)$ then $h\in \Sigma^+$ and must have the smallest domain of
any member of $\Sigma^+$ -- hence $h=f_0$, and similarly with
$\Sigma^-$ we have $h=b_0$. In particular $h$ is unique.

We will prove (1) 
by induction on $\xi(S)- \xi(Y)$. 
If $\xi(S)-\xi(Y)=0$ then $Y=S$ and $\Sigma^+=\Sigma^-=\{g_H\}$, hence
(1) holds.

Let $g\in\Sigma^+(Y)$, and suppose that $f\fsubd g$ for some
$f\in\Sigma^+(Y)$.
We claim that $D(f)$ is uniquely determined by $Y$ and $g$, and in
fact if $\xi(D(f))>\xi(Y)$ then $f$ itself is uniquely determined.
Suppose for a moment that $Y$ is not an annulus.
By definition of $f\fsubd g$, $D(f)$ is a component domain
for  $(D(g),v)$ where the vertex $v$ is $\max\phi_g(D(f))$.
By Corollary \lref{Footprints}, $v$ is also $\max\phi_g(Y)$.
Hence, $D(f) $ is the unique component of $(D(g),v)$ containing
$Y$, which depends only on $Y$ and $g$. 
Now if $\xi(D(f)) > \xi(Y)$ then by induction (2) holds for $D(f)$, so
that $f$ is the unique geodesic in $\Sigma^+$ such that $f\fsubd g$.

If $Y$ is an annulus the same proof goes through verbatim, recalling
that if $D(f)$ is not an annulus it must contain $Y$ as a
nonperipheral annulus (otherwise $\T(f)\rest{Y}$ would be empty,
contradicting $f\in\Sigma^+(Y)$),
and is therefore the unique component domain of
$(D(g),v)$ with this  property. If $D(f)$ is an annulus it must be
equal to $Y$ so again it is uniquely determined.

If  $\Sigma^+\ne\emptyset$, then since any $f\in
\Sigma^+$ is forward-subordinate to $g_H$, 
Lemma \ref{Subordinate Intersection 1} implies that
$g_H\in\Sigma^+$. For any $k\in \Sigma^+$ there exists
some $f$ such that $k\fsubeq f \fsubd g_H$, and $f\in \Sigma^+$
again by Lemma \ref{Subordinate Intersection 1}. By the previous
claim, we know that $D(f)$ is independent of $k$, and so is $f$ if
$D(f)\ne Y$.
Thus, replacing $g_H$ with $f$ and 
repeating this
argument inductively, we obtain a single sequence
$f_1\fsubd\cdots \fsubd g_H$ which accounts for all of $\Sigma^+$
except possibly those geodesics $f$ with $D(f)=Y$.

Now repeating this for $\Sigma^-$ we obtain a sequence
$g_H\bsubd\cdots\bsubd b_1$. If $Y$ does not support any geodesic then
we are done (reindexing both sequences to start with $0$).
If $Y$ supports at least one geodesic $h$, then
$h\in\Sigma^+\intersect \Sigma^-$, and by the same logic as above we
have that $h\fsubd f_1$ and $b_1\bsubd h$. However, by the uniqueness
part of the definition of a hierarchy, there can only be one
such $h$. Setting $f_0=b_0=h$, we are done. 
\end{pf}

The following partial converse of Lemma
\ref{Subordinate Intersection 1} is an immediate corollary of Lemma
\ref{Uniqueness of Descent}:
\begin{lemma+}{Subordinate Intersection 2}
Let $k$ and $h$ be geodesics in a hierarchy $H$.

If $h\in \Sigma^+_H(D(k))$ then
$k\fsubeq h$. Similarly,

If $h\in \Sigma^-_H(D(k))$ then
$h\bsubeq k$. 
\end{lemma+}

Here is an easy corollary of Lemma \ref{Subordinate Intersection 2}. 
\begin{corollary}{easy containment}
If $D(k)$ is properly contained in $D(h)$ then $\phi_h(D(k))$ is non-empty.
\end{corollary}
\begin{pf}
If $h$ is infinite then $h=g_H$ and by definition $h\bsub k$, so
$\phi_h(D(k))\ne \emptyset$. Otherwise $\I(h)$ is defined.
If $\I(h)\rest{D(k)}$ is empty then $\phi_h(D(k))$ contains the
initial vertex. If not, then $h\in\Sigma^-(D(k))$ and,
by Lemma \ref{Subordinate Intersection 2}, $k$ is
backward-subordinate to $h$, and hence $D(k)$ is contained in a
component domain 
for some other vertex, so that again $\phi_h(D(k))$ is non-empty.
\end{pf}

Let us also record the following consequence of these lemmas:
\begin{lemma}{Y direct to top}
Let $Y$ be a domain in $S$ and $h$ in a hierarchy $H$ such that
$Y\fsubd h$. Then $h$ is uniquely determined, and in particular,
writing $\Sigma^+_H(Y) =
\{f_0,\ldots\}$ we have 
either $h=f_1$ and $Y = D(f_0)$, or  $h=f_0$ and $Y$ supports no
geodesic in $H$. 

The corresponding statement holds when $h\bsubd Y$.
\end{lemma}
\begin{pf}
By Lemma \ref{Subordinate Intersection 1}, $h\in\Sigma^+(Y)$. 
If $h=f_0$ then $Y$ cannot support a geodesic because $D(h)$ has the
smallest domain among elements of $\Sigma^+(Y)$.
Suppose $h=f_{i+1}$ where $i\ge 0$. Then by Corollary
\ref{Footprints},
$\max\phi_h(Y) = \max\phi_h(D(f_i))$ and hence both $Y$ and $D(f_i)$
are component domains for the same vertex of $h$. As in the proof of
Lemma \ref{Uniqueness of Descent} we conclude $Y=D(f_i)$ and so $i=0$
since there can be no smaller domain in $\Sigma^+(Y)$.
The case where $h\bsubd Y$ is similar.
\end{pf}

\subsection{Time order}
\label{time order section}
The vertices of any geodesic admit a linear order from initial to
terminal, and the relations $\fsub$ and $\bsub$ are, by definition,
partial orders. It turns out that these can be combined to define a
useful partial order $\tprec$ on a hierarchy, and a
related partial order $\pprec$ on the set of ``pointed geodesics'' of
a hierarchy. In this section we define these relations and study their
basic properties. Let a hierarchy $H$ be fixed througout this section.

\begin{definition}{time order def}
For any $h,h'\in H$, 
we say that {\em $h$ precedes $h'$ in time order}, or
$$h\tprec h',$$ 
if there exists a geodesic $m\in H$
such that $D(h), D(h')\subset D(m)$, and
$$\max\phi_m(D(h)) < \min\phi_m(D(h')).$$
(In particular $\phi_m(D(h))$ and
$\phi_m(D(h'))$ are disjoint.)
\end{definition}

Note that if this occurs then automatically $h\fsub m$ and $m \bsub h'$:
since $\phi_m(D(h))$
must miss the terminal vertex of $m$,  $m$ is in $\Sigma^+(D(h))$ and we
may apply Lemma \ref{Subordinate Intersection 2}, and similarly for
$h'$ using $\Sigma^-(D(h'))$.

We call $m$ the {\em geodesic used to compare $h$ and $h'$}, and note that
it is unique: If some $m'$ is also used to obtain $h\tprec h'$, then
both $m$ and $m'$ appear in the forward sequence of $h$ and the
backward sequence of $h'$. In particular either $m\fsub m'$ or
$m'\fsub m$; suppose the first, without loss of generality. Then
$\phi_{m'}(D(m))$ is non-empty, and by (\ref{footprint containment})
is contained in both $\phi_{m'}(D(h))$ and $\phi_{m'}(D(h'))$,
contradicting the assumption that they are disjoint.

If either $h\tprec h'$ or $h'\tprec h$ then we say $h$ and $h'$ are {\em
time-ordered}. Note, we have not yet shown that these two
possibilities are mutually exclusive,
or indeed that $\tprec$ is a partial order. Before we do
that let us define a more general relation.

\bfheading{Partial order on pointed geodesics.}
Let $k$ be a tight geodesic with
vertices $v_0,\ldots,v_N$. We generalize slightly the notion of vertex
to a {\em position} on $k$, which is either a vertex
$v_i$, or $\I(k)$ or $\T(k)$. The linear order $v_i<v_j$ when $i<j$
extends to an order on positions where we say $\I(k) < v_0$ if the two
are not the same,
and similarly $v_N<\T(k)$ if the two are not the same.
We can now discuss {\em pointed geodesics}, which are
pairs $(k,v)$ where $v$ is a position in $k$.

We extend the notion of footprint slightly as follows:
Given a pointed geodesic $(k,v)$ and a geodesic $h$ with
$D(k)\subseteq D(h)$, we define
\begin{equation*}
\hat\phi_h(k,v) = \begin{cases}
                \phi_h(D(k))    & \text{if}\ \ D(k) \subset D(h), \\
                \{v\}               & \text{if}\ \ k=h.
              \end{cases}
\end{equation*}
Note that $\hat\phi_h(k,v)$ could be $\{\I(k)\}$ or $\{\T(k)\}$ in the second
case, in contrast with regular footprints which can only consist of
vertices. 
We now define a relation $\pprec$ on pairs $(k,v)$:
\begin{definition}{pair order def}
We say
$$(k,v) \pprec (h,w)$$ 
if and only if there exists a geodesic $m$ such that
$k\fsubeq m \bsubeq h$, and 
$$\max\hat\phi_m(k,v) < \min\hat\phi_m(h,w).$$
\end{definition}

Again, it is clear that $m$ is unique. It is also immediate from the
definitions that
\begin{equation}
\label{pprec defines tprec}
k \tprec h \iff (k,\T(k)) \pprec (h,\I(h)).
\end{equation}
Indeed, $(k,v)\pprec(h,w)$ breaks up into four 
possible, mutually exclusive, cases:
\begin{itemize}
\item $k\tprec h$,
\item $k=h$ and $v<w$,
\item $k\fsub h$ and $\max \phi_h(D(k)) < w$
\item $k\bsub h$ and $v < \min \phi_k(D(h))$.
\end{itemize}

We now verify that these relations are partial orders, together with a
number of other properties. The following lemma holds for infinite as
well as finite hierarchies.

\begin{lemma+}{Time Order}
\mbox{}

\begin{enumerate}
\item
If $D(h)\subseteq D(h')$ then $h$ and $h'$ are not time-ordered.
\item
On the other hand if $D(h)\intersect D(h') \ne \emptyset$ and neither
domain is contained in the other, then $h$ and $h'$ are time-ordered.
\item
Suppose $b\bsubeq k\fsubeq f$.
Then either $b=f$, $b\fsub f$, 
$b\bsub f$, or $b\tprec f$. 
\item
Suppose that  $k_1\tprec k_2$.  If $h\fsubeq k_1$ then 
$h\tprec k_2$. Similarly if $k_2\bsubeq g$ then $k_1\tprec g$.
\item
The relation $\pprec$ is a (strict) partial order.
\item
The relation $\tprec$ is a (strict) partial order.
\end{enumerate}
\end{lemma+}

\begin{pf}
To prove part (1),
suppose $D(h)\subseteq D(h')$. Then in any geodesic $m$ such that $D(m)$
contains $D(h')$, 
$\phi_m(D(h'))\subseteq \phi_m(D(h))$. In particular the footprints can
never be disjoint, and hence neither $h\tprec h'$
nor $h'\tprec h$ can hold.

\medskip

Next let us prove part (2). Suppose $D(h)\intersect D(h')\ne \emptyset$.
Consider the following assertion:  
If $m$ is a geodesic such that $D(m)$ contains $D(h)\union D(h')$, and
also  $m\bsub h$, then either
$D(h')\subseteq D(h)$, or $D(h)\subseteq D(h')$, or $h$ and $h'$ are
time-ordered. 
We shall prove this by induction on
$\xi(D(m))$.

If $\xi(D(m))=2$, 
then $m=h=h'$ and we are done. More generally, 
if $\xi(D(m)) = \xi(D(h'))$ then $D(h)\subseteq D(h')$, and again we
are done.
Otherwise, we must have
$\xi(D(m)) > \xi(D(h'))$, so consider the footprints $\phi_m(D(h))$
and $\phi_m(D(h'))$ (the former is non-empty since $m\bsub h$, and the
latter is non-empty by Corollary \ref{easy containment}). 
If they are disjoint then $h$ and $h'$ are time-ordered and
again we are done. If
they overlap then, since each is an interval of  contiguous
vertices, the minimum of one must be contained in the other.

If the minimum $v$ of $\phi_m(D(h))$ is contained in $\phi_m(D(h'))$, let
$m'$ be the geodesic such that $m\bsubd m' \bsubeq h$.
Then $D(m')$ is a 
component domain of $(D(m),v)$. Since $v\in \phi_m(D(h'))$, 
$v$ does not intersect $D(h')$,
and since $D(h)$ and $D(h')$ intersect, $D(h')$
must be in the same component domain of $(D(m),v)$, namely $D(m')$.
If $h=m'$ then we are done, with $D(h')\subset D(h)$.  If not, we may
apply the inductive assumption since $\xi(D(m'))<\xi(D(m))$, and again
we are done.

If $v$ is not in $\phi_m(D(h'))$, then  $v'= \min\phi_m(D(h'))$ 
must be in $\phi_m(D(h))$, and furthermore $v'$ is {\em
not} the first vertex of $m$ since $v$ lies to its left. Thus
$m\in\Sigma^-(D(h'))$ and it follows by
Lemma \ref{Subordinate Intersection 2} that $m\bsub h'$,
and therefore we may reverse the roles of $h$ and $h'$ and apply the
previous paragraph. This concludes the proof of the assertion, and (2)
follows by applying the assertion when $m$ is the main geodesic $g_H$.

\medskip

To prove (3), we may 
suppose that $b\bsub k\fsub f$ and $b\ne f$, as the cases of equality
are trivial.

By Lemma \ref{Subordinate Intersection 1}, $f\in\Sigma^+(D(k))$.
If $D(b)\subset D(f)$ then 
$f\in\Sigma^+(D(b))$ as well, and by Lemma
\ref{Subordinate Intersection 2}, $b\fsub f$. Similarly if
$D(f)\subset D(b)$, then $b\bsub f$.

Now suppose that neither $D(b)\subseteq D(f)$ nor
$D(f)\subseteq D(b)$. 

Since $D(k) \subset D(b)\intersect D(f)$, part (2) shows that
$b$ and $f$ are time-ordered, so let $D(m)$ contain both
$D(b)$ and $D(f)$ such that their footprints on $m$ are disjoint. 

We claim that $\min\phi_m(D(b))=\min \phi_m(D(k))$. If $b$ is
backward-subordinate to $m$ then 
this is just Corollary \lref{Footprints}. If $b$ is not backward subordinate
to $m$ then neither is $k$, so that by Lemma 
\ref{Subordinate Intersection 2} both 
minima are equal to the first vertex of $m$.

Similarly, $\max\phi_m(D(f))=\max\phi_m(D(k))$. It follows that
$\max\phi_m(D(b))< \min\phi_m(D(f))$ 
(since we already know they are
disjoint). Thus $b\tprec f$, as desired.

\medskip

Next we prove (4). Since $k_1\tprec k_2$, let $m$ be the geodesic used
to compare them. Since $h\fsubeq k_1\fsub m$,
$\max\phi_m(D(h)) = \max\phi_m(D(k_1))$ by Corollary \lref{Footprints}.
Thus $\max\phi_m(D(h))< \min\phi_m(D(k_2)$, so $h\tprec k_2$. The proof 
that  $k_1\tprec g$ is similar.

\medskip

To prove (5),  we must in particular show that $\pprec$
is transitive. Suppose $(k,v) \pprec (k',v')$, and 
$(k',v') \pprec (k'',v'')$. 
Let $b$ be the geodesic used to compare $(k,v)$ and $(k',v')$,  and $f$ be the
geodesic used to compare $(k',v')$ and $(k'',v'')$. Then in particular
$k\fsubeq b \bsubeq k' \fsubeq f \bsubeq k''$.

By (3), either $b=f$, $b\fsub f$, 
$b\bsub f$, or $b\tprec f$. 

If $b=f$ then the footprints $\hat\phi_b(k,v)$,
$\hat\phi_b(k',v')$, and
$\hat\phi_b(k'',v'')$ are disjoint and linearly ordered from left to
right, so $(k,v) \pprec (k'',v'')$ immediately.

If $b\fsub f$ then neither $k$ nor $k'$ can equal $f$. Thus
$\hat\phi_f(k,v) = \phi_f(D(k))$ and 
$\hat\phi_f(k',v') = \phi_f(D(k'))$. By Corollary \lref{Footprints} we
have
$\max\phi_f(D(k)) = \max\phi_f(D(b))$, and by 
(\ref{footprint containment}) we have  
$\phi_f(D(b)) \subset \phi_f(D(k'))$. Then since $\max\phi_f(D(k')) <
\min \hat\phi_f(k'',v'')$, we conclude $\max\hat\phi_f(k,v) < \min
\hat\phi_f(k'',v'')$ so that $(k,v) \pprec (k'',v'')$.
The case where $b\bsub f$ follows similarly.

If $b\tprec f$ then since $k\fsubeq b$
part (4) gives $k\tprec f$, and since $f\bsubeq k''$
part (4) gives $k\tprec k''$. It follows that $(k,v) \pprec (k'',v'')$.

We have therefore proved that $\pprec$ is transitive. 
It follows from the definition that $(h,v)\pprec (h,v)$ can never
hold, and so transitivity implies $(h,v)\pprec(k,w)$ and
$(k,w)\pprec(h,v)$ are mutually exclusive.
Thus $\pprec$ is a strict partial order.

\medskip

Part (6) follows immediately from the relation (\ref{pprec defines
tprec}) between $\pprec$ and $\tprec$. It is also easy to see it
directly by an argument very similar to the above.
\end{pf}

The next lemma gives a sufficient condition for two
geodesics with disjoint domains {\em not} to be time-ordered.
\begin{lemma}{diagonalizable}
Let $v$ be a vertex of $h$ and suppose that
$D(k)$ and $D(k')$ lie in different component domains of $(D(h),v)$.
Then $k$ and $k'$ are not time-ordered.
\end{lemma}

Remark: we expect that there {\em will} be geodesics
with disjoint domains which are nevertheless time-ordered. This is in
fact one of the serious difficulties in applications of hierarchies.

\begin{pf}
Suppose by way of contradiction that
$k\tprec k'$, and let $m$ be the geodesic
used to compare them. 

Note first  $m\ne h$ since the footprints of $k$ and
$k'$ on $m$ are disjoint, and on $h$ they both contain $v$. 
 
Suppose $D(m)\subset D(h)$. Then $\min\phi_h(D(m))=\min\phi_h(D(k'))$
since $m\bsub k'$, just as we argued in the proof of Lemma
\ref{Time Order}, part (3). Similarly $\max\phi_h(D(m))=\max\phi_h(D(k))$. It
follows that $\phi_h(D(m))$ 
contains $\phi_h(D(k))\intersect \phi_h(D(k'))$, which in particular
contains $v$. Thus $D(m)$ is disjoint from $v$, and since  $D(m)$ is
connected, it must lie in one component domain of $(D(h),v)$. This
contradicts the assumption that $D(k)$ and $D(k')$ lie in different
components.

Now suppose $D(h)\subset D(m)$. By Corollary \ref{easy containment},
$\phi_m(D(h))$ is nonempty, and since $D(h)$ contains both $D(k)$ and
$D(k')$, $\phi_m(D(h))\subset \phi_m(D(k))\intersect \phi_m(D(k'))$.
This contradicts the disjointness of footprints of $k$ and $k'$  in
$m$. 

Finally if neither $D(h)$ nor $D(m)$ is contained in the other, their
intersection is still non-empty since both contain $D(k)\union D(k')$.
Thus by Lemma \ref{Time Order} part (2), $h$ and $m$ are time-ordered.
Suppose $h \tprec m$. Since $m\bsub k'$, by
Lemma \ref{Time Order} part (4) we have $h\tprec k'$. However
$D(k')\subset D(h)$ so this contradicts Lemma \ref{Time Order} part
(1). Similarly if $m\tprec h$ then since $k\fsub m$
we have $k\tprec h$, again a contradiction.
\end{pf}

\subsection{Complete hierarchies}
\label{complete hierarchy section}
A hierarchy $H$ is {\em complete} if, for every domain $Y$ with $\xi(Y)\ne
3$, which is a
component domain in some geodesic $k\in H$, there is a geodesic $h\in
H$ with $Y=D(h)$.
In this section we will prove:

\begin{theorem+}{Completeness}
If the markings $\I(H)$ and $\T(H)$  (where defined) are complete, 
then $H$ is complete.
\end{theorem+}

This will require Lemma \ref{Subordinate Intersection 3} below,
which is also the last and perhaps trickiest piece needed to prove
Theorem \lref{Structure of Sigma}.
This lemma addresses the issue of when a component domain appearing in
a hierarchy is the support of a geodesic in the hierarchy.

If $v$ is a marking in a non-annular $W$, a {\em component domain of
$(W,v)$}  is defined to be any component domain of
$(W,\base(v))$. This slight generalization
will be used below for component domains defined by {\em
  positions} in geodesics, including the initial or terminal markings.

\begin{lemma+}{Subordinate Intersection 3}
Let $H$ be a hierarchy, and let
$Y$ be a component domain of $(D(k), v)$, where $k\in H$ and
$v$ is a position in $k$. Assume that $\xi(Y) \ne 3$.

If $f\in \Sigma^+(Y)$ then either $D(f)= Y$ or $Y\fsub f$.

Similarly if $b\in \Sigma^-(Y)$ then either $D(b) = Y$ or $b\bsub Y$.
\end{lemma+}

Note the similarity of this to Lemma
\lref{Subordinate Intersection 2}, the main difference being that $Y$
is not required to be the support of a geodesic. In fact, 
the conclusions $Y\fsub f$ and $b\bsub Y$ together imply that $Y$ is
the support of a geodesic, so in particular the
lemma gives a sufficient condition for this to occur.

\begin{pf}
We will prove the forward case. Note that it suffices to show that, 
if $\Sigma^+_H(Y) \ne \emptyset$ then there exists $h\in H$ such that
$Y\fsubd h$. Lemma \ref{Y direct to top} then implies that, 
writing $\Sigma^+(Y) = \{f_0,f_1,\ldots\}$, 
either $h=f_0$ and $Y$ is not the support
of any domain, or $Y=D(f_0)$ and $h=f_1$.
The lemma follows from this.

We argue by induction using the partial order $\pprec$.
In particular,
let us write our data as $(Y,k,v)$ where $v$ is a position on $k$ such that
$Y$ is a component domain of $(D(k),v)$.
Let $N(Y,k,v)$ be the number of such triples $(Y',k',v')$
in $H$ for which $Y\subseteq Y'$ and $(k,v) \pprec (k',v')$. 
Note that this number is finite even if the hierarchy is infinite,
because $\phi_{g_H}(D(k'))\subset \phi_{g_H}(Y)$, so the candidate triples
are limited to a finite subset of $H$.
Clearly
if $Y\subseteq Y'$ and $(k,v)\pprec (k',v')$, then $N(Y,k,v)>N(Y',k',v')$,
so we will induct on $N$.

Let us first show that either $Y\fsubd k$, in which case we are done,
or we can find a certain $(Y',k',v')$ with $Y\subseteq Y'$ and $(k,v)\pprec
(k',v')$ to which we can apply the inductive hypothesis. In particular
this will take care of the base case $N(Y,k,v)=0$. We will then use $Y'$
to deduce the desired conclusion for $Y$.

For a position $w < \T(k)$ in $k$ let $succ(w)$ denote the next
position in the linear order. The following two cases occur:

\bfheading{1:} If $v<\T(k)$, let $v'=succ(v)$. 
If $v'\rest Y$ is nonempty then $\T(Y,k)$ is nonempty and
$Y\fsubd k$, as desired.
If not, we define
$(Y',k',v')$ by letting
$k'=k$ and letting $Y'$ be the component domain of
$(D(k),v')$ containing $Y$. 

(Remarks: If $v=\I(k)$ or $v'=\T(k)$ then $v$ and $v'$ will share some
base curves. In this case it is possible that $Y=Y'$. 
If $Y$ is an annulus then
either $Y'=Y$ or $Y'$ contains $Y$ essentially.)

\bfheading{2:}
If $v=\T(k)$ (including the case $\T(k)$ is equal to the last vertex),
we let $k'$ be
the geodesic such that $k\fsubd k'$.
Then $D(k)$ is a component domain of some vertex $w$ in $k'$, with $w<\T(k')$.
Letting $w'=succ(w)$, we have 
$\T(k)  = w'\rest {D(k)}$.
Let $v'=w'$ and let $Y'$ be the component domain of $(D(k'),w')$
containing $Y$.

\medskip

Note that in each case, $(k,v) \pprec (k',v')$: When $k=k'$ we had
$v<v'$, and when $k\fsubd k'$ we had $\phi_{k'}(D(k)) < v'$.

Recall, the assumption that $\Sigma^+(Y) \ne\emptyset$ is equivalent by Lemma
\ref{Uniqueness of Descent} to  $ \T(g_H)\rest Y\ne \emptyset$
or $g_H$ infinite in the forward direction.
Thus, since $Y\subseteq Y'$, $\Sigma^+(Y)\ne\emptyset$
implies $\Sigma^+(Y')\ne\emptyset$. Now we can apply the inductive
assumption to  $(Y',k',v')$ to conclude that
$Y'\fsubd h'$ for some $h'$.

Thus  if $Y=Y'$, we are done. 

From now on let us assume $Y$ is a proper subdomain of $Y'$.
Thus, the relative boundary $\boundary_{Y'}(Y)$ is nonempty.
We claim that $\boundary_{Y'}(Y)$ consists of (base)
curves of $v$ in case 1, and of $w$ in case 2:

In case 1, $\boundary_{D(k)}(Y)$ is in $\base(v)$, and so the claim is
immediate. In case 2, $\boundary_{D(k)}(Y)$ is in $\base(\T(k))$ which is
part of $w'$, and $\boundary_{D(k')}(D(k))$ is in $w$. It follows that
$\boundary_{D(k')}(Y)$ is in $w\union w'$. Since no curve of $w'$
is an essential curve of $Y'$, $\boundary_{Y'}(Y)$ must be in $w$.

Noting also that $\I(Y',k') = v\rest {Y'}$ in case 1
and $\I(Y',k') = w\rest {Y'}$ in case 2, we conclude in each case that
$\I(Y',k')\ne  \emptyset$. Thus,
$k'\bsubd Y'$. Since we already have $Y'\fsubd h'$,
by definition of a hierarchy there is
a geodesic $h\in H$ whose support is $Y'$, and
$k'\bsubd h$. In fact $(k,v) \pprec
(h,\I(h))$: when $k=k'$ this is because the footprint of $h$ in $k$
has minimum at $v'$. When $k\fsubd k'$, it is because $\max\phi_{k'}(D(k))
= w < w' = \min\phi_{k'}(D(h))$, so in fact $k\tprec h$.
Also, $Y$ is a component domain of $(D(h),\I(h))$:
In case 1 this is because $\I(h) = v\rest {Y'}$ and
$\boundary_{Y'}(Y)\subseteq v$, and in case 2 it is because $\I(h) =
w\rest {Y'}$ and $\boundary_{Y'}(Y)\subseteq w$.
Thus, the lemma holds for $(Y,h,\I(h))$ by induction, and we are done.
\end{pf}

\begin{proof}[Proof of Completeness Theorem]
Let $Y$ be any component domain that is not a thrice-punctured sphere.
If $\I(H)$ is defined then it is complete, so $\I(H)\rest Y$ is
nonempty and hence $g_H\in\Sigma^-(Y)$. If $\I(H)$ is undefined then
$g_H$ is infinite in the backward direction so  $g_H\in\Sigma^-(Y)$ by
definition. Similarly $g_H\in\Sigma^+(Y)$, so that by
Lemma \ref{Subordinate Intersection 3} 
$g_H\bsub Y \fsub g_H$. In particular $Y$ must support a geodesic by
definition of a hierarchy.
\end{proof}

Note that, even if $\I(H)$ and $\T(H)$ are pants decompositions with
no transverals at all, we get a hierarchy which is complete except
for the annular domains whose cores are curves of $\I(H)$ or $\T(H)$.
If $H$ has a bi-infinite main geodesic then it is automatically
complete with no further conditions.

\subsection{Proof of the Structural Theorem}
\label{structure sigma proof}
We now have all the ingredients in place to put together a proof of
Theorem \lref{Structure of Sigma}. Parts (1) and (4) were already
shown in Lemma \lref{Uniqueness of Descent}. 
For Part (3), one direction of 
$$
f\in\Sigma^+(Y) \ \iff \ Y\fsub f
$$
follows from
Lemma \lref{Subordinate Intersection 1}, and the other from
Lemma \lref{Subordinate Intersection 3}; similarly for $\Sigma^-$.
If $\Sigma^\pm(Y)$ are both nonempty then, again by Lemma
\ref{Subordinate Intersection 3}, there must be $b$ and $f$
such that $b\bsubd Y \fsubd f$ and so by definition of a hierarchy $Y$
is the support of a geodesic, which must then be $b_0=f_0$.

It remains to show Part (2), that if
$\Sigma^\pm(Y)$ are both nonempty for any domain $Y$,
then $b_0=f_0$, and $\phi_{f_0}(Y) = \emptyset$.
If $Y=D(h)$ for
some geodesic $h$ then we already know $h=b_0=f_0$,
and then $\phi_{f_0}(Y) = \emptyset$ automatically.

In general, if $\phi_{f_0}(Y)\ne
\emptyset$, let $v=\max\phi_{f_0}(Y)$. Let 
$W$ be the component domain of $(D(f_0),v)$ containing $Y$. Recall that
Lemma \ref{Uniqueness of Descent} implies for any domain $X$ that
$\Sigma^+(X)\ne \emptyset$ is equivalent to either $\T(H)\rest X\ne\emptyset$
or  $g_H$ infinite in the forward direction, and
similarly for $\Sigma^-$ and $\I(H)$. Thus
$\Sigma^\pm(Y)\ne \emptyset$ implies $\Sigma^\pm(W)\ne \emptyset$.
Since $W$ is also a component domain, Lemma \ref{Subordinate
  Intersection 3} then implies that $W$ supports some geodesic
$h$. Letting $v'$ be the  successor of $v$, we note that $Y$ intersects
$v'$ nontrivially (this is true even when $v$ is the last vertex of
$f_0$ and $v'=\T(f_0)$). But $ v'\rest W$ is just $\T(h)$ by
definition, and so $h\in\Sigma^+(Y)$. This contradicts the fact that
$f_0$ has the smallest domain in $\Sigma^+$.
Thus $\phi_{f_0}(Y)=\emptyset$.

In particular, it follows that  $f_0\in \Sigma^-$ since 
$\I(f_0)\rest Y\ne \emptyset$ (or $f_0$ is infinite in the backward
direction). In fact $f_0$ must
be $b_0$ since $Y$ has nonempty footprint on any other geodesic in
$\Sigma^-$.  This concludes the proof.

%% file: slices.tex
\section{Slices, resolutions, and markings}
\label{resolution}

In this section we will discuss how to resolve a hierarchy $H$ into a
sequence of markings connecting $\I(H)$
to $\T(H)$, so that successive markings are related by 
elementary moves. Essentially we must somehow combine the vertex
sequences of the various geodesics in $H$, and their partial orders,
into one large linearly ordered sequence. This process is by no means
unique.

Along the way we will need to develop the notion of a {\em slice},
which roughly speaking is a marking pieced together from variously
nested geodesics in the hierarchy, together with additional
organizational structure. These slices will admit a certain partial
order, and we will then describe an elementary move on slices, which
moves a slice forward in the partial order. 

The resulting sequence of slices can then be transformed into a
sequence of clean markings of the surface (in a slightly non-unique
fashion), and we will prove a lemma bounding the length of this
sequence in terms of the size of $H$. 

\bfheading{Slices.}
Let us assume from now on that the hierarchy $H$ is complete.
A {\em slice} in $H$ is a set $\tau$ of pairs  $(h,v)$ where
$h\in H$ and $v$ is a vertex of
$h$, satisfying the following conditions:
\begin{itemize}
\item[S1:] A geodesic $h$ appears in at most one pair in $\tau$.
\item[S2:] There is a distinguished pair $(h_\tau,v_\tau)$ appearing in $\tau$,
  called the bottom pair of $\tau$. We call $h_\tau$ the bottom geodesic. 
\item[S3:] For every $(k,w)\in \tau$ other than the bottom pair, $D(k) $
  is a component domain of $(D(h),v)$ for some $(h,v)\in\tau$.
\end{itemize}

If in addition this fourth condition holds, we call the slice {\em
complete}:
\begin{itemize}
\item[S4:] Given $(h,v)\in\tau$, for every component domain $Y$ of
  $(D(h),v)$ there is a pair $(k,w)\in\tau$ with $D(k)=Y$.
\end{itemize}

Most often $h_\tau$ will just be the main geodesic of $H$.

A slice $\tau$ is called {\em initial} if, for each $(h,v)\in\tau$,
$v$ is the first vertex of $h$. Note that a complete initial slice is
uniquely determined by its bottom geodesic.
The complete initial slice with bottom geodesic $g_H$ is called
{\em the} initial slice of $H$. We similarly
define terminal slices.

\bfheading{Markings associated to a slice:}
To any slice $\tau$ we associate a unique marking $\mu_\tau$ as follows.
It is easy to see by induction that the vertices $v$ appearing in
non-annular geodesics in 
$\tau$ are all disjoint and distinct, and hence form a simplex in $\CC(S)$. 
We let this be $\base(\mu_\tau)$.
For each 
base curve $v$, if $\tau$ contains some $(k,t)$ with $D(k)$ the annulus whose
core is $v$, then we let $t$ be the transversal of $v$ in $\mu_\tau$.
In particular a complete slice determines a complete marking.

Typically $\mu_\tau$ is not clean, so
let us say that a clean marking $\mu'$ is {\em compatible with $\tau$}
if it is compatible with $\mu_\tau$ in the sense of Lemma \ref{clean
  markings}. 
Lemma \ref{clean markings} then shows that such a $\mu'$ exists,
there are at most $\ncompat^b$ possibilities where $b=\#\base(\mu)$,
and any two  differ by  a bounded number of Twist elementary 
moves.

\medskip

Note that, if $\tau$ is the initial slice of $H$, then if the marking
$\I(H)$ is clean it is 
compatible with $\tau$. The same is true for the terminal slice
and $\T(H)$. 

\bfheading{Partial order on slices:}
Consider now the set $V(H)$ of  complete slices whose bottom geodesic
equals the main geodesic of $H$. This set admits a partial order $\sprec$
as follows.
For $\tau, \tau'\in V(H)$, say that $\tau \sprec
\tau'$ iff $\tau \ne \tau'$ and, for any $(h,v) \in \tau$,
either $(h,v)\in\tau'$ or
there is some $(h',v')\in \tau'$
such that $(h,v) \pprec (h',v')$. 

\begin{lemma}{sprec partial order}
Let $H$ be a complete hierarchy. The relation $\sprec$ is a strict
partial order on $V(H)$. 
\end{lemma}

\begin{pf}
Let us first note the following facts:
\begin{enumerate}
\item $\pprec$ is a strict partial order,
\item Any two elements $(h,v)$ and $(k,w)$ of a slice $\tau$ are not
$\pprec$-comparable,
\item If $\tau\subseteq \tau'$ for slices $\tau,\tau'\in V(H)$ then
  $\tau = \tau'$.
\end{enumerate}

Fact (1) is Lemma \ref{Time Order} part (5) , and 
Fact (2) is an application of Lemma \ref{diagonalizable} and 
\ref{Time Order} part (1). 
Fact (3) follows from the fact that slices in $V(H)$ are {\em
  complete}.

For $\sprec$ to be a strict partial order it suffices to show that it
is transitive, since by definition it is never reflexive. Let 
$\tau_1 \sprec \tau_2 \sprec \tau_3$  for $\tau_i\in V(H)$.

By definition of $\sprec$, given any $p_i\in\tau_i$ (where $i=1,2$
and $p_i$ denotes some pair $(h_i,v_i)$), 
there
exists $p_{i+1} \in \tau_{i+1}$ such that either $p_i\pprec p_{i+1}$
or $p_i = p_{i+1}$.
By fact (1) this implies either $p_1 \pprec p_3$ or
$p_1 = p_3$. Thus either $\tau_1 \sprec \tau_3$ or $\tau_1=\tau_3$. To
rule out the latter, note that there is at least one $p_1\in\tau_1$
which is not in $\tau_2$ (by Fact (3)). Thus
$p_1\pprec p_2$ so $p_1\pprec p_3$, and $p_3$ cannot lie in $\tau_1$
by Fact (2).
\end{pf}

\bfheading{Forward elementary moves:}
Roughly, an elementary move on a slice $\tau$ consists of incrementing
the vertex $v$ of some $(h,v)$ in $\tau$, and making certain
adjustments to the other pairs to obtain a new slice $\tau'$. 

To begin, let $h\in H$ and let $v$ be a vertex of $h$, not the last,
with successor $v'$.  These will determine two slices $\sigma$ and
$\sigma'$, not necessarily complete, called the {\em transition
  slices} for $v$ and $v'$. The slices $\sigma,\sigma'$ will have the
property that (at least when $\xi(D(h)) > 4$) $\mu_\sigma = \mu_{\sigma'}
= v\union v'$. After constructing these we will extend them to
complete slices $\tau,\tau'$, which will constitute our elementary
move.

Define $\sigma$ as the smallest slice with bottom pair $(h,v)$ such that,
for any $(k,w)\in\sigma$ and $Y$ a component domain of $(D(k),w)$,
\begin{itemize}
\item[E1:] if $v'\rest Y \ne \emptyset$ and $Y$ supports a geodesic
  $q$ then $(q,u)\in \sigma$ where $u$ is the {\em last} vertex of $q$.

\item[E2:] if $ v'\rest Y = \emptyset$ then no geodesic in $Y$ is
  included in $\sigma$.
\end{itemize}
Note that $\sigma$ is easily built inductively from E1 and E2, and is
uniquely determined. It is also easy to check that it satisfies the
slice properties (S1--S3).
We call the domains appearing in E2 ``unused
domains'' for $\sigma$.
Similarly, define $\sigma'$ as the smallest slice with bottom
pair $(h,v')$, such that for any
$(k,w)\in\sigma'$ and $Y$ a component domain of $(D(k),w))$,
\begin{itemize}
\item[E1':] if $ v\rest Y \ne \emptyset$ and $Y$ supports a geodesic
  $q$ then $(q,u)\in \sigma'$ where $u$ is the {\em first} vertex of $q$.

\item[E2':] if $ v\rest Y = \emptyset$ then no geodesic in $Y$ is
  included in $\sigma$.
\end{itemize}
Before continuing let us consider this construction in several special
cases. 

\begin{enumerate}
\item $D(h)$ is an annulus. Here $v$ and $v'$ are arcs 
 in the closed   annulus with disjoint   interiors, and 
  $\sigma = \{(h,v)\}$, $\sigma' =  \{(h,v')\}$.
\item $D(h)$ is a once-punctured torus. Now $v$ and $v'$ are curves
  intersecting once in $D(h)$. Let $k$ be the geodesic supported in the
  annulus $Y$ whose core is $v$, and let $k'$ be the geodesic supported in
  the annulus $Y'$ whose core is $v'$. Then
$$\sigma = \{(h,v),(k,\pi_Y(v'))\},\ \ \  \sigma' =
\{(h,v'),(k',\pi_{Y'}(v))\}.$$ 
(If $D(h)$ is a 4-holed sphere then $v$ and $v'$ intersect twice so
$\pi_Y(v')=\T(k)$  has two components, only one of which appears as the
last vertex of $k$; and similarly for $\pi_{Y'}(v)$.)
\item $\xi(D(h)) = 5$. Now $v$ and $v'$ are disjoint one-component curves. The
  complementary domain $Y$ of $v$ with $\xi(Y)=4$ must contain $v'$, and the
  complementary domain $Y'$ of $v'$ with $\xi(Y')=4$ must contain
  $v$. Let $k$ and $k'$ be the geodesics supported in $Y$ and $Y'$
  respectively. Then 
$$\sigma = \{(h,v),(k,v')\}, \ \ \ \sigma' = \{(h,v'),(k',v)\}.$$
Note that the annuli with cores $v$ and $v'$ are not included in these
slices, by E2 and E2'. In particular we observe that $\mu_\sigma =
\mu_{\sigma'} = v\union v'$. The general $\xi>4$ case will be treated
in Lemma \ref{transition slices} below.
\end{enumerate}

\begin{lemma}{transition slices}
Let $v,v'$ be successive vertices in a geodesic $h\in H$, where $H$ is
a complete hierarchy. Let $\sigma,\sigma'$ be the transition slices
associated to $v,v'$. If $\xi(D(h)) > 4$ then 
no geodesics in $\sigma$ and $\sigma'$ have annular domains, the associated
markings $\mu_\sigma$ and $\mu_{\sigma'}$ have no transversals and
are both equal to $v\union v'$,
and the unused domains in $\sigma $ and $\sigma'$
are exactly the component domains of $(D(h),v\union v')$.
\end{lemma}
Thus, the move from $\sigma $ to $\sigma'$ in this case involves only
a ``reorganization'', and the underlying curve system is not changed.

\begin{pf}
Since $\xi(D(h))>4$, $v$ and $v'$ are disjoint curve systems. 
Consider first a component domain $Y$ of $(D(h),v)$. If $Y$ misses
$v'$ then it is an unused domain of $\sigma$ (case E2) and is also clearly a
component domain of $(D(h),v\union v')$. If $Y$ doesn't miss $v'$,
then by E1 (and completeness of $H$) we have
$(q,u)$ in $\sigma$ where
$D(q) = Y$ and $u$ is the last vertex in $q$. 

Since $Y\fsubd h$,
it follows from Lemma \ref{Y direct to top} that $q\fsubd h$
with $\T(q) = \T(Y,h)=v'\rest Y$.  Hence in particular $u\subseteq
v'$. Note that $u$ need not be all of $v'\rest Y$.

Now let $Z$ be  any component domain of $(D(q),u)$. By the above, the
relative boundary $\boundary_{D(h)}(Z)$ consists of some subset of 
$v\union v'$.
Again if $Z$ misses $v'$ it is unused in $\sigma$ and a component
domain of $(D(h),v\union v')$, 
and if $v'\rest Z\ne\emptyset$
then $Z\fsub q$ with $\T(Z,q) =  \T(q)\rest Z =  
 v'\rest Z$. Thus the same argument works inductively. The process
terminates in an unused domain exactly when this domain is a component
domain of $(D(h),v\union v')$.

The same argument applies to $\sigma'$ as well, reversing directions
as usual.

Every annulus whose core is a component of $v\union v'$ does not have
essential intersection with either $v$ or $v'$. Thus it is unused, so
that the slices $\sigma,\sigma'$ have no annulus-domain geodesics, and
their markings have no transversals.
\end{pf}

We can now define our elementary moves. 
Let there be given two slices $\tau$ and $\tau'$, and let $h$ be a
geodesic in $H$ with two successive vertices $v,v'$. 
We say that $\tau'$
is related 
to $\tau$ by a {\em forward elementary move} along $h$ from $v$ to
$v'$ (or $\tau \to
\tau'$ for short) provided the following holds: 
Letting $\sigma,\sigma'$ be the transition slices for $v,v'$, we have
$\sigma\subset \tau$ and $\sigma'\subset \tau'$, and $\tau\setminus\sigma =
\tau'\setminus\sigma'$. The next lemma checks that a forward move in
$V(H)$ really moves forward in terms of the partial order:

\begin{lemma}{forward means forward}
Suppose $\tau$ and $\tau'$ are in $V(H)$, and are related by an
elementary move $\tau\to\tau'$. Then
$\tau \sprec \tau'$.
\end{lemma}

\begin{pf}
First, $\tau \ne \tau'$ since $\sigma $ and $\sigma'$ differ in their
bottom pair. Now
let $(k,w)\in \tau$, such that $(k,w)\notin\tau'$. Then
$(k,w)\in\sigma$ and hence $D(k)\subseteq D(h)$, and 
$v'\rest {D(k)}\ne \emptyset$, by construction of $\sigma$. 
If $k=h$ then $(k,w) = (h,v) \pprec (h,v')\in \sigma'$ and we are done.
If not then $\phi_h(D(k))$ contains $v$ and not $v'$, so that
$\max\phi_h(D(k)) = v < v'$. We therefore have
$(k,w) \pprec (h,v')$, and again we are done.
\end{pf}

Next, we should show that in fact a sequence of elementary moves does
exist connecting the initial to the terminal slice of $H$, and
furthermore give a bound for its length.
Let $|H|$ denote the size of the hierarchy $H$, defined as the sum 
$\Sigma_{h\in H} |h|$ of the lengths of its geodesics.

\begin{proposition}{Resolution into slices}
Any complete finite hierarchy $H$ admits a sequence of forward elementary
moves 
$\tau_0\to\cdots\to\tau_N$ where $\tau_0$ is its initial slice, 
$\tau_N$ is its terminal slice, and 
$$
N\le |H|.
$$
Such a sequence is called a {\em resolution} of $H$.
\end{proposition}

\begin{pf}
Let us first show that, if $\tau$ is not the terminal slice of $H$,
then there exists some $\tau'$ such that $\tau\to\tau'$.
Indeed, there is at least one $(h,v)\in\tau$ for which $v$ is not the
last vertex 
of $h$. Choose $h$ minimal in the sense that if $(k,w)\in\tau$ and
$D(k)\subset D(h)$ then $w$ is the last vertex of $k$.
Let $v'$ denote the successor of $v$ in $h$.
The subset
$$
\sigma = \{(k,w)\in\tau: D(k)\subseteq D(h), 
v'\rest{D(k)}\ne\emptyset\}
$$
satisfies conditions (E1,E2), by the minimal choice of $h$ and the
fact that $\tau$ is complete.
Construct 
$\sigma'$ via (E1') and (E2'), thus obtaining the transition slices
for $v,v'$, and let $\tau'  =\sigma'\union(\tau\setminus \sigma)$.
It is easy to check that $\tau'$ satisfies conditions (S1--S3) and is
hence a slice, with $h_{\tau'} = h_\tau = g_H$. To see that it is a
complete slice (S4), consider  
any $(k,y)\in \tau'$ and let $Y$ be a component domain of
$(D(k),y)$. If $(k,y)$ is not in $\sigma'$ then by definition it is in
$\tau$, and since $\tau$ is complete, it contains some $(l,z)$ with
$D(l)=Y$. If $(l,z)\notin\sigma$ then again by definition $(l,z)\in
\tau'$ and we are done. If $(l,z)\in\sigma$ then, since it is a
component domain of a pair outside of $\sigma$, it can only be the bottom
pair $(h,v)$ of $\sigma$. But then $Y=D(h)$ and we know that $(h,v')$
is a pair in $\tau'$, so again we are done. Now suppose that
$(k,y)\in\sigma'$. If $Y$ is a 
used domain of $\sigma'$ then by definition it supports some geodesic
$l'$ appearing in $\tau'$. If $Y$ is an unused domain then by Lemma
\ref{transition slices} it is also an unused domain of $\sigma$, and
hence supports a 
geodesic $q$ appearing in $\tau\setminus \sigma$. Again we conclude
$q$ appears in $\tau'$ and we are done.

We thus have a slice $\tau'$ in $V(H)$, and an elementary forward move
$\tau\to\tau'$. Note that $\tau'$ is not uniquely determined by
$\tau$, as there may have been more than one $h$ to choose from.

Now if $\tau_0$ is the initial slice of $H$, the above implies that
there is a sequence 
$\tau_0 \to \tau_1 \to \tau_2 \to \cdots$, which does not terminate at
$\tau_i$ as
long as $\tau_i$ is not the terminal slice.
On the other hand by Lemma \ref{forward means forward}, 
the sequence is strictly increasing in $\sprec$. Since the set of
slices is finite, it must terminate for some $\tau_N$, which must then
be the terminal slice of $H$.

All we have left to prove is the bound on the length $N$ of the
resolution. 

Suppose a pair $(h,v)$ appears in $\tau_n$ and $(h,w)$ appears in
$\tau_{m}$ for $n<m$. Then $\tau_n \sprec \tau_m$ as we have seen, 
and therefore it must be that  $v\le w$.
For if not we would have $(h,w)\pprec(h,v)$, but by definition of
$\sprec$ there is some $(k,u)\in\tau_m$ such that
$(h,v)\pprec(k,u)$. 
Hence $(h,w) \pprec (k,u)$, but 
this contradicts the fact that all pairs in a
given slice are not $\pprec$-comparable (see proof of Lemma \ref{sprec
  partial order}).

By the definition, a forward move $\tau_n \to \tau_{n+1}$ advances exactly one
geodesic exactly one step, 
erases certain pairs of the form $(k,u)$ where $u$ is
the last vertex, and creates certain pairs of the form $(k',u')$ where $u'$ is
the  first vertex, and keeps the rest of the pairs unchanged.
Since by the previous paragraph
no vertices in any geodesic can be repeated once they have been
incremented or erased,
it follows that the number of forward moves is bounded by
$\Sigma_{h\in H} |h|$, which is $|H|$.
\end{pf}


\bfheading{Conversion to a sequence of clean markings.}
Given a resolution $\tau_0\to\cdots\to\tau_N$ of $H$ into slices, 
we may obtain a sequence of clean markings $\mu_0,\ldots,\mu_N$ by
requiring that each $\mu_i$ be compatible with $\tau_i$. Recall that
there may be a finite number of choices for each $\mu_i$. 
For convenience we also assume that $\I(H)$ and $\T(H)$ are clean, and
$\mu_0 = \I(H)$ and $\mu_N=\T(H)$.

What is left to check is the relationship between $\mu_i$ and
$\mu_{i+1}$. Recall from \S\ref{markings} the definition of the
elementary moves Flip and Twist on clean markings.
We can now establish:
\begin{lemma}{marking move bound}
Let $(\tau_i)$ be a resolution of a complete finite hierarchy $H$, and
let $(\mu_i)$ be 
a sequence of complete clean markings compatible with
$(\tau_i)$.
There exists $B>0$ depending only on the topological type of $S$, such
that
$\mu_i$ and $\mu_{i+1}$ differ by at most $B$ elementary
moves.

In particular, assuming $\I(H)$ and $\T(H)$ are clean, 
there is a sequence of clean markings $(\hat\mu_j)_{j=0}^M$, successive
ones separated 
by elementary moves, such that $\hat\mu_0 = \I(H)$, $\hat \mu_M=\T(H)$, and
$M\le B|H|$.
\end{lemma}

\begin{proof}
We have already seen in the beginning of the section that two clean
markings compatible with the 
same $\tau_i$ differ by a bounded number
of Twist elementary moves.

Now, recall that $\tau_i\to\tau_{i+1}$ is determined by a
transition $v\to v'$ along some geodesic $h$. If $D(h)$ is an annulus,
$v$ and $v'$ differ by distance one in the annular  complex $\CC(D(h))$,
so a bounded number of Twist moves applied to $\mu_i$ yields a
marking $\mu'_{i+1}$ 
which is compatible with $\tau_{i+1}$. Then as above $\mu'_{i+1}$ and
$\mu_{i+1}$ are related by a bounded number of Twist moves.

Suppose that $\xi(D(h)) = 4$. Then recall that the transition slices
$\sigma_i$ and $\sigma_{i+1}$ can be written as
$\{(h,v),(k,t)\}$ and $\{(h,v'),(k',t')\}$ where $k,k'$ are the
geodesics in the complexes of the annuli $Y$ and $Y'$ with cores $v$
and $v'$ respectively, and $t$ and $t'$ are vertices of $\pi_Y(v')$
and $\pi_{Y'}(v)$ respectively.  (If $D(h)$ is a 1-holed torus then
$t=\pi_Y(v')$ and $t'=\pi_{Y'}(v)$).
Thus a clean marking $\mu'_i$ can be constructed compatible with $\tau_i$ and
containing a pair $(v,\pi_Y(v'))$. Now a Flip move on this
marking yields a marking $\mu'_{i+1}$ with the pair $(v',\pi_{Y'}(v))$,
with all other base curves the same, and transversals at distance at most
$\dcompat$  from those of $\mu'_i$ by Lemma \ref{clean markings}.
It follows that, using a bounded number of Twist moves on each
base curve, 
 $\mu'_{i+1}$ can be made into $\mu''_{i+1}$ which is
compatible with $\tau_{i+1}$.
Since the previous discussion bounds the number of moves to get from
$\mu_i$ to $\mu'_i$ and from 
$\mu_{i+1}$ to $\mu''_{i+1}$, we again have a bound on the number of
moves needed to get from $\mu_i$ to $\mu_{i+1}$.

Finally when $\xi(D(h)) > 4$, $\tau_i$ and $\tau_{i+1}$ have exactly
the same base curves, and the positions on their annulus geodesics are
the same. It follows that any marking compatible with $\tau_i$ is also
compatible with $\tau_{i+1}$, and hence again $\mu_i$ and $\mu_{i+1}$
differ by a bounded number of Twist moves.
\end{proof}

We remark that explicit bounds for this lemma are straightforward, but
somewhat tedious, to compute, so we have elected to leave them out.

%% file: control.tex
\section{Comparison and control of hierarchies}
\label{large link etc}

In this section, we combine the structural results of the previous two
sections with Theorem \lref{Bounded Geodesic Image}, to
prove a number of basic results that
allow us to control the higher-order structure of hierarchies, and
to compare hierarchies whose main geodesics are close.

As applications we prove Theorem \lref{Efficiency of Hierarchies}, 
which shows that hierarchies give rise to sequences of markings 
separated by elementary moves which are close to shortest
possible. These will be used to produce quasi-geodesics in $\Mod(S)$ in
Section \ref{MCG}. We also prove Theorem \lref{Convergence of
  Hierarchies} which will allow us to obtain infinite hierarchies as
limits of finite ones. Corollary \lref{Finite Geodesics} states that
between any two points in $\CC(S)$ there are only finitely many tight
geodesics. 

Our basic technical result will be Lemma \lref{Sigma Projection},
which simplifies and generalizes the ``short cut and projection''
argument used in the motivating examples in \S\ref{example}. Recall
how we
showed that a large link (long geodesic) in one hierarchy forces a
similar large link in a fellow-traveling hierarchy, by
producing paths forward and backwards
from the given link to its main geodesic, and projecting these back to
the domain of the link. The forward and backward sequences
$\Sigma^\pm$ provide the framework for making this argument work in
general.

Lemmas \lref{Large Link} and \lref{Common Links} will be
straightforward applications of Lemma \ref{Sigma Projection}, and will
generalize what we did in the motivating examples. Lemma \lref{Slice
  Comparison} is a more delicate comparison between nearby hierarchies
and requires more work. 

\subsection{The forward and backward paths}

The ``forward path'' for a domain $Y$ is built roughly as follows:
Starting on the top geodesic in $\Sigma^+(Y)$ we move forward until
it ends, at which point we have arrived at the position following the
footprint of $Y$ on the next geodesic in $\Sigma^+(Y)$, and we
continue in this way until we get to the bottom geodesic $g_H$. A
``backward path'' is constructed the same way from $\Sigma^-(Y)$.

More precisely: Let
$\sigma$ denote the set of all pairs $(k,v)$ 
where $k\in\Sigma^\pm(Y)$, and $v$ is a position on $k$ such that
$v\rest Y \ne \emptyset$.
We claim that the partial order $\pprec$
restricts to a linear order on $\sigma$, making it into a
sequence: 

Indeed, each $f_i\in\Sigma^+(Y)$ for $i>0$ contributes
a segment
$\sigma^+_i = \{(f_i,v_i)\pprec\cdots\pprec(f_i,\T(f_i))\}$, where
$v_i$ is
the position immediately following $\max\phi_{f_i}(Y)$ (if 
$\max\phi_{f_i}(Y)$ is the last vertex then $\sigma^+_i =
\{(f_i,\T(f_i))\}$). Since 
$\max\phi_{f_i}(Y) = \max\phi_{f_i}(D(f_{i-1}))$ (Corollary
\ref{Footprints}), we also have $(f_{i-1},\T(f_{i-1})) \pprec
(f_i,v_i)$.
Thus the union of all $\sigma_i^+$ are linearly ordered.
The same
holds for $\sigma^-_i$, defined as 
$\{(b_i,\I(b_i))\pprec\cdots\pprec(b_i,u_i)\}$, where $u_i$ is the
last position before $\min\phi_{b_i}(Y)$.
Note that the same geodesic may appear in
$\Sigma^+$ and
$\Sigma^-$, in which case it can contribute both a $\sigma_j^+$
and a
$\sigma_i^-$, one on each side of the footprint.

The top geodesic $h=b_0=f_0$ has 
empty $\phi_{h}(Y)$ by Theorem \ref{Structure of Sigma} part (2),
and so all its positions are included in $\sigma$, and
they
follow all the $\sigma_i^-$ and 
precede all the $\sigma_i^+$ pairs, for $i>0$. 
We denote the sequence of positions of the top geodesic by $\sigma^0$.

We let $\sigma^+$ be the concatenation
$\sigma_1^+ \union \cdots \union \sigma_n^+$ (with the same linear
order), and similarly 
$\sigma^- = \sigma_m^- \union \cdots \union \sigma_1^- $.
In case clarification is needed we write $\sigma^+(Y),
\sigma^-(Y,H)$,
etc.

By the definition, for each $(k,v)\in\sigma$ the projection
$\pi_Y(v)$ is nonempty. 
Let $\pi_Y(\sigma)$ denote the union of these projections, and
similarly for $\pi_Y(\sigma^+)$ and $\pi_Y(\sigma^-)$.
The following property of $\sigma$ forms the basis of all the
proofs in this section:

\begin{lemma+}{Sigma Projection}
There exist constants $M_1,M_2$ depending only on $S$ such that,
for any hierarchy  $H$ and domain $Y$ in $S$, 

$$\diam_{Y}(\pi_Y(\sigma^+(Y,H))) \le M_1$$
and similarly for $\sigma^-$. 

Furthermore, if $Y$ is properly contained in the top domain of
$\Sigma(Y)$, then
$$\diam_{Y}(\pi_Y(\sigma(Y,H))) \le M_2.$$
\end{lemma+}

\begin{pf}
Theorem \lref{Bounded Geodesic Image} bounds the diameter of the
projection to $Y$ of each $\sigma^\pm_i$, and of $\sigma^0$ in the
case where $Y$ is properly contained in the top domain. The
transition from the last position of $\sigma^+_i$ to the first of
$\sigma^+_{i+1}$ just
consists of adding curves to the marking and so projects to a bounded
step in $\CC(Y)$ by Lemma \ref{Lipschitz Projection}.
The same holds for the other
transitions between segments of $\sigma$. Finally, 
the number of segments in each of
$\sigma^\pm_i$ is bounded by $\xi(S)-\xi(Y)$. These facts together
give the desired diameter bounds.
\end{pf}

\subsection{Large links}
The following is an almost immediate consequence of Lemma
\ref{Sigma Projection}:

\begin{lemma+}{Large Link}
If $Y$ is any domain in $S$ and 
$$d_Y(\I(H),\T(H)) > M_2$$
then $Y$ is the support of a geodesic $h$ in $H$.

Conversely if  $h\in H$ is any geodesic with $Y=D(h)$, 
$$
\left| |h| - d_Y(\I(H),\T(H)) \right | \le 2M_1.
$$
\end{lemma+}

\begin{proof}
The top geodesic $k=b_0=f_0$ of $\Sigma(Y)$ has domain $Z=D(k)$ which
either equals $Y$ or contains it.
If $Y$ does not support a geodesic then
$Z$ properly contains $Y$, and Lemma
\ref{Sigma Projection} implies
$$d_Y(\I(H),\T(H)) \le \diam_{Y}(\pi_Y(\sigma)) \le M_2.$$
This proves the first part. 

For the second part, if $Y=D(h)$ then by Theorem \ref{Structure of
  Sigma} we must have $Z=Y$ and $h=k$. Since $\sigma^+$
contains both $\T(h)$ and $\T(H)$, and $\sigma^-$ contains $\I(h)$
and $\I(H)$, Lemma \lref{Sigma Projection} implies that
$$d_Y(\I(h),\I(H)) \le M_1$$ and
$$d_Y(\T(h),\T(H))) \le M_1.$$
The second statement of the lemma follows.
\end{proof}

%
%

\subsection{Fellow traveling}
In a $\delta$-hyperbolic metric space, geodesics whose endpoints are
near each other must stay together for their whole lengths.  Our
hierarchies have some similar properties. Before we state them we need
some definitions.

\begin{definition}{separation}
We say that two hierarchies $H$ and $H'$ are {\em $K$-separated at the
  ends}
if the markings $\I(H)$ and $\I(H')$ are complete and clean, and are
separated by at most $K$
elementary moves, and similarly for $\T(H)$ and $\T(H')$.
\end{definition}

\begin{definition}{K R parallel}
Given two geodesics $g_1$ and $g_2$ with the same domain, 
and $x_i$ a vertex in $g_i$ for $i=1,2$,
we say that
$g_1$ and $g_2$ are {\em $(K,R)$-parallel} at $x_1$ and $x_2$ 
provided $d(x_1,x_2)\le K$ and for at least one of $i=1$ or $2$,
$x_i$ is the midpoint of a segment $L_i$ of 
radius $R$ in $g_i$, and $L_i$ lies in a $K$-neighborhood of $g_{3-i}$.
\end{definition}

\begin{definition}{pseudo-parallel}
We say a hierarchy $H$ is {\em $(K,M)$-pseudo-parallel} to a hierarchy
$H'$ if, for any geodesic $h\in H$ with $|h|\ge M$ there is a geodesic
$h'\in H$ such 
that $D(h) = D(h')$, and $h$ is contained in a $K$-neighborhood of $h'$
in $\CC(D(h))$.
\end{definition}
(Note that the pseudo-parallel relation is not symmetric)

The following lemma is a generalization of Farb's Bounded Coset
Penetration Property. 

\begin{lemma+}{Common Links}
Given $K$ there exist $K', M$ such that, 
if two hierarchies $H$ and $H'$ are $K$-separated at the ends
then each of them is $(K',M)$-pseudo-parallel to the other.
\end{lemma+}

\begin{proof}
Let $h$ be any geodesic in $H$.
By Lemma \ref{Large Link}, the hypothesis, and Lemma
\ref{Elementary Move Projections}, we have 
$$|h|-2M_1-4K\leq d_Y(\I(H'),\T(H'))\leq   |h|+2M_1+4K.$$
If we assume $|h|>M_2 + 2M_1 + 4K$, then 
the left hand side is greater  than $M_2$,
so Lemma \ref{Large Link} implies
that there is a geodesic $h'\in H'$ with $D(h')=D(h)$.
A bound of $2M_1 + 8K$
on $d_Y(\I(h),\I(h'))$ and $d_Y(\T(h),\T(h'))$ follows
from Lemmas \ref{Sigma Projection} and \ref{Elementary Move Projections}.
It follows by hyperbolicity of $\CC(Y)$ that $h$ and $h'$ 
remain a bounded distance apart along their whole length.
\end{proof}

In the next lemma we show how to compare slices in a pair of
hierarchies that are $K$-separated at the ends or have parallel
segments. The idea is that two 
such slices can be joined by a hierarchy that only has long geodesics
when these are parallel to segments in the original two hierarchies. 
This is the closest one can come to saying that two hierarchies are
fellow-travelers. 

\begin{lemma+}{Slice Comparison}
Given $K$ there exist $K'$, $M$ so that the following holds:
Let $\tau $ and $\tau'$ be complete slices in two hierarchies $H$ and $H'$
respectively, with bottom vertices $x\in g_H$ and $x'\in g_{H'}$.
Suppose that either
\begin{enumerate}
\item $H$ and $H'$ are $K$-separated at the ends, or
\item $g_H$ and $g_{H'}$ are $(K,3K+4)$-parallel at $x$ and $x'$.
\end{enumerate}
Let $\mu$ and $\mu'$ be clean markings compatible with $\tau$ and
$\tau'$ respectively.
Then any hierarchy $J$ with $\I(J) = \mu$ and $\T(J)=\mu'$ is
$(K',M)$-pseudo-parallel to both $H$ and $H'$.
\end{lemma+}

Before giving the proof of this lemma we need the following two results.

\begin{lemma}{slices cut}
Let $H$ be a complete hierarchy.
Let $\tau$ be a slice in $V(H)$ and $(k,v)$ a pair where $v$ is a
position in $k\in H$. Then exactly one of the following occurs:
\begin{enumerate}
\item $(k,v)\in\tau$,
\item there exists
$(h,u)\in\tau$ such that $(k,v)\pprec(h,u)$,
\item there exists
$(h,u)\in\tau$ such that $(h,u)\pprec(k,v)$.
\end{enumerate}
Furthermore $h$ may be taken so that $D(k)\subseteq D(h)$.
\end{lemma}

In view of this result, let us write $(k,v) \sprec \tau$ when case (2)
holds, and  $\tau\sprec(k,v)$ when case (3) holds.

\begin{proof}
  Since $\tau\in V(H)$, it is complete and its bottom geodesic is $g_H$.
  We will prove the statement of the lemma inductively for any
  complete slice whose bottom geodesic $g$ satisfies $D(k)\subseteq
  D(g)$. Let $(g,u)$ be the bottom pair of $\tau$. If $g=k$ then the
  statement is immediate -- either $v<u$, $u<v$, or $u=v$.

  Now suppose $D(k)$ is properly contained in $D(g)$. 
By Corollary \ref{easy containment}, $\phi_g(D(k))$ is nonempty.
If $\max\phi_g(D(k)) < u$ or $u<\min\phi_g(D(k))$ then $(k,v)\pprec
  (g,u)$ or $(g,u)\pprec (k,v)$, respectively, and we are done. If not
  then $u\in\phi_g(D(k))$ and there is some component domain $Y$ of $(D(g),u)$
  containing $D(k)$. Since $\tau$ is complete there is a pair
  $(h,w)\in\tau$ with $D(h) = Y$. The slice $\tau'$ consisting of all
  $(h',w')\in\tau$ such that $D(h')\subseteq D(h)$ is itself complete,
  and has bottom pair $(h,w)$. Applying induction to $\tau'$, we have
  the desired statement.  

  The fact that the three possibilities are mutually exclusive follows
  directly from the fact that any two elements of a slice are not
  $\pprec$-comparable (see proof of Lemma \ref{sprec partial order}).
\end{proof}

\begin{lemma}{sigma and slice}
Fix a complete hierarchy $H$ and a slice $\tau\in V(H)$.
Let $Y$ be any domain in $S$. Then the path $\sigma(Y)$
contains a pair $(k,v)$ which is in $\tau$.
\end{lemma}

\begin{proof}
Let $(k,v)$ and $(k',v')$ be succesive pairs in
$\sigma(Y)$. We will show that it is not possible for 
$(k,v) \sprec \tau$ and $\tau\sprec (k',v')$ to hold simultaneously.
 Since the first pair in $\sigma$ is always
$(g_H,\I(g_H))$, for which $(g_H,\I(g_H))\sprec \tau$ holds,
and the last is $(g_H,\T(g_H))$ for which $\tau\sprec(g_H,\T(g_H))$ holds,
the statement of the 
lemma follows from Lemma \ref{slices cut}.

By definition of $\sigma$, there are three possibilities for the
relation between $(k,v)$ and $(k',v')$: 
\begin{enumerate}
\item $k=k'$. Here $v'$ is the position following $v$. 
\item $k\fsubd k'$. Here $v=\T(k)$, and $v'$ is the position following
  $\max\phi_{k'} (D(k))$.
\item $k\bsubd k'$. Here $v'=\I(k')$, and $v$ is the position preceding
  $\max\phi_{k} (D(k'))$.
\end{enumerate}

We will first prove our claim in cases (1) and (2). Suppose that
$(k,v)\sprec\tau$, and 
let $(h,u)\in\tau$ be a pair such that, as in Lemma \ref{slices
  cut}, $(k,v)\pprec(h,u)$ and $D(k)\subseteq D(h)$. 

If $k=h$ then $v<u$ and in particular we must be in case (1) since $v$
  is not the last position of $k$. Thus $k=k'=h$ and $v'$ is the
  successor of $v$, so $v'\le u$. Thus we are done in this case.

If $D(k)$ is  properly contained in $D(h)$ then we note that $k\fsub
h$. In case (1) we still 
have $(k',v') = (k,v') \pprec (h,u)$, so we are done. In case (2),
we either have $k'=h$ or $k'\fsub h$. In the first case, $v'$ is the
successor of $\max\phi_h(D(k))$ and hence $v'\le u$ and we are done. In
the second, we have $\max\phi_h(D(k')) = \max\phi_h(D(k))$ by
Corollary \ref{Footprints}, and so again $(k',v')\pprec (h,u)$.

To prove our claim in case (3), we just note that it is equivalent to
case (2), with the directions and roles of $k$ and $k'$ reversed.
\end{proof}

\begin{proof}[Proof of Slice Comparison Lemma]
Let  $R= 3K+4$. 
Let $m_0$ be any geodesic of $J$ with $|m_0|> M$, where the value of
$M$ will be determined below, and let $Y=D(m_0)$.
We first claim that, up to possibly reversing all the directions in
$H'$ (and interchanging $\I(H')$ with $\T(H')$),
\begin{eqnarray}
\label{TT bound}
d_Y(\T(H),\T(H')) &\le M_3,\\
\label{II bound}
d_Y(\I(H),\I(H')) &\le M_3 
\end{eqnarray}
for appropriate $M_3$.
(If $H$ is infinite then this holds with $\I(H)$ or $\T(H)$ replaced
by any point of $g_H$ on $\sigma^-(Y)$ or $\sigma^+(Y)$, respectively;
and similarly for $H'$).

In case (1) this is true by the hypothesis and Lemma 
\ref{Elementary Move Projections}, provided $M_3\ge 4K$. 

In case (2), up to interchanging $H$ and $H'$ we
may assume there is an interval $L$ of radius $R$ centered on $x$
which lies in a $K$-neighborhood of $g_{H'}$.
If $y,z$ are the endpoints of $L$ and $y<x<z$, let $y',z'$ be points of
$g_{H'}$ closest to $y$ and $z$ respectively. Up to reversing all the
directions in $H'$ we may assume $y'<z'$, and then by the triangle
inequality we have $y'<x'<z'$.

Since $Y$ is a domain of $J$, whose main geodesic has
length at most $K$, and $x$ is a curve in $\I(J)$,
we have $d_S(\boundary Y,x) \le K+2$ (thinking of $\boundary Y$ as a
simplex in $\CC(S)$).
We claim that any curve $w$
on a geodesic from $z$ to $z'$ intersects $\boundary Y$. 
For if not, $d_S(w,\boundary Y) \le 1$, and so
$d(x,z)\leq d(x,\boundary Y)+d(\boundary Y,w)+d(w,z)\leq 2K+3 < R$, a
contradiction. It follows that we can project $[z,z']$ into $\CC(Y)$ and
conclude $d_Y(z,z') \le 2K$ by Lemma \ref{Lipschitz Projection}. 

If $\phi_{g_H}(Y)$ is nonempty,
the triangle inequality similarly gives
$d_S(x,\phi_{g_H}(Y)) \le K+3 < R = d(x,z) $ and hence
$z$  lies to the right of $\phi_{g_H}(Y)$. Similarly
$d_S(x',\phi_{g_{H'}}(Y)) \le K+3 < R-2K \le d(x',z')$ and hence 
$z'$  lies to the right of $\phi_{g_{H'}}(Y)$, if that is nonempty.

It follows that $z$ can be connected to $\T(H)$ 
by a path lying in
$\sigma^+(Y,H)$, and similarly for $z'$ and $\T(H')$.
Lemma \lref{Sigma Projection} then gives a
bound of $M_1$ on
$d_Y(z,\T(H))$ and $d_Y(z',\T(H'))$.
Putting
these together with the bound on $d_Y(z,z')$
gives (\ref{TT bound}), with $M_3=2M_1 + 2K$.
The same argument with $y$ and $y'$ gives (\ref{II bound}).

Next we claim that, for $M_4= 2M_2 + 4M_1 + 4$,
\begin{equation}
\label{mu between}
d_Y(\mu,\I(H)) + d_Y(\mu,\T(H)) \le d_Y(\I(H),\T(H)) +M_4
\end{equation}
and similarly for $\mu'$ and $H'$.  

Begin by observing that, 
by Lemma \ref{sigma and slice}, $\sigma(Y,H)$ contains a pair
$(k,v)\in\tau$.  If $D(k)$ is an annulus then $D(k)=Y$ and $v$ is a
transversal of $\mu$ -- otherwise it is in $\base(\mu)$.
Lemma \ref{Lipschitz   Projection} then implies that $\pi_Y(\mu)$ is within
distance 2 of $\pi_Y(\sigma(Y,H))$.

By Lemma \lref{Sigma Projection},
if $Y$ does not support a geodesic in $H$ then $\diam_Y(\sigma(Y,H)) \le
M_2$. Hence the left side of 
(\ref{mu between}) is at most $2M_2+4$ and the inequality
follows by choice of $M_4$.

If $Y$ supports a geodesic $h\in H$ then Lemma \ref{Sigma Projection}
implies that $\pi_Y(\sigma(Y,H))$ is Hausdorff distance $M_1$ from
$h$ (i.e. each is in an $M_1$-neighborhood of the other). We therefore
have, for $v$ as above, 
$d_Y(v,v_0)  + d_Y(v,v_{|h|}) \le |h| + 2M_1$,
where $v_0$ and $v_{|h|}$ are the first and last vertices of $h$,
and since $\pi_Y(\I(H))$ and $\pi_Y(\T(H))$ are distance $M_1$ from the
respective endpoints of $|h|$, (\ref{mu between}) again follows
with a bound of $6M_1 + 4$. This is at most $M_4$ since $M_1\le M_2$.

Now by the triangle inequality $d_Y(\mu,\mu') $ is bounded by
both 
$$d_Y(\mu,\T(H)) + d_Y(\T(H),\T(H')) + d_Y(\mu',\T(H'))$$
and
$$d_Y(\mu,\I(H)) + d_Y(\I(H),\I(H')) + d_Y(\mu',\I(H')).$$
Adding these two estimates together and applying (\ref{TT
  bound},\ref{II bound}) and (\ref{mu between})), we have 
$$
2d_Y(\mu,\mu')\leq
d_Y(\I(H),\T(H))+d_Y(\I(H'),\T(H')) + 2M_3 + 2M_4.
$$

Now again applying (\ref{TT bound},\ref{II bound}) and the triangle
inequality, we find that 
$d_Y(\I(H),\T(H))$ and $d_Y(\I(H'),\T(H'))$ differ by at most $2M_3$.
This gives 
$$d_Y(\mu,\mu')\leq d_Y(\I(H),\T(H)) + 4M_3  + 2M_4,$$
and the same inequality for $H'$.

Since $|m_0|> M$, Lemma \lref{Large Link} gives
$d_Y(\mu,\mu')> M - 2M_1$, so   
$d_Y(\I(H),\T(H))>M-2M_1-4M_3 - 2M_4$. 
If we set $M=2M_1 + 4M_3+2M_4 + M_2$, 
Lemma \ref{Large Link} again guarantees that
$Y$ is the domain of a geodesic
$m\in H$. Furthermore, $\T(m_0)$ is within $M_1$ of 
$\pi_Y(\T(J))=\pi_Y(\mu') $, which is within 4 of
$\pi_Y(\sigma(Y,H))$ by Lemma \ref{sigma and slice}, as above.
Applying Lemma \ref{Sigma Projection} again we find that this is
within $M_1$ of $m$. A similar estimate holds for $\I(m_0)$, and by
$\delta$-hyperblicity  all of $m_0$ lies within
$K'=2M_1 + 4 + 2\delta$ of $m$. This establishes
that $J$ is $(K',M)$-pseudo-parallel to $H$. 

The corresponding statement holds for $H'$, and the lemma is proved.
\end{proof}

\subsection{Efficiency}
\label{efficiency section}
In Section \ref{resolution} we saw that a hierarchy $H$ can be
resolved into a sequence of markings of length bounded by its
size $|H|$. Here we will obtain an estimate in the opposite direction.
Let $\til \MM$ be the graph whose vertices are complete, clean
markings in $S$, and whose edges  represent elementary moves. Giving
edges length 1, we have for two complete clean markings $\mu,\nu$
their {\em elementary move distance} $d_{\til \MM}(\mu,\nu)$ in this graph.
Proposition \ref{Resolution into slices}
and Lemma \ref{marking move bound} imply that 
this graph is connected, but this fact is already well known:
it follows for example from a similar connectedness result
for the graph of pants decompositions in Hatcher-Thurston
\cite{hatcher-thurston} (and see proof in Hatcher \cite{hatcher:pants}).

\begin{theorem+}{Efficiency of Hierarchies}
There are constants $c_0,c_1>0$ depending only on $S$ so that, 
if $\mu$ and $\nu$ are complete clean markings and
$H$ is a hierarchy with $\I(H)=\mu$, $\T(H)=\nu$, then

        $$
 c_0^{-1}|H| -c_1     \le      d_{\til\MM}(\mu,\nu)  \le c_0|H|.
        $$
\end{theorem+}

See Theorem  \lref{Quasigeodesic Words}, in Section \ref{MCG},
for the implication of this to words in the
Mapping Class Group.

\begin{proof}
The second inequality is an immediate consequence of 
Proposition \ref{Resolution into slices} and Lemma
\ref{marking move bound}. 

For the first inequality, 
the idea of the argument is as follows. Consider a 
shortest path $\{\mu=\mu_0,\ldots,\mu_N=\nu\}$
from $\mu$ to $\nu$ in $\til\MM$.
Each long geodesic $h\in H$ 
imposes a lower bound of the form $N \ge c|h|$, because the projection
$\pi_{D(h)}(\mu_j)$ 
moves at bounded speed in $\CC(D(h))$ (Lemma \ref{Elementary Move Projections})
as $j$ goes from 
$0$ to $N$, and by Lemma \lref{Large Link} it must travel a distance
proportional to $|h|$. To obtain our desired statement
we must show that the projections of $\mu_j$ cannot
move far in many different domains at once, and hence the lower bounds
for the different geodesics will add. This will be done using
Lemma \ref{Order and projections} below, which relates projections to
time-order. 

Let $M_5= 2M_1 +5$ and $M_6 = 4(M_1+M_5+4)$, and let
$\GG$ be the set of geodesics $h\in H$
satisfying $|h|\geq M_6$. Let 
$|\GG| = \Sigma_{h\in\GG} |h|$. Then we have
\begin{equation}
  \label{G is enough}
  |\GG| \ge d_0|H| - d_1
\end{equation}
for $d_0,d_1$ depending only on $S$ (and the choice of $M_6$). The
proof is a simple counting 
argument, using the fact that the number of component domains of any
geodesic is bounded by a constant times its length.
Thus the main point will be to bound $N$ below in terms of $|\GG|$.

For any $h\in\GG$ let us isolate an interval in $[0,N]$ in which the
projections to $Y=D(h)$ of the 
$\mu_j$ are ``making the transition'' between being close to
$\pi_Y(\mu_0)$, to being close to $\pi_Y(\mu_N)$.

Let $L=d_Y(\mu_0,\mu_N)$, noting that $L\ge |h|-2M_1 \ge M_6-2M_1$ by
Lemma \ref{Large Link}. The projections of $\mu_j$ to $\CC(Y)$ are a
sequence that moves by bounded jumps 
$d_Y(\mu_j,\mu_{j+1}) \le 4$, by
Lemma \ref{Elementary Move Projections}. 
Therefore there must be some largest value of $j\in[0,N]$ for which
$d_Y(\mu_0,\mu_j)\in[M_5,M_5+4]$. Let $a_Y$ be this value. 
Since $L>2(M_5+4)$, we know that $d_Y(\mu_{a_Y},\mu_N)> M_5+4$.
Therefore there 
is a smallest $j\in[a_Y,N]$ for which
$d_Y(\mu_j,\mu_N)\in[M_5,M_5+4]$. Let this be $b_Y$.

Let $J_Y$ be the interval $[a_Y,b_Y]$. 
These intervals have the following properties:

\begin{enumerate}
\item For any $j\in J_Y$ we have $d_Y(\mu_0,\mu_j)\ge M_5$ and
$d_Y(\mu_j,\mu_N) \ge M_5$.
\item $|J_{D(h)}|\geq |h|/8$ for any $h\in\GG$.
\item If $h,k\in \GG$ are such that $Y=D(h)$, $Z=D(k)$ have
nonempty intersection and neither is contained in the other,
then $J_Y$ and $J_Z$ are disjoint intervals.
\end{enumerate}

(1) follows immediately from the definition and Lemma 
\ref{Elementary Move Projections}

To prove (2), by the triangle inequality we have $d_Y(\mu_0,\mu_{b_Y}) \ge
L-M_5-4$. Again by Lemma \ref{Elementary Move Projections},
$d_Y(\mu_0,\mu_j)$ changes by at most $4$ with each increment
of $j$, so since $d_Y(\mu_0,\mu_{a_Y}) \le M_5+4$
we conclude that $b_Y-a_Y \ge (L-2M_5-8)/4$. This implies (2), by the
choice of constants and the fact that $|h|\ge M_6$.

To prove (3), we will first need the following lemma:
\begin{lemma+}{Order and projections}
Let $H$ be a hierarchy and $h,k\in H$ with $D(h)=Y$ and
$D(k)=Z$. Suppose that $Y\intersect Z \ne \emptyset$, and neither
domain is contained in the other. Then, if $h\tprec k$ then
\begin{equation}\label{order proj 1}
d_Y(\boundary Z, \T(H)) \le M_1+2
\end{equation}
and
\begin{equation}\label{order proj 2}
d_Z(\I(H),\boundary Y) \le M_1+2.
\end{equation}
\end{lemma+}

\begin{proof}[Proof of Lemma \ref{Order and projections}]
Let $m$ be
the geodesic used to compare $h$ and $k$.  
It lies in $\Sigma^+(Y)$.  Let $v\in
\phi_m(Z)$.   Since $v$ lies to the right of $\phi_m (Y)$ the
pair $(m,v)$ is in the sequence $\sigma^+(Y)$. 
Lemma \lref{Sigma Projection} now  implies
$$
d_Y(v,\T(H)) \le M_1.$$
Since $\boundary Z$ intersects $Y$ essentially by the assumption on
$Y$ and $Z$, and since $\boundary Z$ is disjoint from $v$, by applying
Lemma \ref{Lipschitz Projection} we find that
$$
d_Y(\boundary Z,\T(H)) \le M_1+2$$
as desired. The second inequality is proved in the same way.
\end{proof}

Returning to the proof of Theorem \ref{Efficiency of Hierarchies},
suppose that property (3) is false, so that $Y$ and $Z$ intersect
and are non-nested, but $J_Y$ and $J_Z$ overlap.
Let $j\in J_Y\intersect J_Z$. 

Let $H_j$ be a hierarchy such that $\I(H_j) = \mu_0$ and $\T(H_j)=\mu_j$.
By property (1) $d_Y(\mu_0,\mu_j)$ and $d_Z(\mu_0,\mu_j)$ are both at
least $M_5$, so Lemma 
\ref{Large Link} implies that $Y$ and $Z$ support geodesics $h_j$
and $k_j$ in $H_j$.
The condition on $Y$ and $Z$ implies that $h_j$ and
$k_j$ are time-ordered in $H_j$ by Lemma \ref{Time Order}, so suppose
without loss of generality that $h_j\tprec k_j$.
By Lemma \ref{Order and projections}, we have $d_Y(\boundary Z,\mu_j)
\le M_1+2$.

However, $h$ and $k$ must also be time-ordered in $H$, and 
applying Lemma \ref{Order and projections} to the hierarchy
$H$, we have either 
$$d_Y(\boundary Z,\mu_N) \le M_1+2$$
if $h\tprec k$ by (\ref{order proj 1}), or 
$$d_Y(\mu_0,\boundary Z) \le M_1+2$$
if $k\tprec h$ by (\ref{order proj 2}). Thus, either $d_Y(\mu_j,\mu_N) \le
2M_1 + 4$ or 
$d_Y(\mu_0,\mu_j) \le 2M_1+4$. Either one of these contradicts the
assumption that $j\in J_Y$, since $2M_1+4 <  M_5$.
This proves (3).

Thus, the intervals $\{J_{D(h)}: h\in\GG\}$ cover a
subset of $[0,N]$ with multiplicity at most $s$, where
$s$ is the maximal cardinality of a set $D_1,\ldots,D_s$ of domains in
$S$, $\xi(D_i)\ne 3$, for which
any two are
either disjoint or nested. This number depends only on $S$ (in fact it
is easy to show that $s=2\xi(S)-6$).
It follows that 
$$
sN \ge \sum_{h\in\GG} |J_{D(h)}|.
$$
Combining this with  (2) which gives  $\sum_{h\in\GG}|J_{D(h)}| \ge
|\GG|/8$, and then using (\ref{G is enough}), we obtain
$$
N \ge c_0^{-1}|H| - c_1
$$
with suitable constants $c_0,c_1$.
\end{proof}

The following corollary of this theorem can be stated without
any mention of hierarchies. It relates elementary-move distance to the
sum of all ``sufficiently large'' projections to subsurfaces in $S$
(including $S$ itself).

\begin{theorem+}{Move distance and projections}
  There is a constant $M_6(S)$ such that, given $M\ge M_6$, 
there are $e_0,e_1$ for which, if $\mu$ and $\nu$ are any two complete
clean markings then
$$
e_0^{-1} d_{\til \MM}(\mu,\nu) -e_1 \le
\sum_{\substack{
Y\subseteq S\\
d_Y(\mu,\nu)\ge M}}         d_Y(\mu,\nu) \le
e_0 d_{\til \MM}(\mu,\nu) + e_1
$$
\end{theorem+}
The proof is simply a rephrasing of the result of Theorem
\ref{Efficiency of Hierarchies}, together with the inequality
(\ref{G is enough}) to restrict consideration to ``long'' geodesics,
and Lemma \ref{Large Link} to relate this to projection diameters.

\subsection{Finiteness and limits of hierarchies}
\label{limits of hierarchies}
In this section we apply the comparison lemmas to the question of when
a sequence of hierarchies converges to a limiting hierarchy. Let us
first discuss what we mean by convergence.

Fix a point $x_0\in \CC(S)$ and let $B_R=B_R(x_0)$ denote the
$R$-neighborhood of $x_0$ in $\CC(S)$. For a tight geodesic $h$ with
$D(h)\subseteq S$, let $h\intersect B_R$ denote the 
following:
\begin{enumerate}
\item If $D(h)$ is a component domain of $(S,v)$ for some $v$ then, 
if $v\subset B_R$ we let $h\intersect B_R  = h$, and otherwise 
$h\intersect B_R =\emptyset$.

\item If $D(h)=S$ then $h\intersect B_R$ is the set of all positions
  of $h$ that lie in $B_R$.
\end{enumerate}
(In (2) this includes $\I(h)$ and/or $\T(h)$ if their bases are in
$B_R$.)

For a hierarchy $H$, let $H\intersect B_R = \{h\intersect B_R: h\in H\}$.
We say that a sequence $\{H_n\}$ of hierarchies {\em converges to a
  hierarchy $H$} if for all $R>0$, $H_n\intersect B_R = H\intersect
B_R$ for large enough $n$. Clearly if $H$ is a finite hierarchy this
just means that eventually $H_n=H$.

It is also easy to see the following: Suppose that for all $R>0$, the
sets $H_n\intersect B_R$ are eventually constant. Then $\{H_n\}$ converges
to a unique hierarchy $H$. We can now prove the following result:

\begin{theorem+}{Convergence of Hierarchies}
Let $\{H_n\}_{n=1}^\infty$ be a sequence of hierarchies such that either
\begin{enumerate}
\item For a fixed $K$ and any $n,m$, $H_n$ and $H_m$ are $K$-separated
  at the ends, or
\item There exists $K>0$ and a vertex $x_n$ on each $g_{H_n}$
  such that, for each
  $R'>0$, there exists $n=n_{R'}$ so that for all $m\ge n$,
  $g_{H_n}$ and $g_{H_m}$ are $(K,R')$-parallel at $x_n$ and $x_m$.
\end{enumerate}
Then $\{H_n\}$ has a convergent subsequence.
\end{theorem+}

\begin{proof}
  Fix an arbitrary $x_0$ and let $\UU_R$ denote the set of all
  vertices in $H_n\intersect B_R(x_0)$ for all $n>0$. We claim that
  $\UU_R$ is finite for each $R>0$. The theorem follows immediately
  since this implies that $H_n\intersect B_R$ varies in a finite set
  of possibilities for each $R$, and so the usual diagonalization step
  extracts a subsequence $H_{n_k}$ for which $H_{n_k}\intersect B_R$ is
  eventually constant. 

To show $\UU_R$ is finite, consider first case (1). 
Fix a slice $\tau_1$ in $H_1$, and note that
any vertex $v$ in $H_n$ appears in some complete slice $\tau$ of $H_n$. 
Consider a hierarchy $J(v)$
joining clean markings compatible with $\tau$ and $\tau_1$ respectively.
Lemma \lref{Slice Comparison}, case (1), implies that
$J(v)$ is $(K',M)$-pseudo-parallel to
$H_1$, so each geodesic in $J(v)$ has
length bounded either by $M$ or by a constant plus the length of a
geodesic in $H_1$. Since $H_1$ is finite this gives some uniform
bound, so every marking compatible with a slice of $H_n$ can be
transformed to a marking compatible with a slice of $H_1$ in a bounded
number of elementary 
moves.  This means the set of all base curves that occur in such
markings is finite, and this bounds the set of all vertices occurring
in non-annular geodesics in all $H_n$. The annular geodesics are
determined by their initial and terminal markings, up to a finite
number of choices (by the definition of tightness), and hence those
vertices are finite in number as well. (Note we have actually proved
finiteness for all the vertices in all $H_n$, without mention of $R$).

In case (2), the condition implies that there is some bound
$d(x_0,x_n)\le R_0$ for all $n>0$. 
Given $R$, choose $R'=R_0 + R+3K + 8\delta + 5$ and
let $n=n_{R'}$. For $m\ge n$ let $\ell_m$ be the segment of radius $R'$
around $x_m$. 
The $(K,R')$-parallel condition means that
$d(x_n,x_m)\le K$ and
either $\ell_m$ is in a $K$-neighborhood of $g_{H_n}$ or
$\ell_n$ is in a $K$-neighborhood of $g_{H_m}$. 
In either case, the
triangle inequality and $\delta$-hyperbolicity imply that, if $x\in
g_{H_m}$ and $d(x,x_m)\le R_0+R+1$ then $x$ is the center of a segment
in $g_{H_m}$ of radius $6\delta+4$ contained in a
$2\delta$-neighborhood of $\ell_n$ (we should assume that $K>\delta$,
which entails no loss of generality).

Now any vertex of $H_m\intersect B_R$ occurs in some complete slice
$\tau$ of $H_m$ with bottom vertex $x$ in $B_{R+1}$ (the slice will be
complete because $\I(H_m)$ and $\T(H_m)$ are sufficiently far away
from $B_R$
that they have non-trivial restriction to any domain occurring in
$H_m\intersect B_R$ -- so one can apply Lemma \ref{Subordinate
  Intersection 3}). Thus
$d(x,x_m) \le R_0 + R + 1$ by the triangle inequality, and the
previous paragraph implies that, for suitable
$x'\in \ell_n$, $g_{H_m}$ and $g_{H_n}$ are
$(2\delta,6\delta+4)$-parallel at $x$ and $x'$. By case (2) of
Lemma  \ref{Slice Comparison} we can again conclude that a
marking compatible with 
$\tau$ can be connected to some marking compatible with a slice in
$H_n\intersect B_{R'}$ by a 
sequence of elementary moves whose length is bounded only in terms of
$H_n\intersect B_{R'}$. The argument then proceeds as in case (1).
\end{proof}

We have the following immediate consequence of this argument:

\begin{corollary+}{Finite Geodesics}
Given a pair of points $x,y\in \CC_0(S)$ there are only a finite number of
tight geodesics joining them.
\end{corollary+}

\begin{proof}
Fix markings $\I$ and $\T$ containing
$x$ and $y$, respectively.  Each tight geodesic connecting $x$ to $y$
can be extended to a hierarchy connecting $\I$ and $\T$, and the
finiteness argument in case (1) of Theorem \ref{Convergence of
  Hierarchies} implies this set of hierarchies is finite.
\end{proof}

%% file: conjugacy.tex
\section{Conjugacy bounds in the Mapping Class Group}
\label{MCG}
In this section we will apply the hierarchy construction to the
Mapping Class Group $\Mod(S)$. Our main goals will be Theorem
\ref{Quasigeodesic Words}, stating that hierarchies give rise to
quasi-geodesic words in $\Mod(S)$, and Theorem \ref{Conjugacy Bound},
which gives a linear upper bound for the length of the minimal word
conjugating two pseudo-Anosov mapping classes.

We recall first that any generating set for $\Mod(S)$ induces a {\em
word metric} on the group, denoting by $|g|$ the length of the
shortest word in the 
generators representing $g\in\Mod(S)$, and by $|g^{-1}h|$ the distance
between $g$ and $h$. The metrics associated to any two
(finite) generating sets are bilipschitz equivalent.

\subsection{Paths in the mapping class group}
\label{paths in MCG}
Our first step is to show how a resolution of a hierarchy into a
sequence of slices gives rise to a word in the mapping class
group. This is a completely standard procedure involving the
connection between groupoids and groups. Theorem
\lref{Efficiency of Hierarchies} will imply that these words are in fact
quasi-geodesics.

Let $\til \MM$ be the graph of complete clean markings of $S$
connected by elementary moves, as in \S\ref{efficiency section}.
The action of $\Mod(S)$ on
$\til \MM$ is not free  -- a mapping class can permute the components
of a marking or reverse their orientations -- but it has finite
stabilizers. The quotient $\MM=\til \MM/\Mod(S)$ is a finite graph, and we
let $D$ denote its diameter. Fix a marking $\mu_0\in\til \MM$
and let $\Delta_{j}\subset\Mod(S)$ denote the set of elements $g$ such
that $d_{\til \MM}(\mu_0,g(\mu_0)) \le j$. 
Any marking $\mu\in \til \MM$ is at most distance $D$ from some
$\psi(\mu_0)$, with $\psi\in\Mod(S)$ determined up to pre-composition
by elements of $\Delta_{2D}$.

Note that $\Delta_{j}$ is a finite set for any $j$. Now given any
$\psi\in \Mod(S)$, 
we can write it as a word in $\Delta_{2D+1}$ as follows: connect
$\mu_0$ to $\psi(\mu_0)$ by a shortest path
$\mu_0,\mu_1,\ldots,\mu_N=\psi(\mu_0)$ in $\til\MM$.
For
each $i<N$ choose $\psi_i$ such that $\psi_i(\mu_0)$ is within $D$ moves
of $\mu_i$, and let $\psi_N = \psi$. Then $\psi_i$ and $\psi_{i+1}$
differ by pre-composition with 
an element $\delta_{i+1}$ in $\Delta_{2D+1}$. Thus we can write 
$\psi = \delta_1\cdots\delta_{N}$.  

This gives an upper bound on the word length of $\psi$, proportional
to $N=d_{\til\MM}(\mu_0,\psi(\mu_0))$. Of course, the word we obtain
here can be translated to a word in any other finite generating set in
the standard way, and its length will only increase by a bounded
multiple (depending on the generating set).

In the other direction,
suppose that $\psi$ can be written
as $\alpha_1\cdots\alpha_{M}$ with $\alpha_i$ in some fixed finite
generating set. Then the sequence of markings
$\{\mu_j=\alpha_1\cdots\alpha_j(\mu_0)\}$ satisfies the property that $\mu_j$
and $\mu_{j+1}$ are separated by some bounded number of elementary
moves, with the bound depending on the generating set. This bounds 
$d_{\til\MM}(\mu_0,\psi(\mu_0))$ above linearly in terms of $|\psi|$.

Now let $H$ be any
hierarchy with $\I(H)=\mu_0$ and $\T(H) = \psi(\mu_0)$. 
Theorem \lref{Efficiency of Hierarchies} gives upper and lower bounds
on $d_{\til  \MM}(\mu_0,\psi(\mu_0))$ in terms of $|H|$, and we 
immediately obtain:

\begin{theorem+}{Quasigeodesic Words}
Fix a complete clean marking  $\mu_0$ of $S$ and a set of generators for
$\Mod(S)$, and let $|\cdot|$ be the word metric with respect to these
generators. Then there are $c_2,c_3>0$ such that the following holds:

Given any $\psi\in\Mod(S)$ let $H$ be a hierarchy such that $\I(H)=\mu_0$ 
and $\T(H) = \psi(\mu_0)$. 
Then the words in $\Mod(S)$ constructed from resolutions of $H$
via the procedure in this section are quasi-geodesics, and in particular
$$
c_2^{-1} |H| - c_3 \le |\psi| \le c_2|H|.
$$
\end{theorem+}
(We remark that the additive constant $c_3$ can be removed if we
always choose the hierarchy to have length 0 when $\psi={\mathrm
  id}$).

Note also that an estimate on the word length $|\psi|$ can be obtained
purely in terms of the quantities $d_Y(\mu_0,\psi(\mu_0))$, using
Theorem \ref{Move distance and projections}.


\subsection{The conjugacy bound}
\label{conjugacy}
We are now ready to prove the main theorem of this section:

\begin{theorem+}{Conjugacy Bound}  
Fixing a set of generators for $\Mod(S)$,
there exists  a constant $K$ 
such that if $h_1,h_2\in \Mod(S)$ are conjugate pseudo-Anosovs there 
is a
conjugating element $w$ with $|w|\leq K(|h_1|+|h_2|)$.
\end{theorem+}

Let $\delta$ be the hyperbolicity constant for $\CC(S)$.
Say that two geodesics are $c$-fellow travelers
if each is in a $c$-neighborhood of the other (Hausdorff distance $c$)
and their endpoints (if any) can be paired to be within distance $c$
of each other. The following three lemmas are standard for any  
hyperbolic metric space. 
 
\begin{lemma}{Bi-infinite fellow traveling}  For any $K$, if 
$\beta_1,\beta_2$ are two $K$-fellow traveling bi-infinite geodesics,  
then they are actually $2\delta$ fellow travelers. 
\end{lemma}

\begin{pf}  
Let  $\bar\beta_1$ be any segment of $\beta_1$.  Choose  
points $y_1,y_2\in\beta_1\setminus \bar\beta_1$ on either side of  
$\bar \beta_1$ whose distance to $\bar\beta$ is $2\delta+K+1$ and
points $x_i$ on $\beta_2$ such that $d(x_i,y_i)=K$.  The
quadrilateral $[x_1x_2y_2y_1]$ is $2\delta$-thin by hyperbolicity:
that is, each edge is within $2\delta$ of the union of the other three.
By the triangle inequality no point of $\bar\beta_1$ can be within  
$2\delta$
of $[x_i,y_i]$ so each point must be  within  $2\delta$ of   
$\beta_2$.  
Since $\bar\beta_1$ was arbitrary we are done.
\end{pf}

For the next two lemmas, let $\beta$ denote any bi-infinite geodesic,
and let $\pi=\pi_\beta:\CC(S)\to \beta$ be any map which for each $x$
picks out some closest point on $\beta$. (Note that $\pi$ need not be
uniquely defined.) Let 
$[a,b]$ denote any geodesic joining $a$ and $b$, taken to be a segment
of $\beta$ whenever $a$ and $b$ lie on $\beta$.

\begin{lemma}{projection stable}
Let $x\in\CC(S)$ and $z\in \beta$, such that 
$d(x,z) \le d(x,\pi_\beta(x)) + k$ for some $k\ge 0$.
Then $d(\pi(x),z) \le k + 4\delta$.
\end{lemma}

\begin{pf}
We may assume $d(\pi(x),z) > 2\delta$. Let $m\in[\pi(x),z]$ be
distance  $2\delta+\ep$ from $\pi(x)$, for $\ep>0$.
By $\delta$-hyperbolicity,
$m$ is distance at most $\delta$ from either some $m_1\in[x,\pi(x)]$
or some $m_2 \in [x,z]$. The former case cannot occur, since then
$d(m_1,\pi(x)) \le
\delta$ and hence $d(\pi(x),m) \le 2\delta$, a contradiction. 
Thus we have
$d(x,\pi(x)) \le d(x,m_2) + \delta$, and
$d(m,z) \le d(m_2,z) + \delta$.
Adding these together and using the hypothesis, we conclude
that $d(m,z) \le 2\delta + k$ and hence $d(\pi(x),z) \le 4\delta + k +
\ep$. Sending $\ep\to 0$ gives the desired result.
\end{pf}

\begin{lemma}{thin quadrilateral}
Let $x,y$ be any two points in $\CC(S)$, such that
$d(\pi_\beta (x),\pi_\beta( y)) > 8\delta + 2$.
Let $\sigma$ be the subsegment of
$[\pi_\beta(x),\pi_\beta(y)]$ on $\beta$ whose endpoints
are distance $4\delta+1$ from $\pi_\beta(x)$ and $\pi_\beta(y)$ respectively.
Then $\sigma$ is in a $2\delta$ neighborhood of $[x,y]$.
\end{lemma}

\begin{pf}
Form the quadrilateral whose sides are $[x,\pi(x)], [\pi(x),\pi(y)],
[y,\pi(y)]$ and $[x,y]$. By $\delta$-hyperbolicity, any point $z\in
[\pi(x),\pi(y)]$ is at most $2\delta$ from one of the other three
sides. Suppose that this is the side $[x,\pi(x)]$. Then there is some
$m\in[x,\pi(x)]$ such that $d(m,z)\le 2\delta$, and since $d(x,\pi(x))
\le d(x,z)$, we must have $d(m,\pi(x)) \le 2\delta$. It follows that
$d(\pi(x),z)\le 4\delta$. The same argument applies to $[y,\pi(y)]$,
and it follows that if $z\in \sigma$ then it must be distance
$2\delta$ from $[x,y]$.
\end{pf}

\begin{proposition+}{Axis} Let $h$ be a pseudo-Anosov element in  
$Mod(S)$. There
exists a bi-infinite tight geodesic
$\beta$ such that for each $j$, $h^j(\beta)$ and $\beta$ are
$2\delta$ fellow travelers.  Moreover there exists a hierarchy
$H$ with main geodesic $\beta$.  
\end{proposition+}
\begin{proof}
Pick any $x\in \CC(S)$.   Let $\beta_n$ be a tight geodesic
joining $h^{-n}(x)$ and
$h^n(x)$ (extend the endpoints in an arbitrary way to complete
markings $\I(\beta_n)$ and $\T(\beta_n)$),
and let $H_n$ be a hierarchy with main geodesic $\beta_n$.
By Proposition 3.6 of \cite{masur-minsky:complex1}, the sequence
$\{h^k(x),k\in \Z\}$ satisfies 
$d(h^k(x),x)\geq c|k|$ for some $c>0$ (independent of $x$ or $h$) so the
sequence is a $\frac{d_0}{c}$-
quasi-geodesic, where $d_0=d(h(x),x))$.  
By $\delta$-hyperbolicity  there is a 
constant $c'=c'(c,d_0,\delta)$ so that $\beta_n$ and the
sequence $\{h^j(x)\}_{|j|\le n}$ lie in a $c'$-neighborhood of each other.
This implies that there exist $x_n\in\beta_n$, so 
that given any $R$,  for $n,m$ sufficiently large,  $\beta_n$ and $\beta_m$
are $(2c',R)$-parallel at $x_n,x_m$.

We apply Theorem \lref{Convergence of Hierarchies} to find a geodesic
$\beta$ and a  
hierarchy $H$, which is the limit of a subsequence of $H_n$.  For  
each $j$, $h^j(\beta)$ is a $2c'$-fellow traveler to $\beta$.   
Applying Lemma \ref{Bi-infinite fellow traveling}
gives the result. 
\end{proof}

We call $\beta$ a {\em quasi-axis} for $h$.  
We will need to know the following:
\begin{lemma}{definite translation}
Given $A>0$, there is an integer $N>0$,
independent of $h$, such that for any $x\in\CC(S)$ and $n\ge N$, 
$$
  d(\pi(x),\pi(h^n(x))) \ge A.
$$
\end{lemma}
\begin{pf}
We first observe that, 
if $g$ is any power of $h$, and $\beta$ a quasi-axis for $h$, then
\begin{equation}\label{quasicommute}
d(\pi g(x),g\pi(x)) \le 10\delta
\end{equation}
for any $x\in\CC$. The proof will be given below. 

Now using the inequality $d(\pi(x),h^n\pi(x)) \ge
c|n|$, with $c$ independent of $x$ and $h$, from Proposition 3.6 of
\cite{masur-minsky:complex1}, 
we simply choose $N$ so that $cN > A + 10\delta$.

It remains to prove (\ref{quasicommute}):
Since $\beta$ and $g(\beta)$ are $2\delta$-fellow travelers, we have
$d(g(x),g(\beta)) \le d(g(x),\beta) + 2\delta$, or equivalently, since $g$
is an isometry,
$d(g(x),g\pi( x)) \le d(g(x),\pi g (x)) + 2\delta$. Now 
$\pi g \pi (x)$ is on $\beta$, and again by the fellow traveler
property, we have $d(g\pi( x),\pi g \pi (x)) \le 2\delta$. Thus
$d(g(x),\pi g \pi (x)) \le d(g(x),\pi g( x)) + 4\delta$. Applying Lemma
\ref{projection stable} with $k=4\delta$ and $z=\pi g \pi (x)$, we find
that
$d(\pi g (x), \pi g \pi (x)) \le 8 \delta$. 
We conclude that $d(\pi g( x), g\pi( x)) \le 10\delta$, as desired.
\end{pf}

\begin{pf*}{Proof of Conjugacy Theorem}
Suppose that $h_2 = w^{-1}h_1w$.
Lemma \ref{definite translation} guarantees that we can choose $N$
independent of $h_1$ and $h_2$ such that
so that $d(\pi(x),\pi h_i^n(x))) \ge 40\delta + 24$ for all $x\in\CC$ and
$n\ge N$. Let $g_i=h_i^N$ (for $i=1,2$). In the proof below, let
$C_1,C_2,\ldots$ denote positive constants which are independent of
$h_1$ and $h_2$. 

Fix  a complete clean marking $\mu_0$ in $S$.
Let
$H_i$ be a hierarchy such that $\I(H_i) = \mu_0$ and $\T(H_i) = g_i(\mu_0)$.
We may also assume
the main geodesic of $H_i$ is a segment
$[v,g_i(v)]$ for a base curve $v$ of $\mu_0$.
By Theorem \lref{Quasigeodesic Words}, we have
$|H_i|\le C_1|g_i| \le NC_1 |h_i|$.

Since $w$ acts by natural isomorphisms on $\CC(S)$ and and all the subsurface
complexes, we have a hierarchy $w(H_2)$ with
main geodesic $[w(v),w g_2(v)]=[w(v),g_1 w(v)]$, 
and $|w(H_2)|=|H_2|$. 

Form a quasi-axis $\beta$ for $g_1$, together with a hierarchy $H$,
and form the segments $$I_0=[\pi(v),\pi g_1(v)]$$
and 
$$
I'_m=[\pi g_1^m  w(v),\pi g_1^{m+1} w(v)]  
$$
on $\beta$.
Each of these have length at least $40\delta + 24$.
Let 
$\sigma_0\subset I_0$ and $\sigma'_m\subset I'_m$ be the subsegments
obtained by removing $4\delta+1$-neighborhoods of the endpoints. 

The $\{I'_m\}_{m\in\Z}$ tile $\beta$ and therefore the gaps between
$\sigma'_m$ and $\sigma'_{m+1}$ have length $8\delta+2$. It follows
that there exists some $m\in\Z$ such that $\sigma_0$ and $\sigma'_m$
overlap on a segment $\zeta$ of length at least $12\delta + 10$.

Let $w' = g_1^m w$, and note that $w'$  also conjugates $h_1$ and
$h_2$. We will now 
bound the word length  of $w'$, by bounding
$d_{\til \MM}(\mu_0,w'(\mu_0))$.

By Lemma \ref{thin quadrilateral},
the segment $\zeta$ is in a $2\delta$-neighborhood of both
$[v,g_1(v)]$ and $[w'(v),g_1 w'(v)]$.
Let $x$ be a vertex of $\zeta$ nearest its midpoint, so that
$\zeta$ contains an interval $L$ of radius $6\delta+4$ around $x$,
and let $u,u'$ be vertices on 
$[v,g_1(v)]$ and $[w'(v),g_1 w'(v)]$, 
respectively, which are nearest  to $x$.
Thus the main geodesics of $H$ and $H_1$ are
$(2\delta,6\delta+3)$-parallel at $x$ and 
$u$, and similarly for those of $H$ and $w'(H_2)$ at $x$ and $u'$. 
This will allow us to apply Lemma \lref{Slice Comparison} below.

Resolve $H_1$ into a sequence of slices. One of them must have bottom
vertex $u$ (see proof of Proposition \ref{Resolution into slices})
 -- let $\mu_1$ be a clean marking compatible with this slice. 
The resolution gives  a bound $d_{\til\MM}(\mu_0,\mu_1) \le C_2|h_1|$,
by Proposition \ref{Resolution into slices} and Lemma \ref{marking
  move bound}.
Similarly, resolve $w'(H_2)$, find a slice with bottom vertex
$u'$, let $\mu_2$ be a clean marking compatible with this slice, and
conclude $d_{\til\MM}(\mu_2,w'(\mu_0)) \le C_3|h_2|$. 

Let $\mu_3$ be a clean marking associated to a slice of $H$ with
base vertex $x$. Let $W$ be a 
hierarchy with $\I(W)=\mu_1$ and $\T(W)=\mu_3$. 
Case (2) of Lemma  \lref{Slice Comparison} tells us that $W$ is
$(K',M)$-pseudo-parallel to $H_1$ (with $K',M$ depending only on $\delta$).
In particular this means
$|W|\le C_4|H_1|$. Resolving $W$, we obtain 
$d_{\til\MM}(\mu_1,\mu_3) \le C_5|h_1|$. 

Similarly, join $\mu_3$ to
$\mu_2$ by a hierarchy $W'$.
The same argument as for $W$ gives us 
$d_{\til\MM}(\mu_2,\mu_3) \le C_5|h_2|$.

Adding these bounds, we obtain 
$d_{\til\MM}(\mu_0,w'(\mu_0)) \le C_6(|h_1|+|h_2|)$, which as in
\S\ref{paths in MCG} gives the desired bound on $|w'|$.
\end{pf*}

%% file: biblio.tex
\providecommand{\bysame}{\leavevmode\hbox to3em{\hrulefill}\thinspace}